\documentclass[12pt,oneside,reqno]{amsart}
\usepackage{graphicx}
\usepackage{mathrsfs}
\usepackage{stmaryrd}
\usepackage{amsfonts}
\usepackage{cite}
\usepackage{enumerate,amsmath,amssymb,amsthm}
\usepackage{booktabs}
\usepackage{diagbox}
\usepackage{amsfonts,color}

\newcommand{\dif}{\mathrm{d}}

\newcommand{\be}{\begin{eqnarray}}
\newcommand{\ee}{\end{eqnarray}}
\newcommand{\ce}{\begin{eqnarray*}}
\newcommand{\de}{\end{eqnarray*}}
\newtheorem{theorem}{Theorem}[section]
\newtheorem{lemma}[theorem]{Lemma}
\newtheorem{remark}[theorem]{Remark}
\newtheorem{definition}[theorem]{Definition}
\newtheorem{proposition}[theorem]{Proposition}
\newtheorem{Examples}[theorem]{Examples}
\newtheorem{corollary}[theorem]{Corollary}
\def\e{\varepsilon}

\def\a{\alpha}

\def\b{\beta}
\def\d{\delta}
\def\D{\Delta}
\def\p{\partial}
\def\g{\gamma}
\def\s{\sigma}
\def\eps{\varepsilon}

\def\la{\langle}
\def\ra{\rangle}
\def\[{{\Big[}}
\def\]{{\Big]}}
\def\<{{\langle}}
\def\>{{\rangle}}
\def\({{\Big(}}
\def\){{\Big)}}

\def\no{\nonumber}
\def\bt{\begin{theorem}}
\def\et{\end{theorem}}
\def\bl{\begin{lemma}}
\def\el{\end{lemma}}
\def\br{\begin{remark}}
\def\er{\end{remark}}
\def\bx{\begin{Examples}}
\def\ex{\end{Examples}}
\def\bd{\begin{definition}}
\def\ed{\end{definition}}
\def\bp{\begin{proposition}}
\def\ep{\end{proposition}}
\def\bc{\begin{corollary}}
\def\ec{\end{corollary}}

\def\cI{{\mathcal I}}
\def\cJ{{\mathcal J}}

\def\cL{{\mathcal L}}

\def\cP{{\mathcal P}}

\def\mE{{\mathbb E}}

\def\mN{{\mathbb N}}

\def\mP{{\mathbb P}}

\def\mR{{\mathbb R}}

\def\mW{{\mathbb W}}

\def\sB{{\mathscr B}}
\def\sC{{\mathscr C}}

\def\sF{{\mathscr F}}

\def\sL{{\mathscr L}}

\def\sR{{\mathscr R}}

\def\tx{\tilde{x}}
\def\cx{\check{x}}
\def\geq{\geqslant}
\def\leq{\leqslant}

\begin{document}

\allowdisplaybreaks
\title{Averaging principles for nonautonomous multiscale McKean-Vlasov stochastic systems}

\author{Jie Xiang and Huijie Qiao}

\thanks{{\it AMS Subject Classification(2020):} 60H10}

\thanks{{\it Keywords:} Nonautonomous multiscale McKean-Vlasov stochastic systems, averaging principles, the nonautonomous Poisson equation, evolution system of measures}

\thanks{This work was supported by NSF of China (No.12071071) and the Jiangsu Provincial Scientific Research Center of Applied Mathematics (No. BK20233002).}

\thanks{Corresponding author: Huijie Qiao, hjqiaogean@seu.edu.cn}

\subjclass{}

\date{}

\dedicatory{School of Mathematics,
Southeast University,\\
Nanjing, Jiangsu 211189, P.R.China\\
}

\begin{abstract}
This paper investigates a class of nonautonomous multiscale McKean-Vlasov stochastic systems. By leveraging the nonautonomous Poisson equation, we rigorously establish both strong and weak averaging principles, accompanied by explicit convergence rates. Notably, the coefficients of the averaging equations derived in the general case retain dependence on the scaling parameter $\e$. However, under the additional assumptions that the fast-scale coefficients are either asymptotically convergent or time-periodic, we demonstrate that the slow component converges, in the strong or weak sense, to averaging equations with coefficients independent of $\eps$.
\end{abstract}

\maketitle \rm

\section{Introduction}\label{intro}

McKean-Vlasov stochastic differential equations (SDEs), also termed distribution-\\dependent or mean-field SDEs, characterize the limiting evolution of a typical particle in mean-field interacting particle systems as the number of particles tends to infinity. A defining feature of these equations is that their drift and diffusion coefficients depend not only on the current state of the solution process but also explicitly on its probability law, that is, the marginal distribution of the process at each time. The theoretical framework for McKean-Vlasov SDEs originates in the seminal work of McKean \cite{hpm}, which was directly motivated by Kac's program in kinetic theory. Since then, these equations have attracted sustained attention across probability theory, analysis, and mathematical physics, and have been investigated rigorously from multiple perspectives. Notably, substantial advances have been achieved concerning well-posedness and regularity of solutions, stability and asymptotic behavior, rigorous correspondence with nonlinear Fokker-Planck equations, and ergodic properties-including existence, uniqueness, and convergence rates of invariant measures. For comprehensive surveys of these developments, we refer the reader to \cite{br, br1, blpr, carm, dq1, dq2, hw, rw, rz, Wang} and the references therein.

The averaging principle, originally introduced by Khasminskii in his seminal work \cite{rk}, constitutes a foundational analytical framework for SDEs and has become indispensable across numerous applied disciplines. By rigorously approximating the slow dynamics of multiscale stochastic systems via reduced-order averaging equations, whose coefficients are obtained by averaging over the fast-scale variables, the principle enables substantial dimensionality reduction while preserving essential asymptotic behavior. This approximation not only simplifies computational and theoretical analysis but also yields tractable models that retain the core qualitative and quantitative features of the original system. Over the past several decades, the averaging principle has spurred extensive theoretical development and broad application (see, e.g., \cite{hll, lx, lrsx, pix, q, qw, qw2, rsx, sxw, syw, xq1, xq2}).

Nonautonomous multiscale systems arise naturally in diverse scientific domains, including climate-weather modeling, where interactions across diurnal and seasonal time scales play a central role (\cite{bm, csg}), and stochastic neural dynamics driven by time-varying external inputs (\cite{gw}). While averaging principles have been established for several classes of nonautonomous multiscale stochastic systems, such as nonautonomous multiscale SDEs (\cite{hwx, swx0, swx, Uda, gW}) and nonautonomous multiscale SPDEs (\cite{cl, csx, lswx, xw}), these results presuppose coefficients that depend solely on state variables and time, with no distributional dependence. In contrast, nonautonomous multiscale McKean-Vlasov stochastic systems incorporate explicit dependence of coefficients on the law of the solution, a structural feature essential for modeling collective, self-reinforcing phenomena in large interacting populations under time-varying environments. Such systems thus provide a more faithful and mathematically refined framework for complex real-world processes. Motivated by this gap in the literature and the need for rigorous asymptotic analysis, this paper establishes both strong and weak averaging principles with explicit convergence rates for a general class of nonautonomous multiscale McKean-Vlasov stochastic systems.

Concretely speaking, consider the following slow-fast system on $\mathbb{R}^n\times\mathbb{R}^m$:
\be\left\{\begin{array}{l}
\dif X_t^{\eps}=b_1(X_t^{\eps}, \sL_{X_t^{\eps}}, Y_t^{\eps}) \dif t+\s_1(X_t^{\eps}, \sL_{X_t^{\eps}}, Y_t^{\eps}) \dif B_t, \quad X_0^{\eps}=\xi, \\
\dif Y_t^{\eps}=\eps^{-1} b_2(t/\eps,X_t^{\eps}, \sL_{X_t^{\eps}}, Y_t^{\eps})\dif t+\eps^{-1/2} \s_2(t/\eps,X_t^{\eps}, \sL_{X_t^{\eps}}, Y_t^{\eps})\dif W_t, \quad Y_0^{\eps}=\varrho,
\end{array}\right.
\label{orieq}
\ee
where $\left(B_t\right),\left(W_t\right)$ are $d_1$-dimensional and $d_2$-dimensional standard Brownian motions, respectively, defined on a complete filtered probability space $(\Omega,\sF,\{\sF_t\}_{t \in[0, T]}, \mathbb{P})$ and are mutually independent. Moreover, these mappings
$b_1: \mR^n \times \cP_2(\mR^n) \times \mR^m \rightarrow \mR^n$,
$\s_1: \mR^n \times \cP_2(\mR^n) \times \mR^m \rightarrow \mR^{n \times d_1}$,
$b_2: \mR_+\times\mR^n \times \cP_2(\mR^n) \times \mR^m\rightarrow \mR^m$,
$\s_2: \mR_+\times\mR^n \times \cP_2(\mR^n) \times \mR^m\rightarrow \mR^{m \times d_2}$
are all Borel measurable, and $\xi, \varrho$ are two random variables. Due to the time-inhomogeneity of the fast part coefficients in system (\ref{orieq}), constructing a well-defined averaging equation and establishing explicit convergence rates present the primary technical challenge: the associated frozen equation (\ref{frozeq}) generally lacks an invariant measure, precluding standard ergodic averaging arguments. To overcome this obstacle, we employ an evolution system of measures together with a suitably formulated nonautonomous Poisson equation. Notably, the coefficients of the averaging equations obtained under general conditions retain explicit dependence on the scaling parameter $\e$. Consequently, to derive strong and weak averaging principles wherein the averaging coefficients are independent of $\e$, we impose additional structural assumptions on the fast-scale dynamics: specifically, that the fast coefficients either converge asymptotically or are time-periodic. In establishing the three weak averaging principles, a second key difficulty arises-the estimation of solutions to three auxiliary Cauchy problems whose formulations inherently involve the law of the solution. Drawing inspiration from the decoupling methodology introduced in \cite{blpr}, we adapt this technique to handle distributional dependencies.

It is worthwhile to mention that in \cite{slg}, Shi, Liu and Gao proved strong and weak averaging principles for nonautonomous multiscale McKean-Vlasov SDEs with almost periodic fast coefficients by means of tightness arguments. Comparing Theorem \ref{xbarstrper} and \ref{xbarweakper} with Theorem 6.1 and 6.3 in \cite{slg}, we find that explicit convergence rates are presented here.

The remainder of this paper is organized as follows. In Section \ref{noteassu}, we introduce the notations used throughout the paper. Section \ref{main} presents the main results. Section \ref{frozpoi} is devoted to the study of the time-inhomogeneous frozen dynamics and the associated nonautonomous Poisson equation. The proofs of the main results are provided in Sections \ref{geneproo}-\ref{periconvproo}. Finally, in Section \ref{example}, we provide two examples to explain our results.

The following convention will be used throughout the paper: $C$ with or without indices will denote different positive constants whose values may change from one place to another.

\section{Notations}\label{noteassu}

In this section, we will recall some notations.

\subsection{Notations}\label{nn}
In this subsection, we introduce some notations used in the sequel.

We use $|\cdot|$ and $\|\cdot\|$ to denote the norms of vectors and matrices, respectively. Let $\<\cdot,\cdot\>$ represent the scalar product in $\mR^d$, and let $A^*$ denote the transpose of the matrix $A$.

Let $\sB(\mR^n)$ be the Borel $\s$-algebra on $\mR^n$, and $\cP(\mR^n)$ represent the space of all probability measures defined on $\sB(\mR^n)$ equipped with the usual topology of weak convergence. Let $\cP_2(\mR^n)$ be the collection of the probability measures $\mu$ on $\sB(\mR^n)$ satisfying
$$
\mu(|\cdot|^2):=\int_{\mR^n}|x|^2\mu(\dif x)<\infty.
$$
This space is a Polish space under the $L^2$-Wasserstein distance, defined as
$$
\mW_2(\mu_1,\mu_2):=\inf_{\pi\in\sC(\mu_1,\mu_2)}\bigg(\int_{\mR^n\times\mR^n}|x-y|^2\pi(\dif x,\dif y)\bigg)^{\frac{1}{2}},\mu_1,\mu_2 \in\cP_2(\mR^n),
$$
where $\sC(\mu_1,\mu_2)$ is the set of all couplings $\pi$ with marginals distributions $\mu_1$ and $\mu_2.$ Moreover, for any $x\in\mR^n$, the Dirac measure $\delta_x$ belongs to $\cP_2(\mR^n)$, and if $\mu_1=\sL_X$ and $\mu_2=\sL_Y$ are the corresponding distributions of random variables $X$ and $Y$, respectively, then
$$\mW_{2}(\mu_1,\mu_2)\leq(\mathbb{E}|X-Y|^{2})^{\frac12},$$
where $\mathbb{E}$ denotes the expectation with respect to $\mP$.

\subsection{Derivatives for functions on $\cP_2(\mR^d)$}\label{deri}
\quad In this subsection, we recall the definition of the $L$-derivative for functions on $\cP_2(\mR^d)$. This definition was first introduced by Lions in \cite{card} (see also \cite{blpr,carm}), who employed abstract probability spaces to describe the $L$-derivatives. For the sake of clarity, we present a more straightforward formulation here (cf. \cite{rw}). Let $I$ be the identity map on $\mR^n$. Set $\mu(\phi):=\int_{\mR^n}\phi(x)\mu(\dif x)$ for $\mu\in\cP_2(\mR^n)$ and $\phi\in L^2(\mR^n,\sB(\mR^n),\mu;\mR^n)$, where $L^2(\mR^n,\sB(\mR^n),\mu;\mR^n)$ stands for the space of Borel measurable functions $\phi: \mR^n\rightarrow\mR^n$ with $\int_{\mR^n}|\phi(x)|^2\mu(\dif x)<\infty$. Moreover, simple calculations show that $\mu\circ(I+\phi)^{-1}\in\cP_2(\mR^n)$.

\bd\label{Lderidef}
(i) A function $g:\cP_2(\mR^n)\rightarrow\mR$ is called $L$-differentiable at $\mu\in\cP_2(\mR^n)$ if the functional
$$
L^2(\mR^n,\sB(\mR^n),\mu;\mR^n)\ni\phi\rightarrow g(\mu\circ(I+\phi)^{-1})
$$
is Fr\'{e}chet differentiable at $0\in L^2(\mR^n,\sB(\mR^n),\mu;\mR^n)$, that is, there exists a unique $\xi\in L^2(\mR^n,\sB(\mR^n),\mu;\mR^n)$ such that
$$
\lim_{\mu(|\phi|^2)\rightarrow0}\frac{g(\mu\circ(I+\phi)^{-1})-g(\mu)-\mu(\<\xi,\phi\>)}{\sqrt{\mu(|\phi|^2)}}=0.
$$
In this case, we denote $\p_\mu g(\mu)=\xi$ and call it the $L$-derivative of $g$ at $\mu$.\\
(ii) A function $g:\cP_2(\mR^n)\rightarrow\mR$ is called $L$-differentiable on $\cP_2(\mR^n)$ if the $L$-derivative $\p_\mu g(\mu)$ exists for all $\mu\in\cP_2(\mR^n)$.
\ed

We say that a matrix-valued function $h(\mu):=(h_{ij}(\mu))$ is $L$-differentiable at $\mu\in\cP_2(\mR^n)$, if all its components are $L$-differentiable at $\mu$. Moreover, we say that $\p_\mu h(\mu)(\tx)$ is differentiable at $\tx\in \mR^n$, if all its components are differentiable at $\tx$. We denote
$$\p_\mu h(\mu):=(\p_\mu h_{ij}(\mu)), \quad \p_{\tx}\p_\mu h(\mu)(\tx):=(\p_{\tx}\p_\mu h_{ij}(\mu)(\tx))$$

Let $d, l, k\in \mN=\{0,1,2,\ldots\}$ with $k\leq l$. We introduce the following spaces of functions.
\begin{definition}\label{Ckdef}
For a map $h(\cdot):\mR^n\rightarrow\mR$, we say that $h\in C^k(\mR^n,\mR)$ if all partial derivatives of $h$ up to order $k$ exist and are continuous on $\mR^n$. We say that $h\in C_b^k(\mR^n,\mR)$ if $h\in C^k(\mR^n,\mR)$ and all partial derivatives of $h$ up to order $k$ are bounded on $\mR^n$. For a matrix-valued map $h(\cdot):\mR^n\rightarrow\mR^{n\times d_1}$, we say that $h\in C^k(\mR^n,\mR^{n\times d_1})$ (resp.~$C_b^k(\mR^n,\mR^{n\times d_1})$) if all its components belong to $C^k(\mR^n,\mR)$ (resp.~$C_b^k(\mR^n,\mR)$).
\end{definition}

\begin{definition}\label{Ckldef}
For a map $h(\cdot,\cdot):\mR^n\times\mR^m\rightarrow\mR$, we say that $h\in C^{k,l}(\mR^n\times\mR^m,\mR)$, if $h(x,\cdot)\in C^l(\mR^m,\mR)$ for any $x\in\mR^n$ and $h(\cdot,y)\in C^k(\mR^n,\mR)$ for any $y\in\mR^m$, and the derivatives
\ce
\p_x^i\p_y^jh(x,y),\quad 0\leq i\leq k \ \text{and}\ 0\leq j\leq l \ \text{with}\ 1\leq i+j\leq l,
\de
are continuous at any $(x,y)$. We say $h\in C_b^{k,l}(\mR^n\times\mR^m,\mR)$, if $h\in C^{k,l}(\mR^n\times\mR^m,\mR)$, and all the above derivatives are uniformly bounded w.r.t. $(x,y)$.
For a matrix-valued map $h(\cdot,\cdot):\mR^n\times\mR^m\rightarrow\mR^{n \times d_1}$, we say $h\in C^{k,l}(\mR^n\times\mR^m,\mR^{n \times d_1})$ $(resp.~ C_b^{k,l}(\mR^n\times\mR^m,\mR^{n \times d_1}))$ if all the components belong to $C^{k,l}(\mR^n\times\mR^m,\mR)$ $(resp.~ C_b^{k,l}(\mR^n\times\mR^m,\mR))$.
\end{definition}

\begin{definition}\label{C1kdef}
For a map $h(\cdot):\cP_2(\mR^n)\rightarrow\mR$, we say $h\in C^{(1,d)}(\cP_2(\mR^n),\mR)$, if this map is $L$-differentiable at any $\mu\in\cP_2(\mR^n)$, its $L$-derivative $\p_\mu h(\mu)(\cdot):\mR^n\rightarrow\mR^n$ belongs to $C^d(\mR^n,\mR^n)$ for every $\mu\in\cP_2(\mR^n)$. Moreover, the derivatives $\p_{\tx}^i\p_\mu h(\mu)(\tx)$, $0\leq i\leq d$, are continuous at any $(\mu,\tx)$.
We say $h\in C_b^{(1,d)}(\cP_2(\mR^n),\mR)$, if $h\in C^{(1,d)}(\cP_2(\mR^n),\mR)$, and $\sup\limits_{\mu \in \cP_2(\mR^n), \tx\in \mR^n }\|\p_{\tx}^i\p_\mu h(\mu)(\tx)\|<\infty$, $0\leq i\leq d$.
For a matrix-valued map $h(\cdot):\cP_2(\mR^n)\rightarrow\mR^{n \times d_1}$, we say $h\in C^{(1,d)}(\cP_2(\mR^n),\mR^{n \times d_1})$ $(resp.~ C_b^{(1,d)}(\cP_2(\mR^n),\mR^{n \times d_1}))$ if all the components belong
to $C^{(1,d)}(\cP_2(\mR^n),\mR)$ $(resp.~ C_b^{(1,d)}(\cP_2(\mR^n),\mR))$.
\end{definition}

\begin{definition}\label{C22def}
For a map $h(\cdot):\cP_2(\mR^n)\rightarrow\mR$, we say $h\in C^{(2,2)}(\cP_2(\mR^n),\mR)$, if this map is twice $L$-differentiable at any $\mu\in\cP_2(\mR^n)$, its second-order $L$-derivative $\p_\mu^2 h(\mu)(\cdot,\cdot):\mR^n\times\mR^n\rightarrow\mR^{n\times n}$ is twice differentiable with respect to spatial variables, and the derivatives
$$
\p_\mu^2 h(\mu)(\tx,\cx), \quad \p_{\tx}\p_\mu^2 h(\mu)(\tx,\cx),\quad \p_{\cx}\p_\mu^2 h(\mu)(\tx,\cx),\quad \p_{\cx}\p_{\tx}\p_\mu^2 h(\mu)(\tx,\cx)
$$
are continuous at any $(\mu,\tx,\cx)$. We say $h\in C_b^{(2,2)}(\cP_2(\mR^n),\mR)$, if $h\in C^{(2,2)}(\cP_2(\mR^n),\mR)$, and all the above derivatives are uniformly bounded w.r.t. $(\mu,\tx,\cx)$.
For a matrix-valued map $h(\cdot):\cP_2(\mR^n)\rightarrow\mR^{n \times d_1}$, we say $h\in C^{(2,2)}(\cP_2(\mR^n),\mR^{n \times d_1})$ $(resp.~ C_b^{(2,2)}(\cP_2(\mR^n),\mR^{n \times d_1}))$ if all the components belong
to $C^{(2,2)}(\cP_2(\mR^n),\mR)$ $(resp.~ C_b^{(2,2)}(\cP_2(\mR^n),\mR))$.
\end{definition}

\begin{definition}\label{Ck1ddef}
For a matrix-valued map $h(\cdot,\cdot):\mR^n\times\cP_2(\mR^n)\rightarrow\mR^{n \times d_1}$, we say $h\in C^{k,(1,d)}(\mR^n\times\cP_2(\mR^n),\mR^{n \times d_1})$, if $h(x,\cdot)\in C^{(1,d)}(\cP_2(\mR^n),\mR^{n \times d_1})$ for any $x\in\mR^n$ and $h(\cdot,\mu)\in C^k(\mR^n,\mR^{n \times d_1})$ for any $\mu\in\cP_2(\mR^n)$, and the derivatives
$$
\p_x^ih(x,\mu),\quad 1\leq i\leq k,\quad \p_{\tx}^j\p_\mu h(x,\mu)(\tx),\quad 0\leq j\leq d,
$$
are continuous at any $(x,\mu,\tx)$.
We say $h\in C_b^{k,(1,d)}(\mR^n\times\cP_2(\mR^n),\mR^{n \times d_1})$, if $h\in C^{k,(1,d)}(\mR^n\times\cP_2(\mR^n),\mR^{n \times d_1})$, and all the above derivatives are uniformly bounded w.r.t. $(x,\mu,\tx)$.
\end{definition}

\begin{definition}\label{Ck22def}
For a matrix-valued map $h(\cdot,\cdot):\mR^n\times\cP_2(\mR^n)\rightarrow\mR^{n \times d_1}$, we say $h\in C^{k,(2,2)}(\mR^n\times\cP_2(\mR^n),\mR^{n \times d_1})$, if $h(x,\cdot)\in C^{(2,2)}(\cP_2(\mR^n),\mR^{n \times d_1})$ for any $x\in\mR^n$, $h(\cdot,\mu)\in C^k(\mR^n,\mR^{n \times d_1})$ for any $\mu\in\cP_2(\mR^n)$, and the derivatives
$$\p_x^ih(x,\mu),\quad \p_x^i\p_\mu h(x,\mu)(\tx),\quad\p_x^i\p_{\tx}\p_\mu h(x,\mu)(\tx),\quad 1\leq i\leq k,$$
$$\p_\mu^2 h(x,\mu)(\tx,\cx),\quad \p_{\tx}\p_\mu^2 h(x,\mu)(\tx,\cx),\quad \p_{\cx}\p_\mu^2 h(x,\mu)(\tx,\cx),\quad \p_{\cx}\p_{\tx}\p_\mu^2 h(x,\mu)(\tx,\cx)$$
are continuous at any $(x,\mu,\tx,\cx)$. We say $h\in C_b^{k,(2,2)}(\mR^n\times\cP_2(\mR^n),\mR^{n \times d_1})$, if $h\in C^{k,(2,2)}(\mR^n\times\cP_2(\mR^n),\mR^{n \times d_1})$, and all the above derivatives are uniformly bounded w.r.t. $(x,\mu,\tx,\cx)$.
\end{definition}

\begin{definition}\label{Ck1dldef}
For a matrix-valued map $h(\cdot,\cdot,\cdot):\mR^n\times\cP_2(\mR^n)\times\mR^m\rightarrow\mR^{n \times d_1}$, we say $h\in C^{k,(1,d),l}(\mR^n\times\cP_2(\mR^n)\times\mR^m,\mR^{n \times d_1})$, if $h(x,\cdot,y)\in C^{(1,d)}(\cP_2(\mR^n),\mR^{n \times d_1})$ for any $x\in\mR^n,y\in\mR^m$, and $h(\cdot,\mu,\cdot)\in C^{k,l}(\mR^n\times\mR^m,\mR^{n \times d_1})$ for any $\mu\in\cP_2(\mR^n)$, and the derivatives
$$\p_x^i\p_y^jh(x,\mu,y),\quad 0\leq i\leq k \ \text{and}\ 0\leq j\leq l \ \text{with}\ 1\leq i+j\leq l,$$
$$\p_{\tx}^i\p_\mu h(x,\mu,y)(\tx),\quad \p_y\p_{\tx}^i\p_\mu h(x,\mu,y)(\tx)\quad 0\leq i\leq d,$$
are continuous at any $(x,\mu,y,\tx)$. We say $h\in C_b^{k,(1,d),l}(\mR^n\times\cP_2(\mR^n)\times\mR^m,\mR^{n \times d_1})$, if $h\in C^{k,(1,d),l}(\mR^n\times\cP_2(\mR^n)\times\mR^m,\mR^{n \times d_1})$, and all the above derivatives are uniformly bounded w.r.t. $(x,\mu,y,\tx)$.
\end{definition}

\begin{definition}\label{Ck22ldef}
For a matrix-valued map $h(\cdot,\cdot,\cdot):\mR^n\times\cP_2(\mR^n)\times\mR^m\rightarrow\mR^{n \times d_1}$, we say $h\in C^{k,(2,2),l}(\mR^n\times\cP_2(\mR^n)\times\mR^m,\mR^{n \times d_1})$, if $h(x,\cdot,y)\in C^{(2,2)}(\cP_2(\mR^n),\mR^{n \times d_1})$ for any $x\in\mR^n,y\in\mR^m$ and $h(\cdot,\mu,\cdot)\in C^{k,l}(\mR^n\times\mR^m,\mR^{n \times d_1})$ for any $\mu\in\cP_2(\mR^n)$, and the derivatives
$$\p_x^i\p_y^jh(x,\mu,y),\quad \p_x^i\p_y^j\p_\mu h(x,\mu)(\tx),\quad\p_x^i\p_y^j\p_{\tx}\p_\mu h(x,\mu)(\tx),$$
$$\p_\mu^2 h(x,\mu,y)(\tx,\cx),\quad \p_{\tx}\p_\mu^2 h(x,\mu,y)(\tx,\cx),$$
$$\p_{\cx}\p_\mu^2 h(x,\mu,y)(\tx,\cx),\quad \p_{\cx}\p_{\tx}\p_\mu^2 h(x,\mu,y)(\tx,\cx),$$
$$\p_y\p_\mu^2 h(x,\mu,y)(\tx,\cx),\quad \p_y\p_{\tx}\p_\mu^2 h(x,\mu,y)(\tx,\cx),$$
$$\p_y\p_{\cx}\p_\mu^2 h(x,\mu,y)(\tx,\cx),\quad \p_y\p_{\cx}\p_{\tx}\p_\mu^2 h(x,\mu,y)(\tx,\cx)$$
are continuous at any $(x,\mu,y,\tx,\cx)$, where $0\leq i\leq k \ \text{and}\ 0\leq j\leq l \ \text{with}\ 1\leq i+j\leq l$. We say $h\in C_b^{k,(2,2),l}(\mR^n\times\cP_2(\mR^n)\times\mR^m,\mR^{n \times d_1})$, if $h\in C^{k,(2,2),l}(\mR^n\times\cP_2(\mR^n)\times\mR^m,\mR^{n \times d_1})$, and all the above derivatives are uniformly bounded w.r.t. $(x,\mu,y,\tx,\cx)$.
\end{definition}

\section{Main results}\label{main}

In this section, we formulate the main results in this paper.

\subsection{Nonautonomous multiscale McKean-Vlasov SDEs in the general case}
In this subsection, we present two averaging principles for nonautonomous multiscale McKean--Vlasov SDEs in the general setting.

We recall the system (\ref{orieq}):
\ce\left\{\begin{array}{l}
\dif X_t^{\eps}=b_1(X_t^{\eps}, \sL_{X_t^{\eps}}, Y_t^{\eps}) \dif t+\s_1(X_t^{\eps}, \sL_{X_t^{\eps}}, Y_t^{\eps}) \dif B_t, \quad X_0^{\eps}=\xi, \\
\dif Y_t^{\eps}=\eps^{-1} b_2(t/\eps,X_t^{\eps}, \sL_{X_t^{\eps}}, Y_t^{\eps})\dif t+\eps^{-1/2} \s_2(t/\eps,X_t^{\eps}, \sL_{X_t^{\eps}}, Y_t^{\eps})\dif W_t, \quad Y_0^{\eps}=\varrho,
\end{array}\right.
\de
where $\mE|\xi|^4<+\infty$ and $\mE|\varrho|^4<+\infty$.

Assume:
\begin{enumerate}[$(\mathbf{H}^1_{b_{1}, \s_{1}})$]
\item
$b_1\in C_b^{2,(1,1),3}(\mR^n\times\cP_2(\mR^n)\times\mR^m,\mR^n)$ and $\s_1$ is bounded. Moreover, there exists a constant $C>0$ such that for any $x_1,x_2\in \mR^n$, $\mu_1,\mu_2\in\cP_2(\mR^n)$, $y_1,y_2\in \mR^m$,
\be
&&|b_1(x_1,\mu_1,y_1)-b_1(x_2,\mu_2,y_2)|+\|\s_1(x_1,\mu_1,y_1)-\s_1(x_2,\mu_2,y_2)\|\no\\
&\leq&C\left(|x_1-x_2|+\mW_2(\mu_1, \mu_2)+|y_1-y_2|\right).\label{b1s1lip}
\ee
\end{enumerate}
\begin{enumerate}[$(\mathbf{H}^1_{b_{2}, \s_{2}})$]
\item
There exists a constant $C>0$ such that for any $t\in\mR$, $x_1, x_2 \in \mR^n$, $\mu_1, \mu_2 \in \cP_2(\mR^n)$, $y_1, y_2 \in \mR^m$,
\ce
 && 2\< y_1-y_2, b_2(t,x_1,\mu_1, y_1)-b_2(t,x_2,\mu_2, y_2)\>+3\|\s_2(t,x_1,\mu_1, y_1)-\s_2(t,x_2,\mu_2, y_2)\|^2\\
 &\leq& -2\a(t)|y_1-y_2|^2+C\a(t)(|x_1-x_2|^2+\mW_2^2(\mu_1, \mu_2)),
\de
where $\a:\mR\rightarrow(0,+\infty)$ satisfy for some $\g\in(0,1)$,
\be
\min\left\{\int_{-\infty}^0\a(u)\dif u, \int_0^{+\infty}\a(u)\dif u\right\}=+\infty,~ \Upsilon_{\g}(t):=\int_t^{+\infty}e^{-\g\int_t^r\a(u)\dif u}\dif r<+\infty, ~t\in\mR,\label{alpha}
\ee
and $\Upsilon(t):=\Upsilon_1(t)$.
\end{enumerate}
\begin{enumerate}[$(\mathbf{H}^2_{b_{2}, \s_{2}})$]
\item
For any $t\in\mR$, $b_2(t,\cdot,\cdot, \cdot)\in C^{2,(1,1),3}(\mR^n\times\cP_2(\mR^n)\times\mR^m,\mR^m)$, $\s_2(t,\cdot,\cdot, \cdot)\in C^{2,(1,1),3}(\mR^n\times\cP_2(\mR^n)\times\mR^m,\mR^{m \times d_2})$. Moreover, there exists a constant $C>0$ such that for any $t\in\mR$, $x\in \mR^n$, $\mu\in \cP_2(\mR^n)$, $y\in \mR^m$, $\tx\in \mR^n$, $i=0,1,2$ and $j=0,1,2,3$ with $1\leq i+j\leq 3$,
\ce
&&\|\p_x^i\p_y^jb_2(t,x,\mu, y)\|+\|\p_x^i\p_y^j\s_2(t,x,\mu, y)\|^2\leq C\a(t),\\
&&\|\p_{\mu}b_2(t,x,\mu, y)(\tx)\|+\|\p_{\mu}\s_2(t,x,\mu, y)(\tx)\|^2\leq C\a(t),\\
&&\|\p_{\tx}\p_{\mu}b_2(t,x,\mu, y)(\tx)\|+\|\p_{\tx}\p_{\mu}\s_2(t,x,\mu, y)(\tx)\|^2\leq C\a(t),\\
&&\|\p_y\p_{\mu}b_2(t,x,\mu, y)(\tx)\|+\|\p_y\p_{\mu}\s_2(t,x,\mu, y)(\tx)\|^2\leq C\a(t),\\
&&\|\p_y\p_{\tx}\p_{\mu}b_2(t,x,\mu, y)(\tx)\|+\|\p_y\p_{\tx}\p_{\mu}\s_2(t,x,\mu, y)(\tx)\|^2\leq C\a(t),\\
&&|b_2(t,x,\mu, y)|\leq C\a(t)\left(1+|x|^2+\mu(|\cdot|^2)+|y|^2\right)^{1/2},\\
&&\|\s_2(t,x,\mu, y)\|^2\leq C\a(t)\left(1+|x|^2+\mu(|\cdot|^2)+|y|^2\right),
\de
where $\a(t)$ is defined in $(\mathbf{H}^1_{b_{2}, \s_{2}})$.
\end{enumerate}
\begin{enumerate}[$(\mathbf{H}^3_{b_{2}, \s_{2}})$]
\item
$b_2$ and $\s_2$ satisfy $(\mathbf{H}^1_{b_{2}, \s_{2}})$ and $(\mathbf{H}^2_{b_{2}, \s_{2}})$ for any $t\geq0$, $x\in\mathbb{R}^n$, $\mu\in\cP_2(\mR^n)$ and $y\in\mathbb{R}^m$ with $\a:[0,+\infty)\rightarrow(0,+\infty)$ which satisfies for some $\g\in(0,1)$,
\be
\int_0^{+\infty}\a(u)\dif u=+\infty,\quad \int_t^{+\infty}e^{-\g\int_t^r\a(u)\dif u}\dif r<+\infty,\quad t\in[0,+\infty).\label{alphaposi}
\ee
\end{enumerate}

\br
$(i)$ By $(\mathbf{H}_{b_2, \s_2}^1)$ and $(\mathbf{H}_{b_2, \s_2}^2)$, one obtains that for any $\b\in(0,1)$ there exists a constant $C_{\b}>0$ such that for any $t\in\mR$, $x\in \mR^n$, $\mu\in \cP_2(\mR^n)$ and $y\in \mR^m$,
\be
2\la y, b_2(t,x,\mu,y)\ra+3\|\s_2(t,x,\mu,y)\|^2\leq -2\b\a(t)|y|^2+C_{\b}\a(t)\left(1+|x|^2+\mu(|\cdot|^2)\right).\label{b2s2dis}
\ee

$(ii)$ By $(\mathbf{H}^1_{b_{2}, \s_{2}})$, it follows that for any $y,k\in\mR^m$, $t\in\mR$, $x\in \mR^n$, $\mu\in \cP_2(\mR^n)$, and $\d>0$,
\ce
2\la b_2(t,x,\mu,y+\d k)-b_2(t,x,\mu,y), \d k\ra+3\|\s_2(t,x,\mu,y+\d k)-\s_2(t,x,\mu,y)\|^2\leq -2\a(t)|\d k|^2,
\de
and
\ce
&&2\la \d^{-1}(b_2(t,x,\mu,y+\d k)-b_2(t,x,\mu,y)), k\ra+3\|\d^{-1}(\s_2(t,x,\mu,y+\d k)-\s_2(t,x,\mu,y))\|^2\\
&\leq& -2\a(t)|k|^2,
\de
which yields that for any $y,k\in\mR^m$, $t\in\mR$, $x\in \mR^n$, $\mu\in \cP_2(\mR^n)$,
\be
2\la \p_y b_2(t,x,\mu,y)\cdot k, k\ra+3\|\p_y\s_2(t,x,\mu,y)\cdot k\|^2
\leq -2\a(t)|k|^2.
\label{parb2s2}
\ee
\er

Under $(\mathbf{H}^1_{b_{1}, \s_{1}})$ and $(\mathbf{H}^3_{b_{2}, \s_{2}})$, it follows from \cite[Theorem 2.1]{Wang} that the system (\ref{orieq}) admits a unique strong solution $(X_\cdot^{\eps},Y_\cdot^{\eps})$.

Next, we require that $\s_1$ is independent of $y$, define
\ce
\bar{b}_1(t,x,\mu):=\int_{\mR^m}b_1(x,\mu,y)\nu_t^{x,\mu}(\dif y),\quad t\geq0,
\de
where $\nu_t^{x,\mu}$ is from Lemma \ref{nuxmu}, and construct the following equation
\be
\dif \bar{X}_t^\eps=\bar{b}_1(t/\eps,\bar{X}_t^\eps,\sL_{\bar{X}_t^\eps})\dif t+\s_1(\bar{X}_t^\eps,\sL_{\bar{X}_t^\eps})\dif B_t,\quad \bar{X}_0^\eps=\xi.
\label{straveeq}
\ee
So, $\{\bar{X}_t^\eps\}_{t\geq0}$ denotes the unique solution of the above equation.

The following theorem establishes the strong averaging principle for the system (\ref{orieq}) in the general case.

\bt\label{xbarx2}
 Assume that $(\mathbf{H}^1_{b_{1}, \s_{1}})$ and $(\mathbf{H}^3_{b_{2}, \s_{2}})$ hold. Then for $T>0$, there exists $C_T>0$ such that for $\eps>0$,
\be
\mE(\sup_{0\leq t\leq T}|X_t^{\eps}-\bar{X}_t^\eps|^2)&\leq& C_T(1+\mE|\xi|^4+\mE|\varrho|^4)\eps^2\bigg[\sup_{0\leq t\leq T}\Upsilon^2_{\g}(\frac{t}{\eps})+\sup_{0\leq t\leq T}\Upsilon^2(\frac{t}{\eps})\left(\int_0^{\frac{T}{\eps}}\a(s)\dif s\right)^{\frac{1}{2}}\no\\
&&\qquad +\int_0^{\frac{T}{\eps}}\a(s)\Upsilon^2(s)\dif s \bigg].
\label{xebarxel2}
\ee
\et

\br\label{xbarx2rem}
Let $\a(t)=a_0(1+t)^\b$ with $a_0>0$ and $\b\in(-1,\infty)$. By \cite[Lemma 3.6]{swx}, it holds that
\ce
\sup_{0\leq t\leq T}\Upsilon^2_{\g}(t/\eps)\leq C\int_0^{T/\eps}\a(s)\Upsilon^2(s)\dif s.
\de
Combining this estimate with Theorem \ref{xbarx2}, we obtain
\ce
&&\mE(\sup_{0\leq t\leq T}|X_t^{\eps}-\bar{X}_t^\eps|^2)\\
&\leq& C_T(1+\mE|\xi|^4+\mE|\varrho|^4)\eps^2\int_0^{T/\eps}\a(s)\Upsilon^2(s)\dif s\left[1+\left(\int_0^{T/\eps}\a(s)\dif s\right)^{1/2} \right].
\de
Moreover, by \cite[Remark 3.7]{swx}, there exist two positive constants $a_1\leq a_2$ such that
\ce
a_1\int_0^t\a^{-1}(s)\dif s\leq\int_0^t\a(s)\Upsilon^2(s)\dif s\leq a_2\int_0^t\a^{-1}(s)\dif s,\quad t\geq0.
\de
Consequently, we deduce that
\ce
\left[\mE(\sup_{0\leq t\leq T}|X_t^{\eps}-\bar{X}_t^\eps|^2)\right]^{1/2}
\le
\begin{cases}
C\,\varepsilon^{(1+\beta)/4}, & -1 < \beta < 1, \\
C\,\varepsilon^{1/2} \left(\log 1/\varepsilon\right)^{1/2}, & \beta = 1, \\
C\,\varepsilon^{(3-\b)/4}, & \beta > 1.
\end{cases}
\de
This estimate shows that the convergence rate of the strong averaging principle depends explicitly on the function $\a(t)$ appearing in $(\mathbf{H}^3_{b_{2}, \s_{2}})$.
\er

In the following, we admit that $\s_1$ depends on $y$ and study  the weak averaging principle in the general setting. Assume:
\begin{enumerate}[$(\mathbf{H}^2_{b_{1}, \s_{1}})$]
\item $b_1\in(C_b^{4,(1,3),4}\cap C_b^{2,(2,2),2})(\mR^n\times\cP_2(\mR^n)\times\mR^m,\mR^n)$, $\s_1\in(C_b^{4,(1,3),4}\cap C_b^{2,(2,2),2})(\mR^n\times\cP_2(\mR^n)\times\mR^m,\mR^{n\times d_1})$ and $\s_1$ is bounded.
\end{enumerate}
\begin{enumerate}[$(\mathbf{H}^3_{b_{1}, \s_{1}})$]
\item
\ce
\inf_{x\in\mR^n,\mu\in\cP_2(\mR^n),y\in\mR^m,z\in\mR^n\backslash\{0\}}\frac{\<(\s_1\s_1^*)(x,\mu, y)\cdot z,z\>}{|z|^2}>0.
\de
\end{enumerate}
\begin{enumerate}[$(\mathbf{H}^4_{b_{2}, \s_{2}})$]
\item For any $t\geq0$, $b_2(t,\cdot,\cdot, \cdot)\in(C^{4,(1,3),5}\cap C^{2,(2,2),3})(\mR^n\times\cP_2(\mR^n)\times\mR^m,\mR^m)$, $\s_2(t,\cdot,\cdot, \cdot)\in(C^{4,(1,3),5}\cap C^{2,(2,2),3})(\mR^n\times\cP_2(\mR^n)\times\mR^m,\mR^{m \times d_2})$. Moreover, there exists a constant $C>0$ such that for any $t\geq0$, $x\in \mR^n$, $\mu\in \cP_2(\mR^n)$, $y\in \mR^m$, $\tx,\cx\in \mR^n$, the following conditions hold:

$(i)$ for any $0\leq i\leq4$ and $0\leq j\leq5$ with $1\leq i + j\leq 5$,
\ce
\|\p_x^i\p_y^jb_2(t,x,\mu, y)\|+\|\p_x^i\p_y^j\s_2(t,x,\mu, y)\|^2\leq C\a(t);
\de

$(ii)$ for any $0\leq i\leq3$,
\ce
&&\|\p_{\tx}^i\p_{\mu}b_2(t,x,\mu, y)(\tx)\|+\|\p_{\tx}^i\p_{\mu}\s_2(t,x,\mu, y)(\tx)\|^2\leq C\a(t),\\
&&\|\p_y\p_{\tx}^i\p_{\mu}b_2(t,x,\mu, y)(\tx)\|+\|\p_y\p_{\tx}^i\p_{\mu}\s_2(t,x,\mu, y)(\tx)\|^2\leq C\a(t);
\de

$(iii)$ for any $0\leq i\leq2$ and $0\leq j\leq3$ with $1\leq i + j\leq 3$,
\ce
&&\|\p_x^i\p_y^j\p_{\mu}b_2(t,x,\mu, y)(\tx)\|+\|\p_x^i\p_y^j\p_{\mu}\s_2(t,x,\mu, y)(\tx)\|^2\leq C\a(t),\\
&&\|\p_x^i\p_y^j\p_{\tx}\p_{\mu}b_2(t,x,\mu, y)(\tx)\|+\|\p_x^i\p_y^j\p_{\tx}\p_{\mu}\s_2(t,x,\mu, y)(\tx)\|^2\leq C\a(t);
\de

$(iv)$ for any $i=0,1$ and $j=0,1$,
\ce
&&\|\p_{\tx}^i\p_{\cx}^j\p_{\mu}^2b_2(t,x,\mu, y)(\tx,\cx)\|+\|\p_{\tx}^i\p_{\cx}^j\p_{\mu}^2\s_2(t,x,\mu, y)(\tx,\cx)\|^2\leq C\a(t),\\
&&\|\p_y\p_{\tx}^i\p_{\cx}^j\p_{\mu}^2b_2(t,x,\mu, y)(\tx,\cx)\|+\|\p_y\p_{\tx}^i\p_{\cx}^j\p_{\mu}^2\s_2(t,x,\mu, y)(\tx,\cx)\|^2\leq C\a(t),
\de
where $\a(t)$ is defined in $(\mathbf{H}^3_{b_{2}, \s_{2}})$.
\end{enumerate}

Note that $(\mathbf{H}^2_{b_{1}, \s_{1}})$ is stronger than $(\mathbf{H}^1_{b_{1}, \s_{1}})$. Thus, under $(\mathbf{H}^2_{b_{1}, \s_{1}})$ and $(\mathbf{H}^3_{b_{2}, \s_{2}})$, the system (\ref{orieq}) has a unique solution $(X_\cdot^{\eps},Y_\cdot^{\eps})$.

Consider the following averaging equation:
\be
\dif \check{X}_t^\eps=\bar{b}_1(t/\eps,\check{X}_t^\eps,\sL_{\check{X}_t^\eps})\dif t+\bar{\s}_1(t/\eps,\check{X}_t^\eps,\sL_{\check{X}_t^\eps})\dif \bar{B}_t,\quad \check{X}_0^\eps=\xi,
\label{weaveeq}
\ee
where
\ce
\bar{\s}_1(t,x,\mu):=[\overline{\s_1\s_1^*}(t,x,\mu)]^{1/2}:=\Big(\int_{\mR^m}(\s_1\s_1^*)(x,\mu,y)\nu_t^{x,\mu}(\dif y)\Big)^{1/2},\quad t\geq0,
\de
and $\{\bar{B}_t\}_{t\geq0}$ is a $n$-dimensional standard Brownian motion.

We are now in a position to state the weak averaging principle in the general setting.
\bt\label{xbarweak}
Suppose that $(\mathbf{H}^2_{b_{1}, \s_{1}})$, $(\mathbf{H}^3_{b_{1}, \s_{1}})$, $(\mathbf{H}^3_{b_{2}, \s_{2}})$ and $(\mathbf{H}^4_{b_{2}, \s_{2}})$ hold. Then for any $T>0$ and $\varphi\in (C_b^{(1,3)}\cap C_b^{(2,2)})(\cP_2(\mR^n),\mR)$, there exists a constant $C_T>0$ such that for $\eps>0$,
\be
\sup_{0\leq t\leq T}|\varphi(\sL_{X_t^{\eps}})-\varphi(\sL_{\check{X}_t^\eps})|\leq C_T(1+\mE|\xi|^4+\mE|\varrho|^4)\eps\sup_{0\leq t\leq T}\Upsilon_{\g}(t/\eps),\label{lxelbarxe}
\ee
where $\{\check{X}_t^\eps\}_{t\geq0}$ is the unique solution to (\ref{weaveeq}).
\et

\br\label{xbarweakrem}
Let $\a(t)=a_0(1+t)^\b$ with $a_0>0$ and $\b\in(-1,+\infty)$. Under the assumptions of Theorem \ref{xbarweak}, it follows that
\be\label{xbarweakrate}
\sup_{0\leq t\leq T}\big|\varphi(\sL_{X_t^{\eps}})-\varphi(\sL_{\check{X}_t^\eps})\big|
\le
\begin{cases}
C\,\varepsilon^{1+\beta}, & -1<\beta<0,\\
C\,\varepsilon, & \beta\geq 0.
\end{cases}
\ee
This estimate indicates that the convergence rate of the weak averaging principle depends explicitly on the function $\a(t)$ appearing in $(\mathbf{H}^3_{b_{2}, \s_{2}})$.
\er

\br
Note that the coefficients in (\ref{straveeq}) and (\ref{weaveeq}) still depend on the parameter $\eps$. To derive averaging equations whose coefficients are independent of $\eps$, we further impose additional assumptions on the time behavior of the fast coefficients. In particular, we assume that $b_{2}$ and $\sigma_{2}$ exhibit convergent and periodic behaviors with respect to the time variable, namely $(\mathbf{H}^5_{b_{2}, \sigma_{2}})$ and $(\mathbf{H}^6_{b_{2}, \sigma_{2}})$. The corresponding averaging principles for these two cases will be established in Subsections \ref{aveconm} and \ref{aveperm}.
\er

The proofs of Theorem \ref{xbarx2} and \ref{xbarweak} are placed in Section \ref{geneproo}.

\subsection{Nonautonomous multiscale McKean-Vlasov SDEs in the convergent case}\label{aveconm}

In this subsection, we present two averaging principles for nonautonomous multiscale McKean-Vlasov SDEs in the convergent case.

Assume:
\begin{enumerate}[$(\mathbf{H}^5_{b_{2}, \s_{2}})$]
\item $\limsup\limits_{t\rightarrow+\infty}\alpha(t)=\alpha>0$, and there exist functions $\bar{b}_2:\mR^n \times \cP_2(\mR^n) \times \mR^m\rightarrow \mR^m$, and
$\bar{\s}_2: \mR^n \times \cP_2(\mR^n) \times \mR^m\rightarrow \mR^{m \times d_2}$ such that
\be
&&|b_2(t,x,\mu,y)-\bar{b}_2(x,\mu,y)|+\|\s_2(t,x,\mu,y)-\bar{\s}_2(x,\mu,y)\|\no\\
&\leq&\lambda(t)\left(1+|x|+\mu^{1/2}(|\cdot|^2)+|y|\right),\label{conass}
\ee
where $\lambda:[0,+\infty)\rightarrow(0,+\infty)$ is locally bounded and satisfies that $\lim\limits_{t\rightarrow+\infty}\lambda(t)=0$.
\end{enumerate}

We require that $\s_1$ is independent of $y$ and construct the following SDE
\be
\dif \bar{X}_t=\bar{b}_{1,c}(\bar{X}_t,\sL_{\bar{X}_t})\dif t+\s_1(\bar{X}_t,\sL_{\bar{X}_t})\dif B_t,\quad \bar{X}_0=\xi,
\label{constraveeq}
\ee
where
\ce
\bar{b}_{1,c}(x,\mu):=\int_{\mR^m}b_1(x,\mu,y)\nu^{x,\mu}(\dif y),
\de
and $\nu^{x,\mu}$ is the unique invariant measure of Eq.(\ref{confrozeq}). And $\{\bar{X}_t\}_{t\geq0}$ is the solution to Eq.(\ref{constraveeq}).

The following theorem establishes the strong averaging principle in the convergent case.
\bt\label{xbarx2con}
Suppose that $(\mathbf{H}^1_{b_{1}, \s_{1}})$, $(\mathbf{H}^3_{b_{2}, \s_{2}})$ and $(\mathbf{H}^5_{b_{2}, \s_{2}})$ hold. Then for $T>0$ and $\beta\in(0,1)$, there exists $C_{T,\beta}>0$ such that for $\eps>0$,
\be
\mE(\sup_{0\leq t\leq T}|X_t^{\eps}-\bar{X}_t|^2)&\leq& C_{T,\beta}(1+\mE|\xi|^4+\mE|\varrho|^4)\eps^2\Bigg[\sup_{0\leq t\leq T}\Upsilon^2_{\g}(\frac{t}{\eps})+\sup_{0\leq t\leq T}\Upsilon^2(\frac{t}{\eps})\left(\int_0^{\frac{T}{\eps}}\a(s)\dif s\right)^{\frac{1}{2}}\no\\
&&+\int_0^{\frac{T}{\eps}}\a(s)\Upsilon^2(s)\dif s+\Bigg(1+\int_0^{\frac{T}{\eps}}\Big(\int_0^ue^{-2\beta\alpha(u-r)}\lambda^2(r)\dif r\Big)^{1/2}\dif u\Bigg)^2\Bigg].\no\\
\label{xebarxl2con}
\ee
\et

Next, we admit that $\s_1$ depends on $y$ and consider the following averaging equation:
\be
\dif \check{X}_t=\bar{b}_{1,c}(\check{X}_t,\sL_{\check{X}_t})\dif t+\bar{\s}_{1,c}(\check{X}_t,\sL_{\check{X}_t})\dif \bar{B}_t,\quad \check{X}_0=\xi,
\label{conweaveeq}
\ee
where
\ce
\bar{\s}_{1,c}(x,\mu):=[(\overline{\s_1\s_1^*})_c(x,\mu)]^{1/2}:=\Big(\int_{\mR^m}(\s_1\s_1^*)(x,\mu,y)\nu^{x,\mu}(\dif y)\Big)^{1/2}.
\de

The weak averaging principle for the convergent case is given below.
\bt\label{xbarweakcon}
Suppose that $(\mathbf{H}^2_{b_{1}, \s_{1}})$, $(\mathbf{H}^3_{b_{1}, \s_{1}})$, $(\mathbf{H}^3_{b_{2}, \s_{2}})$, $(\mathbf{H}^4_{b_{2}, \s_{2}})$ and $(\mathbf{H}^5_{b_{2}, \s_{2}})$ hold. Then for any $T>0$ and $\varphi\in (C_b^{(1,3)}\cap C_b^{(2,2)})(\cP_2(\mR^n),\mR)$, there exists a constant $C_{T,\beta}>0$ such that for $\eps>0$,
\be
&&\sup_{0\leq t\leq T}|\varphi(\sL_{X_t^{\eps}})-\varphi(\sL_{\check{X}_t})|\no\\
&\leq& C_{T,\beta}(1+\mE|\xi|^4+\mE|\varrho|^4)\eps\Bigg[1+\sup_{0\leq t\leq T}\Upsilon_{\g}(t/\eps)+\int_0^{T/\eps}\Big(\int_0^ue^{-2\beta\alpha(u-r)}\lambda^2(r)\dif r\Big)^{1/2}\dif u\Bigg],\no\\
\label{lxelbarxcon}
\ee
where $\{\check{X}_t\}_{t\geq0}$ is the solution to Eq.(\ref{conweaveeq}).
\et

The proofs of Theorem \ref{xbarx2con} and \ref{xbarweakcon} are postponed to Section \ref{convcoefproo}.

\subsection{Nonautonomous multiscale McKean-Vlasov SDEs in the periodic case}\label{aveperm}

In this subsection, we present two averaging principles for nonautonomous multiscale McKean-Vlasov SDEs in the periodic case.

Assume:
\begin{enumerate}[$(\mathbf{H}^6_{b_{2}, \s_{2}})$]
\item
$b_2(\cdot,x,\mu,y)$ and $\s_2(\cdot,x,\mu,y)$ are $\tau$-periodic for some $\tau>0$, that is, for any $t\in\mR, x\in\mR^n, \mu\in\cP_2(\mR^n)$ and $y\in\mR^m$,
\be
b_2(t+\tau,x,\mu,y)=b_2(t,x,\mu,y),\quad \s_2(t+\tau,x,\mu,y)=\s_2(t,x,\mu,y).\label{perass}
\ee
\end{enumerate}

We require that $\s_1$ is independent of $y$ and consider the following averaging SDE:
\be
\dif \tilde{X}_t=\bar{b}_{1,p}(\tilde{X}_t,\sL_{\tilde{X}_t})\dif t+\s_1(\tilde{X}_t,\sL_{\tilde{X}_t})\dif B_t,\quad \tilde{X}_0=\xi,
\label{perstraveeq}
\ee
where
\ce
\bar{b}_{1,p}(x,\mu):=\frac1\tau\int_0^\tau\bar{b}_1(t,x,\mu)\dif t=\frac1\tau\int_0^\tau\int_{\mR^m}b_1(x,\mu,y)\nu_t^{x,\mu}(\dif y)\dif t.
\de
Then $\{\tilde{X}_t\}_{t\geq0}$ denotes the unique solution of the above equation.

The following theorem establishes the strong averaging principle in the periodic case.
\bt\label{xbarstrper}
Suppose that $(\mathbf{H}^1_{b_{1}, \s_{1}})$, $(\mathbf{H}^3_{b_{2}, \s_{2}})$ and $(\mathbf{H}^6_{b_{2}, \s_{2}})$ hold. Then for $T>0$ and $\beta\in(0,1)$, there exists $C_{T}>0$ such that for $\eps>0$,
\be
\mE(\sup_{0\leq t\leq T}|X_t^{\eps}-\tilde{X}_t|^2)&\leq& C_T(1+\mE|\xi|^4+\mE|\varrho|^4)\eps^2\bigg[\eps^{-4/3}+\sup_{0\leq t\leq T}\Upsilon^2_{\g}(t/\eps)\no\\
&&+\sup_{0\leq t\leq T}\Upsilon^2(\frac{t}{\eps})\left(\int_0^{T/\eps}\a(s)\dif s\right)^{1/2}+\int_0^{T/\eps}\a(s)\Upsilon^2(s)\dif s \bigg].
\label{xebarxl2per}
\ee
\et

Finally, we admit that $\s_1$ depends on $y$ and construct the following equation
\be
\dif \hat{X}_t=\bar{b}_{1,p}(\hat{X}_t,\sL_{\hat{X}_t})\dif t+\bar{\s}_{1,p}(\hat{X}_t,\sL_{\hat{X}_t})\dif \bar{B}_t,\quad \hat{X}_0=\xi,
\label{perweaveeq}
\ee
where
\ce
\bar{\s}_{1,p}(x,\mu):=[(\overline{\s_1\s_1^*})_p(x,\mu)]^{1/2}:=\Big(\frac1\tau\int_0^\tau\int_{\mR^m}(\s_1\s_1^*)(x,\mu,y)\nu_t^{x,\mu}(\dif y)\dif t\Big)^{1/2},
\de
and present the weak averaging principle for the periodic case.

\bt\label{xbarweakper}
Suppose that $(\mathbf{H}^2_{b_{1}, \s_{1}})$, $(\mathbf{H}^3_{b_{1}, \s_{1}})$, $(\mathbf{H}^3_{b_{2}, \s_{2}})$, $(\mathbf{H}^4_{b_{2}, \s_{2}})$ and $(\mathbf{H}^6_{b_{2}, \s_{2}})$ hold. Then for any $T>0$ and $\varphi\in (C_b^{(1,3)}\cap C_b^{(2,2)})(\cP_2(\mR^n),\mR)$, there exists a constant $C_T>0$ such that for $\eps>0$,
\be
\sup_{0\leq t\leq T}|\varphi(\sL_{X_t^{\eps}})-\varphi(\sL_{\hat{X}_t})|\leq C_T(1+\mE|\xi|^4+\mE|\varrho|^4)\eps\left[\sup_{0\leq t\leq T}\Upsilon_{\g}(t/\eps)+\eps^{-2/3}\right],
\label{lxelbarxper}
\ee
where $\{\hat{X}_t\}_{t\geq0}$ is the solution of Eq.(\ref{perweaveeq}).
\et

The proofs of Theorem \ref{xbarstrper} and \ref{xbarweakper} are placed in Section \ref{periconvproo}.

\section{Time-inhomogeneous frozen SDEs and nonautonomous Poisson equations}\label{frozpoi}
This section is devoted to the study of the time-inhomogeneous frozen dynamics and the associated nonautonomous Poisson equation. In Subsection \ref{frozsde}, we establish regularity estimates for the frozen SDE. Subsection \ref{evomeas} is concerned with the existence and uniqueness of the evolution system of measures. In Subsection \ref{nonpoieq}, we investigate the nonautonomous Poisson equation, derive a representation formula for its solution, and establish the regularity estimates required for the averaging analysis.

\subsection{Time-inhomogeneous frozen SDEs}\label{frozsde}
In this subsection, we derive moment and regularity estimates for the frozen SDE.

For fixed $x\in\mR^n$ and $\mu \in \cP_2(\mR^n)$, consider the time-inhomogeneous frozen equation
\be
\dif Y_t^{s,x,\mu,y}=b_2(t,x,\mu,Y_t^{s,x,\mu,y})\dif t+\s_2(t,x,\mu,Y_t^{s,x,\mu,y})\dif W_t, \quad Y_s^{s,x,\mu,y}=y.\label{frozeq}
\ee
$(\mathbf{H}^2_{b_{2}, \s_{2}})$ implies that the coefficients $b_2(t,x,\mu,\cdot)$ and $\sigma_2(t,x,\mu,\cdot)$ are globally Lipschitz in $y$ and satisfy the linear growth condition. Then Eq.(\ref{frozeq}) has a unique strong solution $Y_\cdot^{s,x,\mu,y}$.

\bl\label{yt}
Suppose that $(\mathbf{H}^1_{b_{2}, \s_{2}})$ and $(\mathbf{H}^2_{b_{2}, \s_{2}})$ hold. Then there exists a positive constant $C>0$ such that

$(i)$ for $t\geq s$, $x\in \mR^n$, $\mu\in \cP_2(\mR^n)$ and $y\in \mR^m$,
\be
\mE|Y_t^{s,x,\mu,y}|^4\leq e^{-4\g\int_s^t\a(u)\dif u}|y|^4+C\left(1+|x|^2+\mu(|\cdot|^2)\right)^2;\label{Yl4}
\ee

$(ii)$ for $t\geq s$, $x_1, x_2 \in \mR^n$, $\mu_1, \mu_2 \in \cP_2(\mR^n)$, $y_1, y_2 \in \mR^m$,
\be
\mE|Y_t^{s,x_1,\mu_1,y_1}-Y_t^{s,x_2,\mu_2,y_2}|^4\leq e^{-4\g\int_s^t\a(u)\dif u}|y_1-y_2|^4+C\left(|x_1-x_2|^4+\mW_2^4(\mu_1, \mu_2)\right);\label{Y1Y2l4}
\ee
where $\g\in(0,1)$ is the constant in (\ref{alpha}).
\el
\begin{proof}
Applying the It\^{o} formula to $|Y_t^{s,x,\mu,y}|^4$ and taking expectation yields
\ce
\mE|Y_t^{s,x,\mu,y}|^4
&=&|y|^4+4\mE\int_s^t|Y_r^{s,x,\mu,y}|^2\<Y_r^{s,x,\mu,y}, b_2(r,x,\mu,Y_r^{s,x,\mu,y})\>\dif r\\
&&+4\mE\int_s^t\|\s_2(r,x,\mu,Y_r^{s,x,\mu,y})Y_r^{s,x,\mu,y}\|^2\dif r\\
&&+2\mE\int_s^t|Y_r^{s,x,\mu,y}|^2\|\s_2(r,x,\mu,Y_r^{s,x,\mu,y})\|^2\dif r.
\de
By (\ref{b2s2dis}), it follows that
\ce
&&\frac{\dif}{\dif t}\mE|Y_t^{s,x,\mu,y}|^4\\
&\leq&4\mE[|Y_t^{s,x,\mu,y}|^2\<Y_t^{s,x,\mu,y}, b_2(t,x,\mu,Y_t^{s,x,\mu,y})\>]+6\mE[|Y_t^{s,x,\mu,y}|^2\|\s_2(t,x,\mu,Y_t^{s,x,\mu,y})\|^2]\\
&\leq&-4\g\a(t)\mE|Y_t^{s,x,\mu,y}|^4+C\a(t)\left(1+|x|^2+\mu(|\cdot|^2)\right)^2.
\de
The comparison theorem implies
\ce
\mE|Y_t^{s,x,\mu,y}|^4\leq e^{-4\g\int_s^t\a(u)\dif u}|y|^4+C\left(1+|x|^2+\mu(|\cdot|^2)\right)^2.
\de
Thus, $(i)$ is proved.

Next, we deal with $(ii)$. To estimate $\mE|Y_t^{s,x_1,\mu_1,y_1}-Y_t^{s,x_2,\mu_2,y_2}|^4$, we apply the It\^{o} formula to obtain
\ce
&&\mE|Y_t^{s,x_1,\mu_1,y_1}-Y_t^{s,x_2,\mu_2,y_2}|^4\\
&=&|y_1-y_2|^4+4\mE\int_s^t|Y_r^{s,x_1,\mu_1,y_1}-Y_r^{s,x_2,\mu_2,y_2}|^2\<Y_r^{s,x_1,\mu_1,y_1}-Y_r^{s,x_2,\mu_2,y_2}, \\
&&\qquad\qquad\qquad\qquad b_2(r,x_1,\mu_1,Y_r^{s,x_1,\mu_1,y_1})-b_2(r,x_2,\mu_2,Y_r^{s,x_2,\mu_2,y_2})\>\dif r\\
&&+4\mE\int_s^t\|(\s_2(r,x_1,\mu_1,Y_r^{s,x_1,\mu_1,y_1})-\s_2(r,x_2,\mu_2,Y_r^{s,x_2,\mu_2,y_2}))(Y_r^{s,x_1,\mu_1,y_1}-Y_r^{s,x_2,\mu_2,y_2})\|^2\dif r\\
&&+2\mE\int_s^t|Y_r^{s,x_1,\mu_1,y_1}-Y_r^{s,x_2,\mu_2,y_2}|^2\|\s_2(r,x_1,\mu_1,Y_r^{s,x_1,\mu_1,y_1})-\s_2(r,x_2,\mu_2,Y_r^{s,x_2,\mu_2,y_2})\|^2\dif r.
\de
By $(\mathbf{H}^1_{b_{2}, \s_{2}})$, we have
\ce
&&\frac{\dif}{\dif t}\mE|Y_t^{s,x_1,\mu_1,y_1}-Y_t^{s,x_2,\mu_2,y_2}|^4\\
&\leq&4\mE[|Y_t^{s,x_1,\mu_1,y_1}-Y_t^{s,x_2,\mu_2,y_2}|^2\\
&&\qquad \cdot\<Y_t^{s,x_1,\mu_1,y_1}-Y_t^{s,x_2,\mu_2,y_2}, b_2(t,x_1,\mu_1,Y_t^{s,x_1,\mu_1,y_1})-b_2(t,x_2,\mu_2,Y_t^{s,x_2,\mu_2,y_2})\>]\\
&&+6\mE[|Y_t^{s,x_1,\mu_1,y_1}-Y_t^{s,x_2,\mu_2,y_2}|^2\|\s_2(t,x_1,\mu_1,Y_t^{s,x_1,\mu_1,y_1})-\s_2(t,x_2,\mu_2,Y_t^{s,x_2,\mu_2,y_2})\|^2]\\
&\leq&2\mE\big[|Y_t^{s,x_1,\mu_1,y_1}-Y_t^{s,x_2,\mu_2,y_2}|^2\cdot\big(-2\a(t)|Y_t^{s,x_1,\mu_1,y_1}-Y_t^{s,x_2,\mu_2,y_2}|^2\\
&&\qquad\qquad\qquad\qquad\qquad\qquad\qquad\qquad+C\a(t)(|x_1-x_2|^2+\mW_2^2(\mu_1, \mu_2))\big)\big]\\
&\leq&-4\g\a(t)\mE|Y_t^{s,x_1,\mu_1,y_1}-Y_t^{s,x_2,\mu_2,y_2}|^4+C\a(t)(|x_1-x_2|^4+\mW_2^4(\mu_1, \mu_2)).
\de
The comparison theorem yields
\ce
\mE|Y_t^{s,x_1,\mu_1,y_1}-Y_t^{s,x_2,\mu_2,y_2}|^4\leq e^{-4\g\int_s^t\a(u)\dif u}|y_1-y_2|^4+C(|x_1-x_2|^4+\mW_2^4(\mu_1, \mu_2)).
\de
This completes the proof.
\end{proof}

\bl\label{pary}
Suppose that $(\mathbf{H}^1_{b_{2}, \s_{2}})$ and $(\mathbf{H}^2_{b_{2}, \s_{2}})$ hold. Then for any $t\geq s$,
\be\label{parYl4}
&&\mE\|\p_yY_t^{s,x,\mu,y}\|^4\leq e^{-4\int_s^t\a(u)\dif u},\quad \mE\|\p_y^2Y_t^{s,x,\mu,y}\|^2\leq e^{-2\g\int_s^t\a(u)\dif u},\no\\
&&\sup_{t\geq s}\mE\|\p_xY_t^{s,x,\mu,y}\|^4\leq C,\quad \sup_{t\geq s}\mE\|\p_x^2Y_t^{s,x,\mu,y}\|^2\leq C,\no\\
&&\sup_{t\geq s}\mE\|\p_{\mu}Y_t^{s,x,\mu,y}(\tx)\|^2\leq C,\quad \sup_{t\geq s}\mE\|\p_{\tx}\p_{\mu}Y_t^{s,x,\mu,y}(\tx)\|^2\leq C,
\ee
where $\g\in(0,1)$ is the constant in (\ref{alpha}).
\el
\begin{proof}
The estimates for derivatives with respect to $x$ and $y$ follow from arguments similar to those in \cite[Lemma 2.3]{swx}. It remains to estimate the derivative with respect to the measure variable. Note that $\p_{\mu}Y_t^{s,x,\mu,y}(\tx)$ satisfies
\ce\left\{\begin{array}{l}
\dif \p_{\mu}Y_t^{s,x,\mu,y}(\tx)=[\p_{\mu}b_2(t,x,\mu,Y_t^{s,x,\mu,y})(\tx)+\p_yb_2(t,x,\mu,Y_t^{s,x,\mu,y})\cdot\p_{\mu}Y_t^{s,x,\mu,y}(\tx)]\dif t\\
\qquad\qquad\qquad\qquad+[\p_{\mu}\s_2(t,x,\mu,Y_t^{s,x,\mu,y})(\tx)+\p_y\s_2(t,x,\mu,Y_t^{s,x,\mu,y})\cdot\p_{\mu}Y_t^{s,x,\mu,y}(\tx)]\dif W_t,\\
\p_{\mu}Y_s^{s,x,\mu,y}(\tx)=0.
\end{array}\right.
\de
Applying the It\^{o} formula, we deduce that
\ce
&&\frac{\dif}{\dif t}\mE\|\p_{\mu}Y_t^{s,x,\mu,y}(\tx)\|^2\\
&=&2\mE\<\p_{\mu}Y_t^{s,x,\mu,y}(\tx), \p_{\mu}b_2(t,x,\mu,Y_t^{s,x,\mu,y})(\tx)+\p_yb_2(t,x,\mu,Y_t^{s,x,\mu,y})\cdot\p_{\mu}Y_t^{s,x,\mu,y}(\tx)\>\\
&&+\mE\|\p_{\mu}\s_2(t,x,\mu,Y_t^{s,x,\mu,y})(\tx)+\p_y\s_2(t,x,\mu,Y_t^{s,x,\mu,y})\cdot\p_{\mu}Y_t^{s,x,\mu,y}(\tx)\|^2\\
&\leq&2\mE\<\p_{\mu}Y_t^{s,x,\mu,y}(\tx), \p_{\mu}b_2(t,x,\mu,Y_t^{s,x,\mu,y})(\tx)\>\\
&&+2\mE\<\p_{\mu}Y_t^{s,x,\mu,y}(\tx), \p_yb_2(t,x,\mu,Y_t^{s,x,\mu,y})\cdot\p_{\mu}Y_t^{s,x,\mu,y}(\tx)\>\\
&&+2\mE\|\p_{\mu}\s_2(t,x,\mu,Y_t^{s,x,\mu,y})(\tx)\|^2+2\mE\|\p_y\s_2(t,x,\mu,Y_t^{s,x,\mu,y})\cdot\p_{\mu}Y_t^{s,x,\mu,y}(\tx)\|^2\\
&\leq&-2\g\a(t)\mE\|\p_{\mu}Y_t^{s,x,\mu,y}(\tx)\|^2+C\a(t).
\de
The comparison theorem implies
\ce
\sup_{t\geq s}\mE\|\p_{\mu}Y_t^{s,x,\mu,y}(\tx)\|^2\leq C.
\de
By analogous arguments, we further obtain
\ce
\sup_{t\geq s}\mE\|\p_{\tx}\p_{\mu}Y_t^{s,x,\mu,y}(\tx)\|^2\leq C.
\de
The proof is complete.
\end{proof}

\subsection{Evolution system of measures}\label{evomeas}
In this subsection, we establish the existence and uniqueness of the evolution system of measures.
\bl\label{yvs}
Suppose that $(\mathbf{H}^1_{b_{2}, \s_{2}})$ and $(\mathbf{H}^2_{b_{2}, \s_{2}})$ hold. Then for any $t\in\mR$, $x\in \mR^n$ and $\mu\in \cP_2(\mR^n)$, there exists $\varsigma_t^{x,\mu}\in L^2(\Omega,\sF_t,\mP;\mR^m)$ such that for all $t\geq s$, $x\in \mR^n$, $\mu\in \cP_2(\mR^n)$ and $y\in \mR^m$,
\be\label{yvsl2}
&&\mE|Y_t^{s,x,\mu,y}-\varsigma_t^{x,\mu}|^2\leq C\left(1+|x|^2+\mu(|\cdot|^2)+|y|^2\right)e^{-2\int_s^t\a(u)\dif u},\no\\
&&\sup_{t\in\mR}\mE|\varsigma_t^{x,\mu}|^2\leq C\left(1+|x|^2+\mu(|\cdot|^2)\right),
\ee
and $\{\varsigma_t^{x,\mu}\}_{t\in\mR}$ solves
\be\label{vseq}
\varsigma_t^{x,\mu}=\varsigma_s^{x,\mu}+\int_s^tb_2(r,x,\mu,\varsigma_r^{x,\mu})\dif r+\int_s^t\s_2(r,x,\mu,\varsigma_r^{x,\mu})\dif W_r,\quad t\geq s.
\ee
Moreover, for all $t\geq s$, $x\in \mR^n$, $\mu\in \cP_2(\mR^n)$, $y\in \mR^m$ and $\tx\in\mR^n$,
\be\label{pyvs}
&&\mE\|\p_xY_t^{s,x,\mu,y}-\p_x\varsigma_t^{x,\mu}\|^2\leq C\left(1+|x|^2+\mu(|\cdot|^2)+|y|^2\right)e^{-2\g\int_s^t\a(u)\dif u},\no\\
&&\mE\|\p_x^2Y_t^{s,x,\mu,y}-\p_x^2\varsigma_t^{x,\mu}\|^2\leq C\left(1+|x|^2+\mu(|\cdot|^2)+|y|^2\right)e^{-2\g\int_s^t\a(u)\dif u},\no\\
&&\mE\|\p_{\mu}Y_t^{s,x,\mu,y}(\tx)-\p_{\mu}\varsigma_t^{x,\mu}(\tx)\|^2\leq C\left(1+|x|^2+\mu(|\cdot|^2)+|y|^2\right)e^{-2\g\int_s^t\a(u)\dif u},\no\\
&&\mE\|\p_{\tx}\p_{\mu}Y_t^{s,x,\mu,y}(\tx)-\p_{\tx}\p_{\mu}\varsigma_t^{x,\mu}(\tx)\|^2\leq C\left(1+|x|^2+\mu(|\cdot|^2)+|y|^2\right)e^{-2\g\int_s^t\a(u)\dif u},\no\\
&&\sup_{t\in\mR}\mE\|\p_x\varsigma_t^{x,\mu}\|^2\leq C, \quad \sup_{t\in\mR}\mE\|\p_x^2\varsigma_t^{x,\mu}\|^2\leq C,\no\\
&&\sup_{t\in\mR}\mE\|\p_{\mu}\varsigma_t^{x,\mu}(\tx)\|^2\leq C, \quad \sup_{t\in\mR}\mE\|\p_{\tx}\p_{\mu}\varsigma_t^{x,\mu}(\tx)\|^2\leq C,
\ee
where $\g\in(0,1)$ is the constant in (\ref{alpha}).
\el
\begin{proof}
First of all, we prove (\ref{yvsl2}) and (\ref{vseq}).

For any $h>0$, $t\geq s$, $x\in \mR^n$, $\mu\in \cP_2(\mR^n)$ and $y\in \mR^m$, we have
\be\label{yshxmuy}
Y_t^{s-h,x,\mu,y}=Y_s^{s-h,x,\mu,y}+\int_s^tb_2(r,x,\mu,Y_r^{s-h,x,\mu,y})\dif r+\int_s^t\s_2(r,x,\mu,Y_r^{s-h,x,\mu,y})\dif W_r,
\ee
and
\ce
Y_t^{s,x,\mu,y}=y+\int_s^tb_2(r,x,\mu,Y_r^{s,x,\mu,y})\dif r+\int_s^t\s_2(r,x,\mu,Y_r^{s,x,\mu,y})\dif W_r.
\de
Consequently,
\ce
Y_t^{s-h,x,\mu,y}-Y_t^{s,x,\mu,y}
&=&Y_s^{s-h,x,\mu,y}-y+\int_s^t[b_2(r,x,\mu,Y_r^{s-h,x,\mu,y})-b_2(r,x,\mu,Y_r^{s,x,\mu,y})]\dif r\\
&&+\int_s^t[\s_2(r,x,\mu,Y_r^{s-h,x,\mu,y})-\s_2(r,x,\mu,Y_r^{s,x,\mu,y})]\dif W_r,
\de
Applying the It\^{o} formula together with $(\mathbf{H}^1_{b_{2}, \s_{2}})$ yields
\ce
&&\frac{\dif}{\dif t}\mE|Y_t^{s-h,x,\mu,y}-Y_t^{s,x,\mu,y}|^2\\
&=&2\mE\<Y_t^{s-h,x,\mu,y}-Y_t^{s,x,\mu,y}, b_2(t,x,\mu,Y_t^{s-h,x,\mu,y})-b_2(t,x,\mu,Y_t^{s,x,\mu,y})\>\\
&&+\mE\|\s_2(t,x,\mu,Y_t^{s-h,x,\mu,y})-\s_2(t,x,\mu,Y_t^{s,x,\mu,y})\|^2\\
&\leq&-2\a(t)\mE|Y_t^{s-h,x,\mu,y}-Y_t^{s,x,\mu,y}|^2.
\de
Combining this with (\ref{Yl4}) yields that for any $h>0$, $t\geq s$, $x\in \mR^n$, $\mu\in \cP_2(\mR^n)$ and $y\in \mR^m$,
\ce
\mE|Y_t^{s-h,x,\mu,y}-Y_t^{s,x,\mu,y}|^2&\leq& \mE|Y_s^{s-h,x,\mu,y}-y|^2e^{-2\int_s^t\a(u)\dif u}\\
&\leq& C\left(1+|x|^2+\mu(|\cdot|^2)+|y|^2\right)e^{-2\int_s^t\a(u)\dif u}.
\de
Let $s\rightarrow-\infty$, it follows from (\ref{alpha}) that
\ce
\lim_{s\rightarrow-\infty}\mE|Y_t^{s-h,x,\mu,y}-Y_t^{s,x,\mu,y}|^2=0.
\de
Then for any $t\in\mR$, $x\in \mR^n$, $\mu\in \cP_2(\mR^n)$ and $y\in \mR^m$, $\{Y_t^{s,x,\mu,y}\}_{t\geq s}$ is a Cauchy sequence in $L^2(\Omega,\sF_t,\mP;\mR^m)$ as $s\rightarrow-\infty$. Thus there exists a random variable $\varsigma_t^{x,\mu,y}\in L^2(\Omega,\sF_t,\mP;\mR^m)$ such that for all $t\geq s$, $x\in \mR^n$, $\mu\in \cP_2(\mR^n)$ and $y\in \mR^m$,
\be
\mE|Y_t^{s,x,\mu,y}-\varsigma_t^{x,\mu,y}|^2\leq C\left(1+|x|^2+\mu(|\cdot|^2)+|y|^2\right)e^{-2\int_s^t\a(u)\dif u}.\label{yvsyl2}
\ee
We next show that $\varsigma_t^{x,\mu,y}$ is independent of $y$. Indeed, (\ref{Y1Y2l4}) implies that for any $y_1, y_2 \in \mR^m$,
\ce
\lim_{s\rightarrow-\infty}\mE|Y_t^{s,x,\mu,y_1}-Y_t^{s,x,\mu,y_2}|^2=0.
\de
Note that
\ce
\mE|\varsigma_t^{x,\mu,y_1}-\varsigma_t^{x,\mu,y_2}|^2&\leq& 3\mE|\varsigma_t^{x,\mu,y_1}-Y_t^{s,x,\mu,y_1}|^2+3\mE|Y_t^{s,x,\mu,y_1}-Y_t^{s,x,\mu,y_2}|^2\\
&&+3\mE|Y_t^{s,x,\mu,y_2}-\varsigma_t^{x,\mu,y_2}|^2.
\de
Taking $s\rightarrow-\infty$ and using (\ref{yvsyl2}), we get $\mE|\varsigma_t^{x,\mu,y_1}-\varsigma_t^{x,\mu,y_2}|^2=0$. So $\varsigma_t^{x,\mu,y_1}=\varsigma_t^{x,\mu,y_2}$ in $L^2(\Omega,\sF_t,\mP;\mR^m)$. We henceforth write $\varsigma_t^{x,\mu,y}$ by $\varsigma_t^{x,\mu}$.

Combining (\ref{yvsyl2}) and (\ref{Yl4}) yields
\ce
\mE|\varsigma_t^{x,\mu}|^2&\leq& 2\mE|\varsigma_t^{x,\mu}-Y_t^{s,x,\mu,y}|^2+2\mE|Y_t^{s,x,\mu,y}|^2\\
&\leq&C\left(1+|x|^2+\mu(|\cdot|^2)+|y|^2\right)e^{-2\g\int_s^t\a(u)\dif u}+C\left(1+|x|^2+\mu(|\cdot|^2)\right).
\de
Letting $s\rightarrow-\infty$, we deduce that
\ce
\mE|\varsigma_t^{x,\mu}|^2\leq C\left(1+|x|^2+\mu(|\cdot|^2)\right).
\de
Moreover, using (\ref{yvsyl2}) and letting $h\rightarrow+\infty$ in (\ref{yshxmuy}), we conclude that for any $t\geq s$, $x\in \mR^n$ and $\mu\in \cP_2(\mR^n)$,
\ce
\varsigma_t^{x,\mu}=\varsigma_s^{x,\mu}+\int_s^tb_2(r,x,\mu,\varsigma_r^{x,\mu})\dif r+\int_s^t\s_2(r,x,\mu,\varsigma_r^{x,\mu})\dif W_r.
\de
Hence, the proofs of (\ref{yvsl2}) and (\ref{vseq}) are complete.

Next, we prove (\ref{pyvs}). For any $h>0$, $t\geq s$, $x,\tx\in\mR^n$, $\mu\in\mathcal{P}_2(\mR^n)$, and $y\in\mR^m$, we have
\ce
&&\p_{\mu}Y_t^{s-h,x,\mu,y}(\tx)\\
&=&\p_{\mu}Y_s^{s-h,x,\mu,y}(\tx)+\int_s^t[\p_{\mu}b_2(r,x,\mu,Y_r^{s-h,x,\mu,y})(\tx)+\p_yb_2(r,x,\mu,Y_r^{s-h,x,\mu,y})\cdot\p_{\mu}Y_r^{s-h,x,\mu,y}(\tx)]\dif r\\
&&+\int_s^t[\p_{\mu}\s_2(r,x,\mu,Y_r^{s-h,x,\mu,y})(\tx)+\p_y\s_2(r,x,\mu,Y_r^{s-h,x,\mu,y})\cdot\p_{\mu}Y_r^{s-h,x,\mu,y}(\tx)]\dif W_r,
\de
and
\ce
\p_{\mu}Y_t^{s,x,\mu,y}(\tx)&=&\int_s^t[\p_{\mu}b_2(r,x,\mu,Y_r^{s,x,\mu,y})(\tx)+\p_yb_2(r,x,\mu,Y_r^{s,x,\mu,y})\cdot\p_{\mu}Y_r^{s,x,\mu,y}(\tx)]\dif r\\
&&+\int_s^t[\p_{\mu}\s_2(r,x,\mu,Y_r^{s,x,\mu,y})(\tx)+\p_y\s_2(r,x,\mu,Y_r^{s,x,\mu,y})\cdot\p_{\mu}Y_r^{s,x,\mu,y}(\tx)]\dif W_r.
\de
Following the argument in $(i)$ and using the Young inequality, we obtain
\ce
\mE\|\p_{\mu}Y_t^{s-h,x,\mu,y}(\tx)-\p_{\mu}Y_t^{s,x,\mu,y}(\tx)\|^2\leq C\left(1+|x|^2+\mu(|\cdot|^2)+|y|^2\right)e^{-2\g\int_s^t\a(u)\dif u}.
\de
Thus $\{\p_{\mu}Y_t^{s,x,\mu,y}(\tx)\}_{t\geq s}$ is a Cauchy sequence in $L^2(\Omega, \sF_t, \mP; \mR^{m\times n})$ as $s\rightarrow-\infty$. Consequently, there exists $\varpi_t^{x,\mu}(\tx)\in L^2(\Omega, \sF_t, \mP; \mR^{m\times n})$, independent of $y$ such that for all $t\geq s$, $x\in \mR^n$, $\mu\in \cP_2(\mR^n)$ and $y\in \mR^m$,
\be\label{pmuyl2}
\mE\|\p_{\mu}Y_t^{s,x,\mu,y}(\tx)-\varpi_t^{x,\mu}(\tx)\|^2\leq C\left(1+|x|^2+\mu(|\cdot|^2)+|y|^2\right)e^{-2\g\int_s^t\a(u)\dif u}.
\ee

Next, we show that $\varsigma_t^{x,\mu}$ is $L$-differentiable and that $\p_{\mu}\varsigma_t^{x,\mu}(\tx)=\varpi_t^{x,\mu}(\tx)$. Let $\vartheta$ be a random variable with law $\sL_\vartheta=\mu\in \cP_2(\mR^n)$. For all directions $\eta\in L^2(\Omega, \sF, \mP; \mR^n)$, $t\geq s$, $x\in \mR^n$, and $\lambda\in\mR\backslash\{0\}$, we have
\ce
&&\mE\left|\frac1{\lambda}(\varsigma_t^{x,\sL_{\vartheta+\lambda\eta}}-\varsigma_t^{x,\sL_{\vartheta}})-\tilde{E}[\varpi_t^{x,\sL_{\vartheta}}(\tilde{\vartheta})\cdot\tilde{\eta}]\right|^2\\
&\leq&\frac4{\lambda^2}\mE|\varsigma_t^{x,\sL_{\vartheta+\lambda\eta}}-Y_t^{s,x,\sL_{\vartheta+\lambda\eta},y}|^2
+4\mE\left|\frac1{\lambda}(Y_t^{s,x,\sL_{\vartheta+\lambda\eta},y}-Y_t^{s,x,\sL_{\vartheta},y})-\tilde{E}[\p_{\mu}Y_t^{s,x,\sL_{\vartheta},y}(\tilde{\vartheta})\cdot\tilde{\eta}]\right|^2\\
&&+\frac4{\lambda^2}\mE|Y_t^{s,x,\sL_{\vartheta},y}-\varsigma_t^{x,\sL_{\vartheta}}|^2
+4\mE\left|\tilde{E}[\p_{\mu}Y_t^{s,x,\sL_{\vartheta},y}(\tilde{\vartheta})\cdot\tilde{\eta}]-\tilde{E}[\varpi_t^{x,\sL_{\vartheta}}(\tilde{\vartheta})\cdot\tilde{\eta}]\right|^2,
\de
where $(\tilde{\vartheta},\tilde{\eta})$ is an independent copy of $(\vartheta,\eta)$, constructed on a probability space $(\tilde{\Omega},\tilde{\sF},\\\{\tilde{\sF}_t\}_{t \in[0, T]}, \tilde{\mP})$, which is an exact copy of the original probability space $(\Omega,\sF,\{\sF_t\}_{t \in[0, T]}, \mP)$, $\tilde{\mE}$ is the expectation taken with respect to $\tilde{\mP}$. By (\ref{yvsl2}), (\ref{pmuyl2}), and the $L$-differentiability of $Y_t^{s,x,\mu,y}$, it follows that
\ce
\lim_{\lambda\rightarrow0}\mE\left|\frac1{\lambda}(\varsigma_t^{x,\sL_{\vartheta+\lambda\eta}}-\varsigma_t^{x,\sL_{\vartheta}})-\tilde{E}[\varpi_t^{x,\sL_{\vartheta}}(\tilde{\vartheta})\cdot\tilde{\eta}]\right|^2=0.
\de
This shows that $\tilde{E}[\varpi_t^{x,\sL_{\vartheta}}(\tilde{\vartheta})\cdot\tilde{\eta}]$ is the directional derivative of $\varsigma_t^{x,\sL_{\vartheta}}$ in direction $\eta$. Following \cite{blpr},  $\varsigma_t^{x,\mu}$ is $L$-differentiable with $\p_{\mu}\varsigma_t^{x,\mu}(\tx)=\varpi_t^{x,\mu}(\tx)$. Thus
\ce
\mE\|\p_{\mu}Y_t^{s,x,\mu,y}(\tx)-\p_{\mu}\varsigma_t^{x,\mu}(\tx)\|^2\leq C\left(1+|x|^2+\mu(|\cdot|^2)+|y|^2\right)e^{-2\g\int_s^t\a(u)\dif u}.
\de
The remaining assertions in (\ref{pyvs}) follow from similar arguments and are therefore omitted. The proof is complete.
\end{proof}

Let $\{P_{s,t}^{x,\mu}\}_{t\geq s}$ denote the semigroup associated with the process $\{Y_t^{s,x,\mu,y}\}_{t\geq s}$, i.e., for any bounded measurable function $\psi:\mR^m\rightarrow\mR$,
\ce
P_{s,t}^{x,\mu}\psi(y)=\mathbb{E}\psi(Y_t^{s,x,\mu,y}),\quad y\in\mR^m.
\de
Let $\varsigma_t^{x,\mu}$ be given in Lemma \ref{yvs}. For each $t\in\mR$, we denote by $\nu_t^{x,\mu}$ the law of $\varsigma_t^{x,\mu}$.

\bl\label{nuxmu}
Suppose that $(\mathbf{H}^1_{b_{2}, \s_{2}})$ and $(\mathbf{H}^2_{b_{2}, \s_{2}})$ hold. Then $\{\nu_t^{x,\mu}\}_{t\in\mR}$ is an evolution system of measures for the semigroup $\{P_{s,t}^{x,\mu}\}_{t\geq s}$, i.e., for any $t\geq s$, $x\in\mR^n$, $\mu \in \cP_2(\mR^n)$, and $\psi\in C_b(\mR^m)$,
\be
\int_{\mR^m}P_{s,t}^{x,\mu}\psi(y)\nu_s^{x,\mu}(\dif y)=\int_{\mR^m}\psi(y)\nu_t^{x,\mu}(\dif y).\label{evo}
\ee
Moreover, there exists a constant $C>0$ such that for any $t\geq s$, $x\in\mR^n$, $\mu \in \cP_2(\mR^n)$, $y\in\mR^m$, and Lipschitz continuous function $\varphi$ on $\mR^m$,
\be
\left|P_{s,t}^{x,\mu}\varphi(y)-\int_{\mR^m}\varphi(z)\nu_t^{x,\mu}(\dif z)\right|\leq C Lip(\varphi)\left(1+|x|^2+\mu(|\cdot|^2)+|y|^2\right)^{1/2}e^{-\int_s^t\alpha(u)du},\label{evoesti}
\ee
where $Lip(\varphi):=\sup\limits_{x\neq y}|\varphi(x)-\varphi(y)|/|x-y|$. Furthermore, let $\{\tilde{\nu}_t^{x,\mu}\}_{t\in\mR}$ be another evolution system of measures for $\{P_{s,t}^{x,\mu}\}_{t\geq s}$. If for all $x\in\mR^n$, $\mu \in \cP_2(\mR^n)$,
\ce
\sup_{t\in\mR}\int_{\mR^m}|z|\tilde{\nu}_t^{x,\mu}(\dif z)<\infty,
\de
then $\nu_t^{x,\mu}=\tilde{\nu}_t^{x,\mu}$ for all $t\in\mR$, $x\in\mR^n$ and $\mu \in \cP_2(\mR^n)$.
\el

Since the proof of the above lemma is similar to that of \cite[Proposition 2.6]{swx}, we omit it.

\subsection{Nonautonomous Poisson equations}\label{nonpoieq}

In this subsection, we study the associated nonautonomous Poisson equation and establish the corresponding regularity estimates.

Let $g:\mR\times\mR^n \times \cP_2(\mR^n) \times \mR^m \rightarrow \mR^n$ satisfy the following centering condition:
\be
\int_{\mR^m}g(s,x,\mu,y)\nu_s^{x,\mu}(\dif y)=0,\quad s\in\mR, x\in\mR^n, \mu\in\cP_2(\mR^n).\label{cencon}
\ee
Assume further that for $i=0,1,2$ and $j=1,2,3$ with $1\leq i+j\leq3$, there exists a constant $C>0$ such that
\be\label{pxpypmug}
&&\sup_{t\in\mR,x\in\mR^n,\mu\in\cP_2(\mR^n),y\in\mR^m,\tx\in\mR^n}\max\Big\{\|\p_x^i\p_y^jg(t,x,\mu,y)\|,\|\p_y\p_{\mu}g(t,x,\mu,y)(\tx)\|,\no\\
&&\qquad\qquad\qquad\qquad\qquad\qquad\qquad\qquad\|\p_y\p_{\tx}\p_{\mu}g(t,x,\mu,y)(\tx)\|\Big\}< C.
\ee
For fixed $x\in\mR^n$ and $\mu\in\cP_2(\mR^n)$, we consider the following nonautonomous Poisson equation:
\be
\p_s\Psi(s,x,\mu,y)+\cL^{x,\mu}(s)\Psi(s,x,\mu,\cdot)(y)=-g(s,x,\mu,y),\quad s\in\mR, y\in\mR^m,\label{poieq}
\ee
where $\cL^{x,\mu}(s)$ denotes the generator associated with Eq.(\ref{frozeq}), that is,
\be
\cL^{x,\mu}(s)\psi(y):=b_2(s,x,\mu,y)\cdot\nabla\psi(y)+\frac12Tr[(\s_2\s_2^*)(s,x,\mu,y)\cdot\nabla^2\psi(y)], \quad \psi\in C^2(\mR^m).\label{generator}
\ee

\bp\label{poiprop}
Suppose that $(\mathbf{H}^1_{b_{2}, \s_{2}})$ and $(\mathbf{H}^2_{b_{2}, \s_{2}})$ hold. Assume that $g:\mR\times\mR^n \times \cP_2(\mR^n) \times \mR^m \rightarrow \mR^n$ satisfy conditions (\ref{cencon}) and (\ref{pxpypmug}). Define
\be
\Psi(s,x,\mu,y):=\int_s^{+\infty}\mE g(r,x,\mu,Y_r^{s,x,\mu,y})\dif r.\label{poisol}
\ee
Then $\Psi(s,x,\mu,y)$ is the unique solution to Eq.(\ref{poieq}) satisfying the centering condition (\ref{cencon}) and the following estimates:
\be\label{ppsi}
&&|\Psi(s,x,\mu,y)|\leq C\left(1+|x|^2+\mu(|\cdot|^2)+|y|^2\right)^{1/2}\Upsilon(s),\no\\
&&\|\p_x\Psi(s,x,\mu,y)\|+\|\p_x^2\Psi(s,x,\mu,y)\|\leq C\left(1+|x|^2+\mu(|\cdot|^2)+|y|^2\right)^{1/2}\Upsilon_{\g}(s),\no\\
&&\|\p_\mu\Psi(s,x,\mu,y)(\tx)\|+\|\p_{\tx}\p_\mu\Psi(s,x,\mu,y)(\tx)\|\leq C\left(1+|x|^2+\mu(|\cdot|^2)+|y|^2\right)^{1/2}\Upsilon_{\g}(s),\no\\
&&\|\p_y\Psi(s,x,\mu,y)\|\leq C\Upsilon(s),\quad \|\p_y^2\Psi(s,x,\mu,y)\|\leq C\Upsilon_{\g}(s),
\ee
where $\Upsilon_{\g}(s)$ is defined in (\ref{alpha}).
\ep
\begin{proof}
{\bf Step 1.} We show that $\Psi$ defined by (\ref{poisol}) is a solution to Eq.(\ref{poieq}) and satisfies the centering condition (\ref{cencon}).

For $r\ge s$, we have
\ce
P_{s,r}^{x,\mu}g(r,x,\mu,\cdot)(y)=\mathbb{E}g(r,x,\mu,Y_r^{s,x,\mu,y}),\quad y\in\mR^m.
\de
By the definition of $\Psi(s,x,\mu,y)$, it follows that
\ce
\p_s\Psi(s,x,\mu,y)&=&-g(s,x,\mu,y)+\int_s^{+\infty}\frac{\dif}{\dif s}P_{s,r}^{x,\mu}g(r,x,\mu,\cdot)(y)\dif r\\
&=&-g(s,x,\mu,y)+\int_s^{+\infty}-\cL^{x,\mu}(s)P_{s,r}^{x,\mu}g(r,x,\mu,\cdot)(y)\dif r\\
&=&-g(s,x,\mu,y)-\cL^{x,\mu}(s)\left(\int_s^{+\infty}P_{s,r}^{x,\mu}g(r,x,\mu,\cdot)(y)\dif r\right)\\
&=&-g(s,x,\mu,y)-\cL^{x,\mu}(s)\Psi(s,x,\mu,\cdot)(y),
\de
which shows that $\Psi(s,x,\mu,y)$ satisfies Eq.(\ref{poieq}).

By Fubini's theorem and the fact that $\{\nu_t^{x,\mu}\}_{t\in\mR}$ is an evolution system of measures for the semigroup $\{P_{s,t}^{x,\mu}\}_{t\geq s}$, we have
\ce
\int_{\mR^m}\Psi(s,x,\mu,y)\nu_s^{x,\mu}(\dif y)
&=&\int_{\mR^m}\int_s^{+\infty}\mE g(r,x,\mu,Y_r^{s,x,\mu,y})\dif r\nu_s^{x,\mu}(\dif y)\\
&=&\int_s^{+\infty}\int_{\mR^m}\mE g(r,x,\mu,Y_r^{s,x,\mu,y})\nu_s^{x,\mu}(\dif y)\dif r\\
&=&\int_s^{+\infty}\int_{\mR^m}P_{s,r}^{x,\mu}g(r,x,\mu,\cdot)(y)\nu_s^{x,\mu}(\dif y)\dif r\\
&=&\int_s^{+\infty}\int_{\mR^m}g(r,x,\mu,y)\nu_r^{x,\mu}(\dif y)\dif r=0.
\de

{\bf Step 2.} We estimate the derivative of $\Psi$ with respect to the measure $\mu$.

Observe that
\ce
&&\p_\mu\Psi(s,x,\mu,y)(\tx)\\
&=&\int_s^{+\infty}\p_\mu\mE[g(r,x,\mu,Y_r^{s,x,\mu,y})-g(r,x,\mu,\varsigma_r^{x,\mu})](\tx)\dif r\\
&=&\int_s^{+\infty}\mE[\p_\mu g(r,x,\mu,Y_r^{s,x,\mu,y})(\tx)-\p_\mu g(r,x,\mu,\varsigma_r^{x,\mu})(\tx)]\dif r\\
&&+\int_s^{+\infty}\mE[\p_y g(r,x,\mu,Y_r^{s,x,\mu,y})\cdot\p_\mu Y_r^{s,x,\mu,y}(\tx)-\p_yg(r,x,\mu,\varsigma_r^{x,\mu})\cdot\p_\mu\varsigma_r^{x,\mu}(\tx)]\dif r\\
&=&I_1+I_2.
\de
Applying (\ref{pxpypmug}) together with (\ref{yvsl2}), we obtain
\ce
\|I_1\|\leq C\int_s^{+\infty}\mE|Y_r^{s,x,\mu,y}-\varsigma_r^{x,\mu}|\dif r
\leq C\left(1+|x|^2+\mu(|\cdot|^2)+|y|^2\right)^{1/2}\Upsilon(s).
\de
Similarly, using (\ref{pxpypmug}), (\ref{yvsl2}), (\ref{parYl4}) and (\ref{pyvs}), we deduce that
\ce
\|I_2\|&\leq& \int_s^{+\infty}\mE\left\|\big(\p_y g(r,x,\mu,Y_r^{s,x,\mu,y})-\p_yg(r,x,\mu,\varsigma_r^{x,\mu})\big)\cdot\p_\mu Y_r^{s,x,\mu,y}(\tx)\right\|\dif r\\
&&+\int_s^{+\infty}\mE\left\|\p_yg(r,x,\mu,\varsigma_r^{x,\mu})\cdot\big(\p_\mu Y_r^{s,x,\mu,y}(\tx)-\p_\mu\varsigma_r^{x,\mu}(\tx)\big)\right\|\dif r\\
&\leq&C\int_s^{+\infty}\mE[|Y_r^{s,x,\mu,y}-\varsigma_r^{x,\mu}|\cdot\|\p_\mu Y_r^{s,x,\mu,y}(\tx)\|]\dif r\\
&&+C\int_s^{+\infty}\mE\|\p_\mu Y_r^{s,x,\mu,y}(\tx)-\p_\mu\varsigma_r^{x,\mu}(\tx)\|\dif r\\
&\leq&C\left(1+|x|^2+\mu(|\cdot|^2)+|y|^2\right)^{1/2}\Upsilon_{\g}(s).
\de
Combining the above estimates yields
\ce
\|\p_\mu\Psi(s,x,\mu,y)(\tx)\|\leq C\left(1+|x|^2+\mu(|\cdot|^2)+|y|^2\right)^{1/2}\Upsilon_{\g}(s).
\de
Using similar arguments, one can also derive
\ce
\|\p_{\tx}\p_\mu\Psi(s,x,\mu,y)(\tx)\|\leq C\left(1+|x|^2+\mu(|\cdot|^2)+|y|^2\right)^{1/2}\Upsilon_{\g}(s).
\de

The remaining estimates in \eqref{ppsi} follow from arguments analogous to those in \cite[Proposition 2.7]{swx}.

{\bf Step 3.} We prove the uniqueness of the solution to Eq.(\ref{poieq}) satisfying the centering condition and estimate \eqref{ppsi}.

Let $\tilde{\Psi}$ be another solution to Eq.(\ref{poieq}) satisfying the centering condition and estimate \eqref{ppsi}. Applying the It\^{o} formula to $\tilde{\Psi}(t,x,\mu,Y_t^{s,x,\mu,y})$ and taking expectation on both sides, we obtain
\be\label{tPsiito}
&&\mE\tilde{\Psi}(t,x,\mu,Y_t^{s,x,\mu,y})\no\\
&=&\tilde{\Psi}(s,x,\mu,y)+\mE\int_s^t[\p_r\tilde{\Psi}(r,x,\mu,Y_r^{s,x,\mu,y})+\cL^{x,\mu}(r)\tilde{\Psi}(r,x,\mu,\cdot)(Y_r^{s,x,\mu,y})]\dif r\no\\
&=&\tilde{\Psi}(s,x,\mu,y)-\int_s^t\mE g(r,x,\mu,Y_r^{s,x,\mu,y})\dif r.
\ee
Since $\tilde{\Psi}$ satisfies the centering condition and estimate \eqref{ppsi}, it follows that
\ce
|\mE\tilde{\Psi}(t,x,\mu,Y_t^{s,x,\mu,y})|
&=&\left|\mE\tilde{\Psi}(t,x,\mu,Y_t^{s,x,\mu,y})-\int_{\mR^m}\tilde{\Psi}(t,x,\mu,z)\nu_t^{x,\mu}(\dif z)\right|\\
&=&\left|P_{s,t}^{x,\mu}\tilde{\Psi}(t,x,\mu,\cdot)(y)-\int_{\mR^m}\tilde{\Psi}(t,x,\mu,z)\nu_t^{x,\mu}(\dif z)\right|\\
&\leq&C\|\tilde{\Psi}(t,x,\mu,\cdot)\|_{Lip}\left(1+|x|^2+\mu(|\cdot|^2)+|y|^2\right)^{1/2}e^{-\int_s^t\alpha(u)du}\\
&\leq&C\left(1+|x|^2+\mu(|\cdot|^2)+|y|^2\right)^{1/2}\Upsilon(t)e^{-\int_s^t\alpha(u)du}.
\de
Hence,
$$\lim_{t\rightarrow+\infty}|\mE\tilde{\Psi}(t,x,\mu,Y_t^{s,x,\mu,y})|=0.$$
Letting $t\rightarrow+\infty$ in (\ref{tPsiito}), we conclude that
$$
\tilde{\Psi}(s,x,\mu,y)=\int_s^{+\infty}\mE g(r,x,\mu,Y_r^{s,x,\mu,y})\dif r,
$$
which proves the uniqueness. The proof is complete.
\end{proof}

\section{Proofs of Theorem \ref{xbarx2} and \ref{xbarweak}}\label{geneproo}

In this section, we prove Theorem \ref{xbarx2} and \ref{xbarweak}.

\subsection{Some estimates for the system (\ref{orieq})}\label{estixy}
In this subsection, we establish uniform moment estimates for both the slow and fast components, which will be used throughout the section.

In order to fully utilize the evolution system of measures introduced in Subsection \ref{evomeas}, we extend the fast component $\{Y_t^\varepsilon\}_{t\geq0}$ of the system (\ref{orieq}) to the whole time axis $\mathbb{R}$. Let $\{\hat{W}_t\}_{t\geq0}$ be a $d_2$-dimensional standard Brownian motion independent of $\{W_t\}_{t\geq0}$, and define
\ce
\check{W}_t=
\begin{cases}
W_t, & t\geq0,\\
\hat{W}_{-t}, & t<0.
\end{cases}
\de
For any $x\in\mathbb{R}^n$, $\mu\in\cP_2(\mR^n)$ and $y\in\mathbb{R}^m$, define
\ce
\check{\a}(t)=
\begin{cases}
\alpha(t), & t\geq0,\\
\alpha(-t), & t<0,
\end{cases}
\qquad
\check{b}_2(t,x,\mu,y)=
\begin{cases}
b_2(t,x,\mu,y), & t\geq0,\\
b_2(-t,x,\mu,y), & t<0,
\end{cases}
\de
and
\ce
\check{\s}_2(t,x,\mu,y)=
\begin{cases}
\s_2(t,x,\mu,y), & t\geq0,\\
\s_2(-t,x,\mu,y), & t<0.
\end{cases}
\de
If $b_2,\s_2$ satisfy $(\mathbf{H}^3_{b_{2}, \s_{2}})$, then $\check{b}_2$ and $\check{\s}_2$ satisfy $(\mathbf{H}^1_{b_{2}, \s_{2}})$ and $(\mathbf{H}^2_{b_{2}, \s_{2}})$ on $\mathbb{R}$ with the function $\check{\a}$ in place of $\alpha$. As a consequence, the corresponding time-inhomogeneous frozen SDE
\ce
\dif Y_t=\check{b}_2(t,x,\mu,Y_t)\dif t+\check{\s}_2(t,x,\mu,Y_t)\dif \check{W}_t, \quad Y_s=y,
\de
admits a unique strong solution, denoted by $\{Y_t^{s,x,\mu,y}\}_{t\geq s}$. We denote by $\{P_{s,t}^{x,\mu}\}_{t\geq s}$ the transition semigroup associated with the process $\{Y_t^{s,x,\mu,y}\}_{t\geq s}$. Throughout the remainder of the paper, $\{\nu_t^{x,\mu}\}_{t\in\mathbb{R}}$ stands for the evolution system of measures associated with $\{P_{s,t}^{x,\mu}\}_{t\geq s}$.

\bl\label{xeye}
Suppose that $(\mathbf{H}^1_{b_{1}, \s_{1}})$ and $(\mathbf{H}^3_{b_{2}, \s_{2}})$ hold. Then, for any $T>0$, there exists $C_T>0$ such that
\be
&&\mE(\sup_{0\leq t\leq T}|X_t^\eps|^4)\leq C_T(1+\mE|\xi|^4+\mE|\varrho|^4),\label{xel4}\\
&&\sup_{0\leq t\leq T}\mE|Y_t^\eps|^4\leq C_T(1+\mE|\xi|^4+\mE|\varrho|^4),\label{yel4out}\\
&&\mE(\sup_{0\leq t\leq T}|Y_t^\eps|^4)\leq 2\mE|\varrho|^4+C_T(1+\mE|\xi|^4+\mE|\varrho|^4)\int_0^{T/\eps}\a(r)\dif r.\label{yel4in}
\ee
\el
\begin{proof}
Since the first two estimates follow from \cite[Lemma 3.3]{swx}, we only establish the third estimate.

For $Y_t^\eps$, applying the It\^{o} formula to $|Y_t^\eps|^4$ yields
\ce
|Y_t^\eps|^4
&=&|\varrho|^4+\frac{4}{\eps}\int_0^t|Y_r^\eps|^2\<Y_r^\eps, b_2(r/\eps,X_r^{\eps}, \sL_{X_r^{\eps}}, Y_r^{\eps})\>\dif r\\
&&+\frac{4}{\sqrt\eps}\int_0^t|Y_r^\eps|^2\<Y_r^\eps, \s_2(r/\eps,X_r^{\eps}, \sL_{X_r^{\eps}}, Y_r^{\eps})\dif W_r\>\\
&&+\frac{4}{\eps}\int_0^t\|\s_2(r/\eps,X_r^{\eps}, \sL_{X_r^{\eps}}, Y_r^{\eps})Y_r^\eps\|^2\dif r\\
&&+\frac{2}{\eps}\int_0^t|Y_r^\eps|^2\|\s_2(r/\eps,X_r^{\eps}, \sL_{X_r^{\eps}}, Y_r^{\eps})\|^2\dif r.
\de
Using (\ref{b2s2dis}), we obtain
\ce
|Y_t^\eps|^4
&\leq&|\varrho|^4+\frac{2}{\eps}\int_0^t|Y_r^\eps|^2\left[2\<Y_r^\eps, b_2(r/\eps,X_r^{\eps}, \sL_{X_r^{\eps}}, Y_r^{\eps})\>+3\|\s_2(r/\eps,X_r^{\eps}, \sL_{X_r^{\eps}}, Y_r^{\eps})\|^2\right]\dif r\\
&&+\frac{4}{\sqrt\eps}\int_0^t|Y_r^\eps|^2\<Y_r^\eps, \s_2(r/\eps,X_r^{\eps}, \sL_{X_r^{\eps}}, Y_r^{\eps})\dif W_r\>\\
&\leq&|\varrho|^4+\frac{2}{\eps}\int_0^t|Y_r^\eps|^2\left[-2\g\a(r/\eps)|Y_r^\eps|^2+C\a(r/\eps)(1+|X_r^\eps|^2+\mE|X_r^\eps|^2)\right]\dif r\\
&&+\frac{4}{\sqrt\eps}\int_0^t|Y_r^\eps|^2\<Y_r^\eps, \s_2(r/\eps,X_r^{\eps}, \sL_{X_r^{\eps}}, Y_r^{\eps})\dif W_r\>\\
&\leq&|\varrho|^4+\frac{C}{\eps}\int_0^t\a(r/\eps)(1+|X_r^\eps|^4+\mE|X_r^\eps|^4)\dif r\\
&&+\frac{4}{\sqrt\eps}\int_0^t|Y_r^\eps|^2\<Y_r^\eps, \s_2(r/\eps,X_r^{\eps}, \sL_{X_r^{\eps}}, Y_r^{\eps})\dif W_r\>.
\de
By the Burkholder-Davis-Gundy inequality and the Young inequality, together with (\ref{xel4}), we have
\ce
\mE(\sup_{0\leq t\leq T}|Y_t^\eps|^4)&\leq&\mE|\varrho|^4+\frac{C_T(1+\mE|\xi|^4+\mE|\varrho|^4)}{\eps}\int_0^T\a(r/\eps)\dif r\\
&&+\frac{C}{\sqrt\eps}\mE\left[\int_0^T|Y_r^\eps|^6\|\s_2(r/\eps,X_r^{\eps}, \sL_{X_r^{\eps}}, Y_r^{\eps})\|^2\dif r\right]^{\frac12}\\
&\leq&\mE|\varrho|^4+\frac{C_T(1+\mE|\xi|^4+\mE|\varrho|^4)}{\eps}\int_0^T\a(r/\eps)\dif r\\
&&+\frac{C}{\sqrt\eps}\mE\left[\sup_{0\leq t\leq T}|Y_t^\eps|^4\int_0^T\a(r/\eps)|Y_r^\eps|^2(1+|X_r^\eps|^2+\mE|X_r^\eps|^2+|Y_r^\eps|^2)\dif r\right]^{\frac12}\\
&\leq&\mE|\varrho|^4+\frac{C_T(1+\mE|\xi|^4+\mE|\varrho|^4)}{\eps}\int_0^T\a(r/\eps)\dif r+\frac12\mE(\sup_{0\leq t\leq T}|Y_t^\eps|^4)\\
&&+\frac{C}{\eps}\mE\int_0^T\a(r/\eps)|Y_r^\eps|^2(1+|X_r^\eps|^2+\mE|X_r^\eps|^2+|Y_r^\eps|^2)\dif r\\
&\leq&\mE|\varrho|^4+\frac{C_T(1+\mE|\xi|^4+\mE|\varrho|^4)}{\eps}\int_0^T\a(r/\eps)\dif r+\frac12\mE(\sup_{0\leq t\leq T}|Y_t^\eps|^4),
\de
which implies that
\ce
\mE(\sup_{0\leq t\leq T}|Y_t^\eps|^4)&\leq&2\mE|\varrho|^4+\frac{C_T(1+\mE|\xi|^4+\mE|\varrho|^4)}{\eps}\int_0^T\a(r/\eps)\dif r\\
&\leq&2\mE|\varrho|^4+C_T(1+\mE|\xi|^4+\mE|\varrho|^4)\int_0^{T/\eps}\a(r)\dif r.
\de
The proof is complete.
\end{proof}

\subsection{Proof of Theorem \ref{xbarx2}}\label{stravegen}
In this subsection, we require that $\s_1$ is independent of $y$ and prove Theorem \ref{xbarx2}.

\bl
Suppose that $(\mathbf{H}^1_{b_{1}, \s_{1}})$ and $(\mathbf{H}^3_{b_{2}, \s_{2}})$ hold. Then for any $\eps>0$, Eq.(\ref{straveeq}) has a unique solution, denoted by $\{\bar{X}_t^\eps\}_{t\geq0}$. Moreover, for $T>0$, there exists $C_T>0$ such that
\be
\mE(\sup_{0\leq t\leq T}|\bar{X}_t^\eps|^2)\leq C_T(1+\mE|\xi|^2).\label{barxel2}
\ee
\el

\begin{proof}
By Lemma \ref{nuxmu}, (\ref{Y1Y2l4}) and (\ref{b1s1lip}), we have
\ce
&&|\bar{b}_1(t,x_1,\mu_1)-\bar{b}_1(t,x_2,\mu_2)|\\
&=&\left|\int_{\mR^m}b_1(x_1,\mu_1,y)\nu_t^{x_1,\mu_1}(\dif y)-\int_{\mR^m}b_1(x_2,\mu_2,y)\nu_t^{x_2,\mu_2}(\dif y)\right|\\
&\leq&\left|\int_{\mR^m}b_1(x_1,\mu_1,y)\nu_t^{x_1,\mu_1}(\dif y)-\mE b_1(x_1,\mu_1,Y_t^{s,x_1,\mu_1,0})\right|\\
&&+\left|\mE b_1(x_1,\mu_1,Y_t^{s,x_1,\mu_1,0})-\mE b_1(x_2,\mu_2,Y_t^{s,x_2,\mu_2,0})\right|\\
&&+\left|\mE b_1(x_2,\mu_2,Y_t^{s,x_2,\mu_2,0})-\int_{\mR^m}b_1(x_2,\mu_2,y)\nu_t^{x_2,\mu_2}(\dif y)\right|\\
&\leq&C\left(1+|x_1|^2+\mu_1(|\cdot|^2)+|x_2|^2+\mu_2(|\cdot|^2)\right)^{1/2}e^{-\int_s^t\alpha(u)du}\\
&&+C\left(|x_1-x_2|+\mW_2(\mu_1, \mu_2)+\mE|Y_t^{s,x_1,\mu_1,0}-Y_t^{s,x_2,\mu_2,0}|\right)\\
&\leq&C\left(1+|x_1|^2+\mu_1(|\cdot|^2)+|x_2|^2+\mu_2(|\cdot|^2)\right)^{1/2}e^{-\int_s^t\alpha(u)du}\\
&&+C\left(|x_1-x_2|+\mW_2(\mu_1, \mu_2)\right).
\de
Letting $s\rightarrow-\infty$ in above inequality, and using (\ref{alpha}), we have
\be
|\bar{b}_1(t,x_1,\mu_1)-\bar{b}_1(t,x_2,\mu_2)|\leq C\left(|x_1-x_2|+\mW_2(\mu_1, \mu_2)\right).\label{barb1lip}
\ee
On the other hand, by (\ref{yvsl2}), we get
\be
|\bar{b}_1(t,x,\mu)|^2&\leq&\mE|b_1(x,\mu,\varsigma_t^{x,\mu})|^2\no\\
&\leq&C\left(1+|x|^2+\mu(|\cdot|^2)+\mE|\varsigma_t^{x,\mu}|^2\right)\no\\
&\leq&C\left(1+|x|^2+\mu(|\cdot|^2)\right).\label{barb1grow}
\ee
Combining (\ref{barb1lip}) and (\ref{barb1grow}) with (\ref{b1s1lip}), we can obtain that Eq.(\ref{straveeq}) has a unique solution, and (\ref{barxel2}) holds. The proof is complete.
\end{proof}

{\bf Proof of Theorem \ref{xbarx2}.}
Note that
\ce
X_t^{\eps}-\bar{X}_t^\eps&=&\int_0^t[b_1(X_s^{\eps}, \sL_{X_s^{\eps}}, Y_s^{\eps})-\bar{b}_1(s/\eps,\bar{X}_s^\eps,\sL_{\bar{X}_s^\eps})]\dif s\\
&&+\int_0^t[\s_1(X_s^{\eps}, \sL_{X_s^{\eps}})-\s_1(\bar{X}_s^\eps,\sL_{\bar{X}_s^\eps})]\dif B_s.
\de
From the H\"{o}lder inequality, the Burkholder-Davis-Gundy inequality, (\ref{barb1lip}) and (\ref{b1s1lip}), it follows that for any $t\in[0,T]$,
\ce
\mE(\sup_{0\leq t\leq r}|X_t^{\eps}-\bar{X}_t^\eps|^2)&\leq&3\mE\left(\sup_{0\leq t\leq r}\left|\int_0^t[b_1(X_s^{\eps}, \sL_{X_s^{\eps}}, Y_s^{\eps})-\bar{b}_1(s/\eps,X_s^{\eps}, \sL_{X_s^{\eps}})]\dif s\right|^2\right)\\
&&+3\mE\left(\sup_{0\leq t\leq r}\left|\int_0^t[\bar{b}_1(s/\eps,X_s^{\eps}, \sL_{X_s^{\eps}})-\bar{b}_1(s/\eps,\bar{X}_s^\eps,\sL_{\bar{X}_s^\eps})]\dif s\right|^2\right)\\
&&+3\mE\left(\sup_{0\leq t\leq r}\left|\int_0^t[\s_1(X_s^{\eps}, \sL_{X_s^{\eps}})-\s_1(\bar{X}_s^\eps,\sL_{\bar{X}_s^\eps})]\dif B_s\right|^2\right)\\
&\leq&3\mE\left(\sup_{0\leq t\leq r}\left|\int_0^t[b_1(X_s^{\eps}, \sL_{X_s^{\eps}}, Y_s^{\eps})-\bar{b}_1(s/\eps,X_s^{\eps}, \sL_{X_s^{\eps}})]\dif s\right|^2\right)\\
&&+C_T\int_0^r\mE|X_s^{\eps}-\bar{X}_s^\eps|^2\dif s.
\de
The Gronwall inequality yields that
\ce
\mE(\sup_{0\leq t\leq T}|X_t^{\eps}-\bar{X}_t^\eps|^2)\leq C_T\mE\left(\sup_{0\leq t\leq T}\left|\int_0^t[b_1(X_s^{\eps}, \sL_{X_s^{\eps}}, Y_s^{\eps})-\bar{b}_1(s/\eps,X_s^{\eps}, \sL_{X_s^{\eps}})]\dif s\right|^2\right).
\de

In what follows, we aim to estimate the right-hand side of the above inequality. According to Proposition \ref{poiprop}, there exists a solution $\Psi(s,x,\mu,y)$ such that
\ce
\p_s\Psi(s,x,\mu,y)+\cL^{x,\mu}(s)\Psi(s,x,\mu,\cdot)(y)=-[b_1(x,\mu,y)-\bar{b}_1(s,x,\mu)], \label{b1barb1poieq}
\de
where $\cL^{x,\mu}(s)$ is defined by (\ref{generator}). Applying the It\^{o} formula to $\Psi(t/\eps,X_t^{\eps}, \sL_{X_t^{\eps}}, Y_t^{\eps})$, we have
\ce
&&\Psi(t/\eps,X_t^{\eps}, \sL_{X_t^{\eps}}, Y_t^{\eps})\\
&=&\Psi(0,\xi, \sL_{\xi}, \varrho)+\eps^{-1}\int_0^t\p_s\Psi(s/\eps,X_s^{\eps}, \sL_{X_s^{\eps}}, Y_s^{\eps})\dif s\\
&&+\int_0^tb_1(X_s^{\eps}, \sL_{X_s^{\eps}}, Y_s^{\eps})\cdot\p_x\Psi(s/\eps,X_s^{\eps}, \sL_{X_s^{\eps}}, Y_s^{\eps})\dif s\\
&&+\frac12\int_0^tTr[(\s_1\s_1^*)(X_s^{\eps}, \sL_{X_s^{\eps}})\cdot\p_x^2\Psi(s/\eps,X_s^{\eps}, \sL_{X_s^{\eps}}, Y_s^{\eps})]\dif s\\
&&+\int_0^t\p_x\Psi(s/\eps,X_s^{\eps}, \sL_{X_s^{\eps}}, Y_s^{\eps})\cdot\s_1(X_s^{\eps}, \sL_{X_s^{\eps}})\dif B_s\\
&&+\int_0^t\tilde{\mE}[b_1(\tilde{X}_s^{\eps}, \sL_{X_s^{\eps}}, \tilde{Y}_s^{\eps})\cdot\p_\mu\Psi(s/\eps,X_s^{\eps}, \sL_{X_s^{\eps}}, Y_s^{\eps})(\tilde{X}_s^{\eps})]\dif s\\
&&+\frac12\int_0^t\tilde{\mE}\left[Tr[(\s_1\s_1^*)(\tilde{X}_s^{\eps}, \sL_{X_s^{\eps}})\cdot\p_{\tx}\p_\mu\Psi(s/\eps,X_s^{\eps}, \sL_{X_s^{\eps}}, Y_s^{\eps})(\tilde{X}_s^{\eps})]\right]\dif s\\
&&+\eps^{-1}\int_0^tb_2(s/\eps,X_s^{\eps}, \sL_{X_s^{\eps}}, Y_s^{\eps})\cdot\p_y\Psi(s/\eps,X_s^{\eps}, \sL_{X_s^{\eps}}, Y_s^{\eps})\dif s\\
&&+\frac1{2\eps}\int_0^tTr[(\s_2\s_2^*)(s/\eps,X_s^{\eps}, \sL_{X_s^{\eps}}, Y_s^{\eps})\cdot\p_y^2\Psi(s/\eps,X_s^{\eps}, \sL_{X_s^{\eps}}, Y_s^{\eps})]\dif s\\
&&+\frac1{\sqrt\eps}\int_0^t\p_y\Psi(s/\eps,X_s^{\eps}, \sL_{X_s^{\eps}}, Y_s^{\eps})\cdot\s_2(s/\eps,X_s^{\eps}, \sL_{X_s^{\eps}}, Y_s^{\eps})\dif W_s\\
&=&\Psi(0,\xi, \sL_{\xi}, \varrho)+\eps^{-1}\int_0^t\p_s\Psi(s/\eps,X_s^{\eps}, \sL_{X_s^{\eps}}, Y_s^{\eps})\dif s\\
&&+\int_0^t\cL_1^{\sL_{X_s^{\eps}},Y_s^{\eps}}\Psi(s/\eps,\cdot, \sL_{X_s^{\eps}},Y_s^{\eps}) (X_s^{\eps})\dif s\\
&&+\int_0^t\tilde{\mE}[b_1(\tilde{X}_s^{\eps}, \sL_{X_s^{\eps}}, \tilde{Y}_s^{\eps})\cdot\p_\mu\Psi(s/\eps,X_s^{\eps}, \sL_{X_s^{\eps}}, Y_s^{\eps})(\tilde{X}_s^{\eps})]\dif s\\
&&+\frac12\int_0^t\tilde{\mE}\left[Tr\big[(\s_1\s_1^*)(\tilde{X}_s^{\eps}, \sL_{X_s^{\eps}})\cdot\p_{\tx}\p_\mu\Psi(s/\eps,X_s^{\eps}, \sL_{X_s^{\eps}}, Y_s^{\eps})(\tilde{X}_s^{\eps})\big]\right]\dif s\\
&&+\eps^{-1}\int_0^t\cL^{X_s^{\eps}, \sL_{X_s^{\eps}}}(s/\eps)\Psi(s/\eps,X_s^{\eps}, \sL_{X_s^{\eps}},\cdot) (Y_s^{\eps})\dif s\\
&&+N_t^{\eps,1}+\frac1{\sqrt\eps}N_t^{\eps,2},
\de
where the process $(\tilde{X}_s^{\varepsilon},\tilde{Y}_s^{\varepsilon})$ represents a copy of the original process $(X_s^{\varepsilon},Y_s^{\varepsilon})$, constructed on a probability space $(\tilde{\Omega},\tilde{\mathscr{F}},\{\tilde{\mathscr{F}}_t\}_{t \in[0, T]}, \tilde{\mathbb{P}})$,
which is an exact copy of the original probability space $(\Omega,\mathscr{F},\{\mathscr{F}_t\}_{t \in[0, T]}, \mathbb{P})$, $\tilde{\mE}$ is the expectation taken with respect to $\tilde{\mathbb{P}}$, and
\ce
\cL_1^{\mu,y}\psi(x)&:=&b_1(x,\mu,y)\cdot\nabla\psi(x)+\frac12Tr[(\s_1\s_1^*)(x,\mu)\cdot\nabla^2\psi(x)],\\
N_t^{\eps,1}&:=&\int_0^t\p_x\Psi(s/\eps,X_s^{\eps}, \sL_{X_s^{\eps}}, Y_s^{\eps})\cdot\s_1(X_s^{\eps}, \sL_{X_s^{\eps}})\dif B_s,\\
N_t^{\eps,2}&:=&\int_0^t\p_y\Psi(s/\eps,X_s^{\eps}, \sL_{X_s^{\eps}}, Y_s^{\eps})\cdot\s_2(s/\eps,X_s^{\eps}, \sL_{X_s^{\eps}}, Y_s^{\eps})\dif W_s.
\de
Then, we conclude that
\ce
&&\mE(\sup_{0\leq t\leq T}|X_t^{\eps}-\bar{X}_t^\eps|^2)\\
&\leq&C_T\mE\left(\sup_{0\leq t\leq T}\left|\int_0^t\left[\p_s\Psi(s/\eps,X_s^{\eps}, \sL_{X_s^{\eps}}, Y_s^{\eps})+\cL^{X_s^{\eps}, \sL_{X_s^{\eps}}}(s/\eps)\Psi(s/\eps,X_s^{\eps}, \sL_{X_s^{\eps}},\cdot) (Y_s^{\eps})\right]\dif s\right|^2\right)\\
&\leq&C_T\eps^2\mE\left(\sup_{0\leq t\leq T}\left|\Psi(t/\eps,X_t^{\eps}, \sL_{X_t^{\eps}}, Y_t^{\eps})\right|^2\right)
+C_T\eps^2\mE\left|\Psi(0,\xi, \sL_{\xi}, \varrho)\right|^2\\
&&+C_T\eps^2\mE\left(\sup_{0\leq t\leq T}\left|\int_0^t\cL_1^{\sL_{X_s^{\eps}},Y_s^{\eps}}\Psi(s/\eps,\cdot, \sL_{X_s^{\eps}},Y_s^{\eps}) (X_s^{\eps})\dif s\right|^2\right)\\
&&+C_T\eps^2\mE\left(\sup_{0\leq t\leq T}\left|\int_0^t\tilde{\mE}[b_1(\tilde{X}_s^{\eps}, \sL_{X_s^{\eps}}, \tilde{Y}_s^{\eps})\cdot\p_\mu\Psi(s/\eps,X_s^{\eps}, \sL_{X_s^{\eps}}, Y_s^{\eps})(\tilde{X}_s^{\eps})]\dif s\right|^2\right)\\
&&+C_T\eps^2\mE\left(\sup_{0\leq t\leq T}\left|\int_0^t\tilde{\mE}\left[Tr\big[(\s_1\s_1^*)(\tilde{X}_s^{\eps}, \sL_{X_s^{\eps}})\cdot\p_{\tx}\p_\mu\Psi(s/\eps,X_s^{\eps}, \sL_{X_s^{\eps}}, Y_s^{\eps})(\tilde{X}_s^{\eps})\big]\right]\dif s\right|^2\right)\\
&&+C_T\eps^2\mE\left(\sup_{0\leq t\leq T}\left|N_t^{\eps,1}\right|^2\right)+C_T\eps\mE\left(\sup_{0\leq t\leq T}\left|N_t^{\eps,2}\right|^2\right)\\
&=&\sum_{i=1}^7\cI_i.
\de
For $\cI_1$ and $\cI_2$, by Proposition \ref{poiprop} and Lemma \ref{xeye}, we have
\ce
\cI_1&\leq&C_T\eps^2\left\{1+\mE\left(\sup_{0\leq t\leq T}\left|X_t^{\eps}\right|^2\right)+\mE\left(\sup_{0\leq t\leq T}
\left|Y_t^{\eps}\right|^2\right)\right\}\sup_{0\leq t\leq T}\Upsilon^2(t/\eps)\\
&\leq&C_T\eps^2(1+\mE|\xi|^4+\mE|\varrho|^4)\sup_{0\leq t\leq T}\Upsilon^2(t/\eps)\left[1+\left(\int_0^{T/\eps}\a(s)\dif s\right)^{1/2}\right],
\de
and
\ce
\cI_2\leq C_T\eps^2(1+\mE|\xi|^2+\mE|\varrho|^2)\Upsilon^2(0).
\de
For $\cI_3$, $\cI_4$ and $\cI_5$, based on $(\mathbf{H}^1_{b_{1}, \s_{1}})$, the H\"{o}lder inequality, and Proposition \ref{poiprop}, one can
obtain that
\ce
\cI_3&\leq&C_T\eps^2\mE\int_0^T\left|\cL_1^{\sL_{X_s^{\eps}},Y_s^{\eps}}\Psi(s/\eps,\cdot, \sL_{X_s^{\eps}},Y_s^{\eps}) (X_s^{\eps})\right|^2\dif s\\
&\leq&C_T\eps^2\left\{1+\mE\left(\sup_{0\leq t\leq T}|X_t^{\eps}|^4\right)+\sup_{0\leq t\leq T}\mE|Y_t^{\eps}|^4\right\}\int_0^T\Upsilon^2_{\g}(s/\eps)\dif s,
\de
and
\ce
\cI_4+\cI_5
&\leq&C_T\eps^2\mE\tilde{\mE}\int_0^T\left|b_1(\tilde{X}_s^{\eps},\sL_{X_s^{\eps}},\tilde{Y}_s^{\eps})\cdot\p_\mu\Psi(s/\eps,X_s^{\eps},\sL_{X_s^{\eps}},Y_s^{\eps})(\tilde{X}_s^{\eps})\right|^2\dif s\\
&&+C_T\eps^2\mE\tilde{\mE}\int_0^T\left|Tr[(\s_1\s_1^*)(\tilde{X}_s^{\eps},\sL_{X_s^{\eps}})\cdot\p_{\tx}\p_\mu\Psi(s/\eps,X_s^{\eps},\sL_{X_s^{\eps}},Y_s^{\eps})(\tilde{X}_s^{\eps})]\right|^2\dif s\\
&\leq&C_T\eps^2\left\{1+\mE\left(\sup_{0\leq t\leq T}|X_t^{\eps}|^4\right)+\sup_{0\leq t\leq T}\mE|Y_t^{\eps}|^4\right\}\int_0^T\Upsilon^2_{\g}(s/\eps)\dif s.
\de
For $\cI_6$ and $\cI_7$, by the Burkholder-Davis-Gundy inequality, the H\"{o}lder inequality and
Proposition \ref{poiprop}, it holds that
\ce
\cI_6+\cI_7&\leq&C_T\eps^2\mE\int_0^T\left|\p_x\Psi(s/\eps,X_s^{\eps}, \sL_{X_s^{\eps}}, Y_s^{\eps})\cdot\s_1(X_s^{\eps}, \sL_{X_s^{\eps}})\right|^2\dif s\\
&&+C_T\eps\mE\int_0^T\left|\p_y\Psi(s/\eps,X_s^{\eps}, \sL_{X_s^{\eps}}, Y_s^{\eps})\cdot\s_2(s/\eps,X_s^{\eps}, \sL_{X_s^{\eps}}, Y_s^{\eps})\right|^2\dif s\\
&\leq&C_T\left(\eps^2\int_0^T\Upsilon^2_{\g}(s/\eps)\dif s+\eps\int_0^T\a(s/\eps)\Upsilon^2(s/\eps)\dif s \right)\\
&&\qquad\qquad\cdot\left\{1+\mE\left(\sup_{0\leq t\leq T}|X_t^{\eps}|^2\right)+\sup_{0\leq t\leq T}\mE|Y_t^{\eps}|^2\right\}
\de
Collecting the above deductions and applying Lemma \ref{xeye}, we conclude that
\ce
&&\mE(\sup_{0\leq t\leq T}|X_t^{\eps}-\bar{X}_t^\eps|^2)\\
&\leq&C_T(1+\mE|\xi|^4+\mE|\varrho|^4)\eps^2\sup_{0\leq t\leq T}\Upsilon^2(t/\eps)\left[1+\left(\int_0^{T/\eps}\a(s)\dif s\right)^{1/2}\right]\\
&&+C_T(1+\mE|\xi|^4+\mE|\varrho|^4)\left(\eps^2\int_0^T\Upsilon^2_{\g}(s/\eps)\dif s+\eps\int_0^T\a(s/\eps)\Upsilon^2(s/\eps)\dif s \right)\\
&\leq&C_T(1+\mE|\xi|^4+\mE|\varrho|^4)\eps^2\\
&&\cdot\left[\sup_{0\leq t\leq T}\Upsilon^2_{\g}(t/\eps)+\sup_{0\leq t\leq T}\Upsilon^2(t/\eps)\left(\int_0^{T/\eps}\a(s)\dif s\right)^{1/2}+\int_0^{T/\eps}\a(s)\Upsilon^2(s)\dif s \right].
\de
The proof is complete.

\subsection{Proof of Theorem \ref{xbarweak}}\label{weakavegen}
In this subsection, we prove Theorem \ref{xbarweak}.

It follows from Lemma \ref{yvs} that, under $(\mathbf{H}^3_{b_{2}, \s_{2}})$ and $(\mathbf{H}^4_{b_{2}, \s_{2}})$, the following estimates hold:
\ce
&&\sup_{t\geq0}\sum_{i=1}^4\mE\|\p_x^i\varsigma_t^{x,\mu}\|^4\leq C, \quad \sup_{t\geq0}\sum_{i=0}^3\mE\|\p_{\tx}^i\p_{\mu}\varsigma_t^{x,\mu}(\tx)\|^2\leq C,\\
&&\sup_{t\geq0}\sum_{i=1}^2\mE\|\p_x^i\p_{\mu}\varsigma_t^{x,\mu}(\tx)\|^2\leq C, \quad \sup_{t\geq0}\sum_{i=1}^2\mE\|\p_x^i\p_{\tx}\p_{\mu}\varsigma_t^{x,\mu}(\tx)\|^2\leq C,\\
&&\sup_{t\geq0}\sum_{i,j=0}^1\mE\|\p_{\cx}^i\p_{\tx}^j\p_{\mu}^2\varsigma_t^{x,\mu}(\tx,\cx)\|^2\leq C.
\de
Since $\bar{b}_1(t,x,\mu)=\mE b_1(x,\mu,\varsigma_t^{x,\mu})$ and $\overline{\s_1\s_1^*}(t,x,\mu)=\mE (\s_1\s_1^*)(x,\mu,\varsigma_t^{x,\mu})$, it follows from $(\mathbf{H}^2_{b_{1}, \s_{1}})$, $(\mathbf{H}^3_{b_{2}, \s_{2}})$ and $(\mathbf{H}^4_{b_{2}, \s_{2}})$ that for any $t\geq0$, $\bar{b}_1(t,\cdot,\cdot)\in (C_b^{4,(1,3)}\cap C_b^{2,(2,2)})(\mR^n\times\cP_2(\mR^n),\mR^n)$ and $\overline{\s_1\s_1^*}(t,\cdot,\cdot)\in (C_b^{4,(1,3)}\cap C_b^{2,(2,2)})(\mR^n\times\cP_2(\mR^n),\mR^{n\times n})$. Moreover, by $(\mathbf{H}^3_{b_{1}, \s_{1}})$, $\overline{\s_1\s_1^*}$ is non-degenerate, i.e.,
\ce
\inf_{t\geq0,x\in\mR^n,\mu\in\cP_2(\mR^n),z\in\mR^n\backslash\{0\}}\frac{\<(\overline{\s_1\s_1^*})(t,x,\mu)\cdot z,z\>}{|z|^2}>0.
\de
Since $\overline{\s_1\s_1^*}(t,\cdot,\cdot)\in (C_b^{4,(1,3)}\cap C_b^{2,(2,2)})(\mR^n\times\cP_2(\mR^n),\mR^{n\times n})$, it follows from \cite[Theorem 1.1]{dn} that
$\bar{\s}_1(t,\cdot,\cdot)\in (C_b^{4,(1,3)}\cap C_b^{2,(2,2)})(\mR^n\times\cP_2(\mR^n),\mR^{n\times n})$. Consequently, Eq.(\ref{weaveeq}) admits a unique solution $\{\check{X}_t^\eps\}_{t\geq0}$.

Furthermore, arguing as (\ref{barb1grow}), we have
\ce
\|\bar{\s}_1(t,x,\mu)\|\leq C.
\de
Consequently, together with (\ref{barb1grow}), it follows that for any $0\leq s\leq t\leq T$,
\be
\mE|\check{X}_t^{\eps}-\check{X}_s^{\eps}|^2\leq  C_T(1+\mE|\xi|^2)|t-s|,\label{barxetsweak}
\ee
and
\be
\mE(\sup_{0\leq t\leq T}|\check{X}_t^{\eps}|^2)\leq C_T(1+\mE|\xi|^2).\label{barxeweak}
\ee

Let $\check{X}^{\eps,s,\xi}_t$ be the unique solution to Eq.\eqref{weaveeq} with initial data $\xi\in L^2(\Omega, \sF_s, \mP; \mR^n)$ at time $s$, i.e., for $t\geq s$,
\ce
\left\{\begin{array}{l}
\dif\check{X}^{\eps,s,\xi}_t=\bar{b}_1(t/\eps,\check{X}^{\eps,s,\xi}_t, \sL_{\check{X}^{\eps,s,\xi}_t})\dif t+\bar{\s}_1(t/\eps,\check{X}^{\eps,s,\xi}_t, \sL_{\check{X}^{\eps,s,\xi}_t})\dif \bar{B}_t, \\
\check{X}^{\eps,s,\xi}_s=\xi.
\end{array}
\right.
\de
The corresponding decoupled equation is given by
\ce
\left\{\begin{array}{l}
\dif\check{X}^{\eps,s,x,\sL_{\xi}}_t=\bar{b}_1(t/\eps,\check{X}^{\eps,s,x,\sL_{\xi}}_t, \sL_{\check{X}^{\eps,s,\xi}_t})\dif t+\bar{\s}_1(t/\eps,\check{X}^{\eps,s,x,\sL_{\xi}}_t, \sL_{\check{X}^{\eps,s,\xi}_t})\dif \bar{B}_t, \\
\check{X}^{\eps,s,x,\sL_{\xi}}_s=x.
\end{array}
\right.
\de
It is clear that $\check{X}^{\eps,s,\xi}_t=\check{X}^{\eps,s,x,\sL_{\xi}}_t\big|_{x=\xi}$. Combining the above estimates, we deduce that
\be\label{pxmubarx}
&&\sup_{0\leq s\leq t\leq T}\left(\sum_{i=1}^4\mE\|\p_x^i\check{X}^{\eps,s,x,\sL_{\xi}}_t\|^4+\sum_{i=0}^3\mE\|\p_{\tx}^i\p_\mu\check{X}^{\eps,s,x,\sL_{\xi}}_t(\tx)\|^2\right)\leq C_T,\no\\
&&\sup_{0\leq s\leq t\leq T}\left(\sum_{i=1}^2\mE\|\p_x^i\p_\mu\check{X}^{\eps,s,x,\sL_{\xi}}_t(\tx)\|^2+\sum_{i=1}^2\mE\|\p_x^i\p_{\tx}\p_\mu\check{X}^{\eps,s,x,\sL_{\xi}}_t(\tx)\|^2 \right)\leq C_T,\no\\
&&\sup_{0\leq s\leq t\leq T}\sum_{i,j=0}^1\mE\|\p_{\cx}^i\p_{\tx}^j\p_\mu^2\check{X}^{\eps,s,x,\sL_{\xi}}_t(\tx,\cx)\|^2\leq C_T.
\ee

Next, we introduce a Cauchy problem associated with the averaging dynamics. For any $\varphi\in (C_b^{(1,3)}\cap C_b^{(2,2)})(\cP_2(\mR^n),\mR)$, we consider the following Cauchy problem:
\be
\left\{\begin{array}{l}
\p_su^{\eps,t}(s,\sL_{\xi})+\mE[\bar{b}_1(s/\eps,\xi,\sL_{\xi})\cdot\p_\mu u^{\eps,t}(s,\sL_{\xi})(\xi)]\\
\qquad+\frac12\mE\left[Tr\big(\overline{\s_1\s_1^*}(s/\eps,\xi,\sL_{\xi})\cdot\p_{\tx}\p_\mu u^{\eps,t}(s,\sL_{\xi})(\xi)\big)\right]=0,\\
u^{\eps,t}(t,\sL_{\xi})=\varphi(\sL_{\xi}),\quad t\geq s.
\end{array}
\right.
\label{cauchypro}
\ee
The unique solution to (\ref{cauchypro}) is given by
\be
u^{\eps,t}(s,\sL_{\xi})=\varphi(\sL_{\check{X}^{\eps,s,\xi}_t}),\quad t\geq s.\label{cauchysol}
\ee
According to $\varphi\in (C_b^{(1,3)}\cap C_b^{(2,2)})(\cP_2(\mR^n),\mR)$ and (\ref{pxmubarx}), we know that for any $T>0$, there exists a constant $C_T>0$ such that for any $0\leq s\leq t\leq T$, $u^{\eps,t}(s,\cdot)\in (C_b^{(1,3)}\cap C_b^{(2,2)})(\cP_2(\mR^n),\mR)$, and
\be\label{pspmuuet}
\sup_{0\leq s\leq t\leq T}\|\p_s\p_{\mu}u^{\eps,t}(s,\mu)(\tx)\|+\sup_{0\leq s\leq t\leq T}\|\p_s\p_{\tx}\p_{\mu}u^{\eps,t}(s,\mu)(\tx)\|\leq C_T(1+|\tx|+\mu^{1/2}(|\cdot|^2)),
\ee
where the term $(1+|\tx|+\mu^{1/2}(|\cdot|^2))$ arises from the linear growth of $\bar{b}_1$ (cf. (\ref{barb1grow})).

{\bf Proof of Theorem \ref{xbarweak}.}
Applying the It\^{o} formula to (\ref{cauchysol}), we derive that
\ce
&&\varphi(\sL_{X_t^{\eps}})-\varphi(\sL_{\check{X}_t^\eps})
=u^{\eps,t}(t,\sL_{X_t^{\eps}})-u^{\eps,t}(0,\sL_{\xi})\\
&=&\int_0^t\p_su^{\eps,t}(s,\sL_{X_s^{\eps}})\dif s+\mE\int_0^tb_1(X_s^{\eps}, \sL_{X_s^{\eps}}, Y_s^{\eps})\cdot\p_\mu u^{\eps,t}(s,\sL_{X_s^{\eps}})(X_s^{\eps})\dif s\\
&&+\frac12\mE\int_0^tTr\big((\s_1\s_1^*)(X_s^{\eps}, \sL_{X_s^{\eps}}, Y_s^{\eps})\cdot\p_{\tx}\p_\mu u^{\eps,t}(s,\sL_{X_s^{\eps}})(X_s^{\eps})\big)\dif s\\
&=&\mE\int_0^t\left(b_1(X_s^{\eps}, \sL_{X_s^{\eps}}, Y_s^{\eps})-\bar{b}_1(s/\eps,X_s^{\eps}, \sL_{X_s^{\eps}})\right)\cdot\p_\mu u^{\eps,t}(s,\sL_{X_s^{\eps}})(X_s^{\eps})\dif s\\
&&+\frac12\mE\int_0^tTr\Big(\left[(\s_1\s_1^*)(X_s^{\eps}, \sL_{X_s^{\eps}}, Y_s^{\eps})-\overline{\s_1\s_1^*}(s/\eps,X_s^{\eps}, \sL_{X_s^{\eps}})\right]\\
&&\qquad\qquad\qquad\qquad\qquad\qquad\cdot\p_{\tx}\p_\mu u^{\eps,t}(s,\sL_{X_s^{\eps}})(X_s^{\eps})\Big)\dif s.
\de

For $r\in[0,t], s\geq0, x\in\mR^n, \mu \in \cP_2(\mR^n)$ and $y\in\mR^m$, define
\ce
H^t(r,s,x,\mu,y):&=&\left(b_1(x,\mu,y)-\bar{b}_1(s,x,\mu)\right)\cdot\p_\mu u^{\eps,t}(r,\mu)(x)\\
&&\qquad+\frac12Tr\left(\left[(\s_1\s_1^*)(x,\mu,y)-\overline{\s_1\s_1^*}(s,x,\mu)\right]\cdot\p_{\tx}\p_\mu u^{\eps,t}(r,\mu)(x)\right).
\de
By $(\mathbf{H}^2_{b_{1}, \s_{1}})$, (\ref{pspmuuet}) and $u^{\eps,t}(s,\cdot)\in (C_b^{(1,3)}\cap C_b^{(2,2)})(\cP_2(\mR^n),\mR)$, $H^t$ satisfies (\ref{pxpypmug}) and the centering condition
\ce
\int_{\mR^m}H^t(r,s,x,\mu,y)\nu_s^{x,\mu}(\dif y)=0,\quad r\in[0,t], s\geq0, x\in\mR^n, \mu\in\cP_2(\mR^n),\label{Hcencon}
\de
as well as
\be
\sup_{0\leq r\leq t\leq T,s\geq0,y\in\mR^m}\|\p_r\p_yH^t(r,s,x,\mu,y)\|\leq C_T\left(1+|x|+\mu^{1/2}(|\cdot|^2)\right).\label{prpyht}
\ee
Let $\cL^{x,\mu}(s)$ be the operator given by (\ref{generator}). Consider the nonautonomous Poisson equation
\be
\p_s\tilde{\Psi}^t(r,s,x,\mu,y)+\cL^{x,\mu}(s)\tilde{\Psi}^t(r,s,x,\mu,\cdot)(y)=-H^t(r,s,x,\mu,y).\label{Htpoieq}
\ee
By Proposition \ref{poiprop}, it admits a solution
\ce
\tilde{\Psi}^t(r,s,x,\mu,y):=\int_s^{+\infty}\mE H^t(r,u,x,\mu,Y_u^{s,x,\mu,y})\dif u.
\de
Moreover, for any $T>0$, there exists $C_T>0$ such that for any $s\geq0, x,\tx\in\mR^n, \mu\in\cP_2(\mR^n), y\in\mR^m$,
\be
&&\sup_{0\leq r\leq t\leq T}\max\Big\{|\tilde{\Psi}^t(r,s,x,\mu,y)|,\|\p_x\tilde{\Psi}^t(r,s,x,\mu,y)\|,\|\p_x^2\tilde{\Psi}^t(r,s,x,\mu,y)\|,\no\\
&&\qquad\qquad\qquad\qquad\|\p_\mu\tilde{\Psi}^t(r,s,x,\mu,y)(\tx)\|,\|\p_{\tx}\p_\mu\tilde{\Psi}^t(r,s,x,\mu,y)(\tx)\|\Big\}\no\\
&\leq& C_T\left(1+|x|+\mu^{1/2}(|\cdot|^2)+|y|\right)\Upsilon_{\g}(s).\label{pxpmutpsi}
\ee
Furthermore, by (\ref{prpyht}) and Lemma \ref{yvs}, we obtain
\be
&&\sup_{0\leq r\leq t\leq T}|\p_r\tilde{\Psi}^t(r,s,x,\mu,y)|\no\\
&=&\sup_{0\leq r\leq t\leq T}\left|\int_s^{+\infty}\mE\left[\p_rH^t(r,u,x,\mu,Y_u^{s,x,\mu,y})-\p_rH^t(r,u,x,\mu,\varsigma_u^{x,\mu})\right]\dif u\right|\no\\
&\leq&C_T\left(1+|x|+\mu^{1/2}(|\cdot|^2)\right)\int_s^{+\infty}\mE|Y_u^{s,x,\mu,y}-\varsigma_u^{x,\mu}|\dif u\no\\
&\leq&C_T\left(1+|x|^2+\mu(|\cdot|^2)+|y|^2\right)\Upsilon_{\g}(s).\label{prtpsi}
\ee

Thus,
\ce
&&\sup_{0\leq t\leq T}|\varphi(\sL_{X_t^{\eps}})-\varphi(\sL_{\check{X}_t^\eps})|\\
&=&\sup_{0\leq t\leq T}\left|\mE\int_0^tH^t(s,s/\eps,X_s^{\eps},\sL_{X_s^{\eps}},Y_s^{\eps})\dif s\right|\\
&=&\sup_{0\leq t\leq T}\left|\mE\int_0^t\left[\p_s\tilde{\Psi}^t(s,s/\eps,X_s^{\eps},\sL_{X_s^{\eps}},Y_s^{\eps})+\cL^{X_s^{\eps},\sL_{X_s^{\eps}}}(s/\eps)\tilde{\Psi}^t(s,s/\eps,X_s^{\eps},\sL_{X_s^{\eps}},\cdot)(Y_s^{\eps})\right]\dif s\right|.
\de
Applying It\^{o}'s formula to $\tilde{\Psi}^t$, we get
\ce
&&\mE\tilde{\Psi}^t(t,t/\eps,X_t^{\eps}, \sL_{X_t^{\eps}}, Y_t^{\eps})\\
&=&\mE\tilde{\Psi}^t(0,0,\xi,\sL_{\xi},\varrho)+\mE\int_0^t\p_r\tilde{\Psi}^t(s,s/\eps,X_s^{\eps},\sL_{X_s^{\eps}},Y_s^{\eps})\dif s\\
&&+\mE\int_0^t\cL_1^{\sL_{X_s^{\eps}}, Y_s^{\eps}}\tilde{\Psi}^t(s,s/\eps,\cdot,\sL_{X_s^{\eps}},Y_s^{\eps})(X_s^{\eps})\dif s\\
&&+\mE\tilde{\mE}\int_0^tb_1(\tilde{X}_s^{\eps}, \sL_{X_s^{\eps}}, \tilde{Y}_s^{\eps})\cdot\p_\mu\tilde{\Psi}^t(s,s/\eps,X_s^{\eps},\sL_{X_s^{\eps}},Y_s^{\eps})(\tilde{X}_s^{\eps})\dif s\\
&&+\frac12\mE\tilde{\mE}\int_0^tTr\big[(\s_1\s_1^*)(\tilde{X}_s^{\eps}, \sL_{X_s^{\eps}}, \tilde{Y}_s^{\eps})\cdot\p_{\tx}\p_\mu\tilde{\Psi}^t(s,s/\eps,X_s^{\eps},\sL_{X_s^{\eps}},Y_s^{\eps})(\tilde{X}_s^{\eps})\big]\dif s\\
&&+\eps^{-1}\mE\int_0^t\Big[\p_s\tilde{\Psi}^t(s,s/\eps,X_s^{\eps},\sL_{X_s^{\eps}},Y_s^{\eps})\\
&&\qquad\qquad\qquad+\cL^{X_s^{\eps},\sL_{X_s^{\eps}}}(s/\eps)\tilde{\Psi}^t(s,s/\eps,X_s^{\eps},\sL_{X_s^{\eps}},\cdot)(Y_s^{\eps})\Big]\dif s,
\de
which implies that
\ce
&&-\mE\int_0^t\Big[\p_s\tilde{\Psi}^t(s,s/\eps,X_s^{\eps},\sL_{X_s^{\eps}},Y_s^{\eps})+\cL^{X_s^{\eps},\sL_{X_s^{\eps}}}(s/\eps)\tilde{\Psi}^t(s,s/\eps,X_s^{\eps},\sL_{X_s^{\eps}},\cdot)(Y_s^{\eps})\Big]\dif s\\
&=&\eps\Big\{\mE\tilde{\Psi}^t(0,0,\xi,\sL_{\xi},\varrho)-\mE\tilde{\Psi}^t(t,t/\eps,X_t^{\eps}, \sL_{X_t^{\eps}}, Y_t^{\eps})+\mE\int_0^t\p_r\tilde{\Psi}^t(s,s/\eps,X_s^{\eps},\sL_{X_s^{\eps}},Y_s^{\eps})\dif s\\
&&\qquad+\mE\int_0^t\cL_1^{\sL_{X_s^{\eps}}, Y_s^{\eps}}\tilde{\Psi}^t(s,s/\eps,\cdot,\sL_{X_s^{\eps}},Y_s^{\eps})(X_s^{\eps})\dif s\\
&&\qquad+\mE\tilde{\mE}\int_0^tb_1(\tilde{X}_s^{\eps}, \sL_{X_s^{\eps}}, \tilde{Y}_s^{\eps})\cdot\p_\mu\tilde{\Psi}^t(s,s/\eps,X_s^{\eps},\sL_{X_s^{\eps}},Y_s^{\eps})(\tilde{X}_s^{\eps})\dif s\\
&&\qquad+\frac12\mE\tilde{\mE}\int_0^tTr\big[(\s_1\s_1^*)(\tilde{X}_s^{\eps}, \sL_{X_s^{\eps}},\tilde{Y}_s^{\eps})\cdot\p_{\tx}\p_\mu\tilde{\Psi}^t(s,s/\eps,X_s^{\eps},\sL_{X_s^{\eps}},Y_s^{\eps})(\tilde{X}_s^{\eps})\big]
\dif s\Big\}.
\de
Therefore,
\ce
&&\sup_{0\leq t\leq T}|\varphi(\sL_{X_t^{\eps}})-\varphi(\sL_{\check{X}_t^\eps})|\\
&\leq&\eps\Big\{\sup_{0\leq t\leq T}\mE|\tilde{\Psi}^t(0,0,\xi,\sL_{\xi},\varrho)|+\sup_{0\leq t\leq T}\mE|\tilde{\Psi}^t(t,t/\eps,X_t^{\eps}, \sL_{X_t^{\eps}}, Y_t^{\eps})|\\
&&\qquad+\sup_{0\leq t\leq T}\mE\left|\int_0^t\p_r\tilde{\Psi}^t(s,s/\eps,X_s^{\eps},\sL_{X_s^{\eps}},Y_s^{\eps})\dif s\right|\\
&&\qquad+\sup_{0\leq t\leq T}\mE\left|\int_0^t\cL_1^{\sL_{X_s^{\eps}}, Y_s^{\eps}}\tilde{\Psi}^t(s,s/\eps,\cdot,\sL_{X_s^{\eps}},Y_s^{\eps})(X_s^{\eps})\dif s\right|\\
&&\qquad+\sup_{0\leq t\leq T}\mE\tilde{\mE}\left|\int_0^tb_1(\tilde{X}_s^{\eps}, \sL_{X_s^{\eps}}, \tilde{Y}_s^{\eps})\cdot\p_\mu\tilde{\Psi}^t(s,s/\eps,X_s^{\eps},\sL_{X_s^{\eps}},Y_s^{\eps})(\tilde{X}_s^{\eps})\dif s\right|\\
&&\qquad+\frac12\sup_{0\leq t\leq T}\mE\tilde{\mE}\left|\int_0^tTr\big[(\s_1\s_1^*)(\tilde{X}_s^{\eps}, \sL_{X_s^{\eps}},\tilde{Y}_s^{\eps})\cdot\p_{\tx}\p_\mu\tilde{\Psi}^t(s,s/\eps,X_s^{\eps},\sL_{X_s^{\eps}},Y_s^{\eps})(\tilde{X}_s^{\eps})\big]\dif s\right|\Big\}\\
&=&\eps\sum_{i=1}^6\cJ_i.
\de
For $\cJ_1$ and $\cJ_2$, by (\ref{pxpmutpsi}) and Lemma \ref{xeye}, it holds that
\ce
&&\cJ_1+\cJ_2\\
&\leq&C_T\left(1+(\mE|\xi|^2)^{1/2}+\mE|\varrho|\right)\Upsilon_{\g}(0)+C_T\sup_{0\leq t\leq T}\left[\left(1+(\mE|X_t^{\eps}|^2)^{1/2}+\mE|Y_t^{\eps}|\right)\Upsilon_{\g}(t/\eps)\right]\\
&\leq& C_T(1+\mE|\xi|^4+\mE|\varrho|^4)\sup_{0\leq t\leq T}\Upsilon_{\g}(t/\eps).
\de
For $\cJ_3$, (\ref{prtpsi}) together with Lemma \ref{xeye} implies that
\ce
\cJ_3&\leq&C_T\int_0^T\left(1+\mE|X_s^{\eps}|^2+\mE|Y_s^{\eps}|^2\right)\Upsilon_{\g}(s/\eps)\dif s\\
&\leq& C_T(1+\mE|\xi|^4+\mE|\varrho|^4)\sup_{0\leq t\leq T}\Upsilon_{\g}(t/\eps).
\de
Similarly, we arrive at
\ce
\cJ_4&\leq&\sup_{0\leq t\leq T}\mE\int_0^t\Big|b_1(X_s^{\eps},\sL_{X_s^{\eps}},Y_s^{\eps})\cdot\p_x\tilde{\Psi}^t(s,s/\eps,X_s^{\eps},\sL_{X_s^{\eps}},Y_s^{\eps})\\
&&\qquad\qquad\qquad\qquad+\frac12Tr[(\s_1\s_1^*)(X_s^{\eps},\sL_{X_s^{\eps}},Y_s^{\eps})\cdot\p_x^2\tilde{\Psi}^t(s,s/\eps,X_s^{\eps},\sL_{X_s^{\eps}},Y_s^{\eps})]
\Big|\dif s\\
&\leq&C_T\int_0^T\left(1+\mE|X_s^{\eps}|^2+\mE|Y_s^{\eps}|^2\right)\Upsilon_{\g}(s/\eps)\dif s\\
&\leq& C_T(1+\mE|\xi|^4+\mE|\varrho|^4)\sup_{0\leq t\leq T}\Upsilon_{\g}(t/\eps),
\de
and
\ce
\cJ_5+\cJ_6\leq C_T(1+\mE|\xi|^4+\mE|\varrho|^4)\sup_{0\leq t\leq T}\Upsilon_{\g}(t/\eps).
\de
From the above deduction, it follows that
\ce
\sup_{0\leq t\leq T}|\varphi(\sL_{X_t^{\eps}})-\varphi(\sL_{\check{X}_t^\eps})|\leq C_T(1+\mE|\xi|^4+\mE|\varrho|^4)\eps\sup_{0\leq t\leq T}\Upsilon_{\g}(t/\eps).
\de
The proof is complete.

\section{Proofs of Theorem \ref{xbarx2con} and \ref{xbarweakcon}}\label{convcoefproo}

In this section, we prove Theorem \ref{xbarx2con} and \ref{xbarweakcon}.

\subsection{Asymptotic behavior of evolution system of measures}\label{asyevo}
In this subsection, we study the asymptotic behavior of the evolution system of measures.

According to $(\mathbf{H}^3_{b_{2}, \s_{2}})$ and $(\mathbf{H}^5_{b_{2}, \s_{2}})$, it holds that for all $x_1, x_2\in\mR^n, \mu_1, \mu_2\in\cP_2(\mR^n)$ and $y_1, y_2\in\mR^m$,
\be
&&2\langle y_1-y_2, \bar{b}_2(x_1,\mu_1,y_1)-\bar{b}_2(x_2,\mu_2,y_2)\rangle+3\|\bar{\s}_2(x_1,\mu_1,y_1)-\bar{\s}_2(x_2,\mu_2,y_2)\|^2\no\\
&\leq&-2\alpha|y_1-y_2|^2+C(|x_1-x_2|^2+\mW_2^2(\mu_1,\mu_2)),
\label{conbarb2s2dis}
\ee
and for all $x\in\mR^n, \mu\in\cP_2(\mR^n)$ and $y\in\mR^m$,
\ce
|\bar{b}_2(x,\mu,y)|+\|\bar{\s}_2(x,\mu,y)\|\leq C(1+|x|+\mu^{1/2}(|\cdot|^2)+|y|).
\de
Consequently, the SDE
\be
\dif \bar{Y}_t=\bar{b}_2(x,\mu,\bar{Y}_t)\dif t+\bar{\s}_2(x,\mu,\bar{Y}_t)\dif W_t, \quad \bar{Y}_0=y\in\mR^m,\label{confrozeq}
\ee
admits a unique strong solution $\{\bar{Y}_t^{x,\mu,y}\}_{t\geq0}$. Moreover, by the argument used in the proof of (\ref{Y1Y2l4}), one obtains that for any $t\geq0, x_1, x_2\in\mR^n, \mu_1, \mu_2\in\cP_2(\mR^n)$ and $y_1, y_2\in\mR^m$,
\be
\mE|\bar{Y}_t^{x_1,\mu_1,y_1}-\bar{Y}_t^{x_2,\mu_2,y_2}|^2\leq |y_1-y_2|^2e^{-2\alpha t}+C(|x_1-x_2|^2+\mW_2^2(\mu_1,\mu_2)).
\label{conbarY}
\ee
Let $\{\bar{P}_t^{x,\mu}\}_{t\geq0}$ denote the transition semigroup associated with $\{\bar{Y}_t^{x,\mu,y}\}_{t\geq0}$, i.e., for any bounded measurable function $\psi:\mR^m\rightarrow\mR$,
\ce
\bar{P}_t^{x,\mu}\psi(y)=\mathbb{E}\psi(\bar{Y}_t^{x,\mu,y}),\quad y\in\mR^m, t\geq0.
\de
The semigroup $\{\bar{P}_t^{x,\mu}\}_{t\geq0}$ admits a unique invariant measure $\nu^{x,\mu}$ satisfying
\ce
\int_{\mR^m}|y|^2\nu^{x,\mu}(\dif y)\leq C(1+|x|^2+\mu(|\cdot|^2)).
\de
Combining this with (\ref{conbarY}) yields that for any Lipschitz function $\psi$ on $\mR^m$,
\be
|\mE\psi(\bar{Y}_t^{x,\mu,y})-\nu^{x,\mu}(\psi)|\leq CLip(\psi)(1+|x|+\mu^{1/2}(|\cdot|^2)+|y|)e^{-\alpha t}.\label{conbarYnu}
\ee

\bl\label{nutnu}
Suppose that $(\mathbf{H}^3_{b_{2}, \s_{2}})$ and $(\mathbf{H}^5_{b_{2}, \s_{2}})$ hold. Then for any Lipschitz function $\psi$ and $\beta\in(0,1)$,
\be
|\nu_t^{x,\mu}(\psi)-\nu^{x,\mu}(\psi)|&\leq& CLip(\psi)(1+|x|+\mu^{1/2}(|\cdot|^2))e^{-\beta\alpha t}\no\\
&&+C_{\beta}Lip(\psi)(1+|x|+\mu^{1/2}(|\cdot|^2))\left[\int_0^te^{-2\beta\alpha(t-r)}\lambda^2(r)\dif r\right]^{1/2},\no\\
\label{nutnuesti}
\ee
where $\{\nu_t^{x,\mu}\}_{t\geq0}$ is an evolution system of measures of the semigroup $\{P_{s,t}^{x,\mu}\}_{t\geq s}$ given in Lemma \ref{nuxmu}.
\el
\begin{proof}
From (\ref{confrozeq}) and (\ref{vseq}), it follows that
\ce\left\{\begin{array}{ll}
\dif(\varsigma_t^{x,\mu}-\bar{Y}_t^{x,\mu,y})=(b_2(t,x,\mu,\varsigma_t^{x,\mu})-\bar{b}_2(x,\mu,\bar{Y}_t^{x,\mu,y}))\dif t\\
\qquad\qquad\qquad\qquad +(\s_2(t,x,\mu,\varsigma_t^{x,\mu})-\bar{\s}_2(x,\mu,\bar{Y}_t^{x,\mu,y}))\dif W_t,\\
\varsigma_0^{x,\mu}-\bar{Y}_0^{x,\mu,y}=\varsigma_0^{x,\mu}-y.
\end{array}
\right.
\de
Applying It\^{o}'s formula, together with (\ref{conbarb2s2dis}), (\ref{yvsl2}) and $(\mathbf{H}^5_{b_{2}, \s_{2}})$, we obtain for any $\beta\in(0,1)$,
\ce
&&\frac{\dif}{\dif t}\mE|\varsigma_t^{x,\mu}-\bar{Y}_t^{x,\mu,y}|^2\\
&=&2\mE\langle\varsigma_t^{x,\mu}-\bar{Y}_t^{x,\mu,y}, b_2(t,x,\mu,\varsigma_t^{x,\mu})-\bar{b}_2(x,\mu,\bar{Y}_t^{x,\mu,y})\rangle\\
&&+\mE\|\s_2(t,x,\mu,\varsigma_t^{x,\mu})-\bar{\s}_2(x,\mu,\bar{Y}_t^{x,\mu,y})\|^2\\
&\leq&2\mE\langle\varsigma_t^{x,\mu}-\bar{Y}_t^{x,\mu,y}, b_2(t,x,\mu,\varsigma_t^{x,\mu})-\bar{b}_2(x,\mu,\varsigma_t^{x,\mu})\rangle\\
&&+2\mE\langle\varsigma_t^{x,\mu}-\bar{Y}_t^{x,\mu,y}, \bar{b}_2(x,\mu,\varsigma_t^{x,\mu})-\bar{b}_2(x,\mu,\bar{Y}_t^{x,\mu,y})\rangle\\
&&+2\mE\|\s_2(t,x,\mu,\varsigma_t^{x,\mu})-\bar{\s}_2(x,\mu,\varsigma_t^{x,\mu})\|^2+2\mE\|\bar{\s}_2(x,\mu,\varsigma_t^{x,\mu})-\bar{\s}_2(x,\mu,\bar{Y}_t^{x,\mu,y})\|^2\\
&\leq&-2\beta\alpha\mE|\varsigma_t^{x,\mu}-\bar{Y}_t^{x,\mu,y}|^2+C_{\beta}\lambda^2(t)(1+|x|^2+\mu(|\cdot|^2)+\mE|\varsigma_t^{x,\mu}|^2)\\
&\leq&-2\beta\alpha\mE|\varsigma_t^{x,\mu}-\bar{Y}_t^{x,\mu,y}|^2+C_{\beta}\lambda^2(t)(1+|x|^2+\mu(|\cdot|^2)).
\de
This along with (\ref{yvsl2}) again implies that
\ce
&&\mE|\varsigma_t^{x,\mu}-\bar{Y}_t^{x,\mu,y}|^2\\
&\leq&e^{-2\beta\alpha t}\mE|\varsigma_0^{x,\mu}-y|^2+C_{\beta}(1+|x|^2+\mu(|\cdot|^2))\int_0^te^{-2\beta\alpha(t-r)}\lambda^2(r)\dif r\\
&\leq&C(1+|x|^2+\mu(|\cdot|^2)+|y|^2)e^{-2\beta\alpha t}+C_{\beta}(1+|x|^2+\mu(|\cdot|^2))\int_0^te^{-2\beta\alpha(t-r)}\lambda^2(r)\dif r.
\de
Combining this with (\ref{conbarYnu}), we get that for any Lipschitz function $\psi$,
\ce
&&|\nu_t^{x,\mu}(\psi)-\nu^{x,\mu}(\psi)|\\
&\leq&\left|\mE\psi(\varsigma_t^{x,\mu})-\mE\psi(\bar{Y}_t^{x,\mu,0})\right|+\left|\mE\psi(\bar{Y}_t^{x,\mu,0})-\nu^{x,\mu}(\psi)\right|\\
&\leq&Lip(\psi)\mE|\varsigma_t^{x,\mu}-\bar{Y}_t^{x,\mu,0}|+CLip(\psi)(1+|x|+\mu^{1/2}(|\cdot|^2))e^{-\alpha t}\\
&\leq&CLip(\psi)(1+|x|+\mu^{1/2}(|\cdot|^2))e^{-\beta\alpha t}\\
&&+C_{\beta}Lip(\psi)(1+|x|+\mu^{1/2}(|\cdot|^2))\left[\int_0^te^{-2\beta\alpha(t-r)}\lambda^2(r)\dif r\right]^{1/2}.
\de
The proof is complete.
\end{proof}

\br
Under $(\mathbf{H}^3_{b_{2}, \s_{2}})$ and $(\mathbf{H}^5_{b_{2}, \s_{2}})$, $|\nu_t^{x,\mu}(\psi)-\nu^{x,\mu}(\psi)|\rightarrow0$ as $t\rightarrow+\infty$ for any Lipschitz function $\psi$ (see Lemma \ref{nutnu} and \cite[Remark 4.2]{swx}). Consequently,
\ce
\lim_{t\rightarrow+\infty}\bar{b}_1(t,x,\mu)=\lim_{t\rightarrow+\infty}\int_{\mR^m}b_1(x,\mu,y)\nu_t^{x,\mu}(\dif y)=\int_{\mR^m}b_1(x,\mu,y)\nu^{x,\mu}(\dif y),
\de
and
\ce
\lim_{t\rightarrow+\infty}\overline{\s_1\s_1^*}(t,x,\mu)=\lim_{t\rightarrow+\infty}\int_{\mR^m}(\s_1\s_1^*)(x,\mu,y)\nu_t^{x,\mu}(\dif y)=\int_{\mR^m}(\s_1\s_1^*)(x,\mu,y)\nu^{x,\mu}(\dif y).
\de
This allows us to construct averaging equations whose coefficients are independent of $\eps$.
\er

\subsection{Proof of Theorem \ref{xbarx2con}}\label{stravecon}
In this subsection, we require that $\sigma_1$ is independent of $y$, and prove Theorem \ref{xbarx2con}.

\bl\label{conbarx}
 Under $(\mathbf{H}^1_{b_{1}, \s_{1}})$, $(\mathbf{H}^3_{b_{2}, \s_{2}})$ and $(\mathbf{H}^5_{b_{2}, \s_{2}})$, Eq.(\ref{constraveeq}) admits a unique strong solution $\{\bar{X}_t\}_{t\geq0}$. Moreover, for $T>0$, there exists a positive constant $C_T$ such that
\be
\mE(\sup_{0\leq t\leq T}|\bar{X}_t|^2)\leq C_T(1+\mE|\xi|^2).
\label{conbarxl2}
\ee
\el
\begin{proof}
For $x_1, x_2\in\mR^n, \mu_1, \mu_2\in\cP_2(\mR^n)$, by (\ref{conbarY}) and (\ref{conbarYnu}), one has
\ce
&&|\bar{b}_{1,c}(x_1,\mu_1)-\bar{b}_{1,c}(x_2,\mu_2)|\\
&=&\left|\int_{\mR^m}b_1(x_1,\mu_1,y)\nu^{x_1,\mu_1}(\dif y)-\int_{\mR^m}b_1(x_2,\mu_2,y)\nu^{x_2,\mu_2}(\dif y)\right|\\
&\leq&\left|\int_{\mR^m}b_1(x_1,\mu_1,y)\nu^{x_1,\mu_1}(\dif y)-\mE b_1(x_1,\mu_1,\bar{Y}_t^{x_1,\mu_1,0})\right|\\
&&+\left|\mE b_1(x_1,\mu_1,\bar{Y}_t^{x_1,\mu_1,0})-\mE b_1(x_2,\mu_2,\bar{Y}_t^{x_2,\mu_2,0})\right|\\
&&+\left|\mE b_1(x_2,\mu_2,\bar{Y}_t^{x_2,\mu_2,0})-\int_{\mR^m}b_1(x_2,\mu_2,y)\nu^{x_2,\mu_2}(\dif y)\right|\\
&\leq&C\left(1+|x_1|+\mu_1^{1/2}(|\cdot|^2)+|x_2|+\mu_2^{1/2}(|\cdot|^2)\right)e^{-\alpha t}\\
&&+C(|x_1-x_2|+\mW_2(\mu_1,\mu_2)).
\de
Letting $t\rightarrow+\infty$, we obtain that
\be
|\bar{b}_{1,c}(x_1,\mu_1)-\bar{b}_{1,c}(x_2,\mu_2)|\leq C(|x_1-x_2|+\mW_2(\mu_1,\mu_2)).\label{barb1clip}
\ee
The result follows from standard arguments. The proof is complete.
\end{proof}

{\bf Proof of Theorem \ref{xbarx2con}.}
Note that
\ce
\bar{X}_t^{\eps}-\bar{X}_t&=&\int_0^t\left[\bar{b}_1(s/\eps,\bar{X}_s^\eps,\sL_{\bar{X}_s^\eps})-\bar{b}_{1,c}(\bar{X}_s,\sL_{\bar{X}_s})\right]\dif s\\
&&+\int_0^t[\s_1(\bar{X}_s^\eps,\sL_{\bar{X}_s^\eps})-\s_1(\bar{X}_s,\sL_{\bar{X}_s})]\dif B_s.
\de
(\ref{barb1clip}) implies that for all $t\in[0,v]$,
\ce
&&\mE(\sup_{0\leq t\leq v}|\bar{X}_t^{\eps}-\bar{X}_t|^2)\\
&\leq&C\mE\left(\sup_{0\leq t\leq v}\left|\int_0^t\left[\bar{b}_1(s/\eps,\bar{X}_s^\eps,\sL_{\bar{X}_s^\eps})-\bar{b}_{1,c}(\bar{X}_s^\eps,\sL_{\bar{X}_s^\eps})\right]\dif s\right|^2\right)\\
&&+C\mE\left(\sup_{0\leq t\leq v}\left|\int_0^t\left[\bar{b}_{1,c}(\bar{X}_s^\eps,\sL_{\bar{X}_s^\eps})-\bar{b}_{1,c}(\bar{X}_s,\sL_{\bar{X}_s})\right]\dif s\right|^2\right)\\
&&+C\mE\left(\sup_{0\leq t\leq v}\left|\int_0^t[\s_1(\bar{X}_s^\eps,\sL_{\bar{X}_s^\eps})-\s_1(\bar{X}_s,\sL_{\bar{X}_s})]\dif B_s\right|^2\right)\\
&\leq&C\mE\left(\sup_{0\leq t\leq v}\left|\int_0^t\left[\bar{b}_1(s/\eps,\bar{X}_s^\eps,\sL_{\bar{X}_s^\eps})-\bar{b}_{1,c}(\bar{X}_s^\eps,\sL_{\bar{X}_s^\eps})\right]\dif s\right|^2\right)\\
&&+C_T\mE\int_0^v\left|\bar{b}_{1,c}(\bar{X}_s^\eps,\sL_{\bar{X}_s^\eps})-\bar{b}_{1,c}(\bar{X}_s,\sL_{\bar{X}_s})\right|^2\dif s\\
&&+C_T\mE\int_0^v\left\|\s_1(\bar{X}_s^\eps,\sL_{\bar{X}_s^\eps})-\s_1(\bar{X}_s,\sL_{\bar{X}_s})\right\|^2\dif s\\
&\leq&C\mE\left(\sup_{0\leq t\leq v}\left|\int_0^t\left[\bar{b}_1(s/\eps,\bar{X}_s^\eps,\sL_{\bar{X}_s^\eps})-\bar{b}_{1,c}(\bar{X}_s^\eps,\sL_{\bar{X}_s^\eps})\right]\dif s\right|^2\right)+C_T\int_0^v\mE|\bar{X}_s^{\eps}-\bar{X}_s|^2\dif s.
\de
Applying the Gronwall inequality yields
\ce
\mE(\sup_{0\leq t\leq v}|\bar{X}_t^{\eps}-\bar{X}_t|^2)\leq C_T\mE\left(\sup_{0\leq t\leq v}\left|\int_0^t\left[\bar{b}_1(s/\eps,\bar{X}_s^\eps,\sL_{\bar{X}_s^\eps})-\bar{b}_{1,c}(\bar{X}_s^\eps,\sL_{\bar{X}_s^\eps})\right]\dif s\right|^2\right)
\de
According to Lemma \ref{nutnu}, for any $\beta\in(0,1)$, there exists $C_\beta>0$ such that
\be\label{barb1barb1c}
&&|\bar{b}_1(t,x,\mu)-\bar{b}_{1,c}(x,\mu)|\no\\
&\leq&\left|\int_{\mR^m}b_1(x,\mu,y)\nu_t^{x,\mu}(\dif y)-\int_{\mR^m}b_1(x,\mu,y)\nu^{x,\mu}(\dif y)\right|\no\\
&\leq& C_{\beta}(1+|x|+\mu^{1/2}(|\cdot|^2))\Big[e^{-\beta\alpha t}+\Big(\int_0^te^{-2\beta\alpha(t-r)}\lambda^2(r)\dif r\Big)^{1/2}\Big].
\ee
Consequently,
\ce
&&\mE\left(\sup_{0\leq t\leq v}\left|\int_0^t\left[\bar{b}_1(s/\eps,\bar{X}_s^\eps,\sL_{\bar{X}_s^\eps})-\bar{b}_{1,c}(\bar{X}_s^\eps,\sL_{\bar{X}_s^\eps})\right]\dif s\right|^2\right)\\
&\leq&\mE\left(\int_0^v|\bar{b}_1(s/\eps,\bar{X}_s^\eps,\sL_{\bar{X}_s^\eps})-\bar{b}_{1,c}(\bar{X}_s^\eps,\sL_{\bar{X}_s^\eps})|\dif s\right)^2\\
&\leq&C_T(1+\mE|\xi|^2)\left(\int_0^v\Big[e^{-\beta\alpha s/\eps}+\Big(\int_0^{s/\eps}e^{-2\beta\alpha(s/\eps-r)}\lambda^2(r)\dif r\Big)^{1/2}\Big]\dif s\right)^2\\
&\leq&C_{T,\beta}(1+\mE|\xi|^2)\eps^2\Bigg(1+\int_0^{v/\eps}\Big(\int_0^ue^{-2\beta\alpha(u-r)}\lambda^2(r)\dif r\Big)^{1/2}\dif u\Bigg)^2.
\de
Combining the above estimates, we obtain that
\ce
&&\mE(\sup_{0\leq t\leq T}|\bar{X}_t^{\eps}-\bar{X}_t|^2)\leq C_{T,\beta}(1+\mE|\xi|^2)\eps^2\Bigg(1+\int_0^{T/\eps}\Big(\int_0^ue^{-2\beta\alpha(u-r)}\lambda^2(r)\dif r\Big)^{1/2}\dif u\Bigg)^2.
\de
Together with (\ref{xebarxel2}), this yields that
\ce
&&\mE(\sup_{0\leq t\leq T}|X_t^{\eps}-\bar{X}_t|^2)\\
&\leq& 2\mE(\sup_{0\leq t\leq T}|X_t^{\eps}-\bar{X}_t^{\eps}|^2)+2\mE(\sup_{0\leq t\leq T}|\bar{X}_t^{\eps}-\bar{X}_t|^2)\\
&\leq&C_{T,\beta}(1+\mE|\xi|^4+\mE|\varrho|^4)\eps^2\Bigg[\sup_{0\leq t\leq T}\Upsilon^2_{\g}(t/\eps)+\sup_{0\leq t\leq T}\Upsilon^2(t/\eps)\left(\int_0^{T/\eps}\a(s)\dif s\right)^{1/2}\\
&&\qquad\qquad\qquad+\int_0^{T/\eps}\a(s)\Upsilon^2(s)\dif s+\Bigg(1+\int_0^{T/\eps}\Big(\int_0^ue^{-2\beta\alpha(u-r)}\lambda^2(r)\dif r\Big)^{1/2}\dif u\Bigg)^2\Bigg].
\de
The proof is complete.

\subsection{Proof of Theorem \ref{xbarweakcon}}\label{weakavecon}
In this subsection, we prove Theorem \ref{xbarweakcon}.

In view of the fact that, for any $t\geq0$, $\bar{b}_1(t,\cdot,\cdot)\in (C_b^{4,(1,3)}\cap C_b^{2,(2,2)})(\mR^n\times\cP_2(\mR^n),\mR^n)$ and $\bar{\s}_1(t,\cdot,\cdot)\in (C_b^{4,(1,3)}\cap C_b^{2,(2,2)})(\mR^n\times\cP_2(\mR^n),\mR^{n\times n})$, the coefficients $\bar{b}_{1,c}$ and $\bar{\s}_{1,c}$ satisfy $\bar{b}_{1,c}\in (C_b^{4,(1,3)}\cap C_b^{2,(2,2)})(\mR^n\times\cP_2(\mR^n),\mR^n)$ and $\bar{\s}_{1,c}\in (C_b^{4,(1,3)}\cap C_b^{2,(2,2)})(\mR^n\times\cP_2(\mR^n),\mR^{n\times n})$. Thus Eq.(\ref{conweaveeq}) has a unique solution $\{\check{X}_t\}_{t\geq0}$. Let $\check{X}^{s,\xi}_t$ be the unique solution to Eq.\eqref{conweaveeq} with initial data $\xi\in L^2(\Omega, \sF_s, \mP; \mR^n)$ at time $s$, i.e., for $t\geq s$,
\ce
\dif \check{X}^{s,\xi}_t=\bar{b}_{1,c}(\check{X}^{s,\xi}_t,\sL_{\check{X}^{s,\xi}_t})\dif t+\bar{\s}_{1,c}(\check{X}^{s,\xi}_t,\sL_{\check{X}^{s,\xi}_t})\dif \bar{B}_t,\quad \check{X}^{s,\xi}_s=\xi.
\de
The corresponding decoupled equation takes the form
\ce
\dif \check{X}^{s,x,\sL_{\xi}}_t=\bar{b}_{1,c}(\check{X}^{s,x,\sL_{\xi}}_t,\sL_{\check{X}^{s,\xi}_t})\dif t+\bar{\s}_{1,c}(\check{X}^{s,x,\sL_{\xi}}_t,\sL_{\check{X}^{s,\xi}_t})\dif \bar{B}_t,\quad \check{X}^{s,x,\sL_{\xi}}_s=x.
\de
Then $\check{X}^{s,\xi}_t=\check{X}^{s,x,\sL_{\xi}}_t\big|_{x=\xi}$. Moreover, for any $T>0$,
\be
&&\sup_{0\leq s\leq t\leq T}\left(\sum_{i=1}^4\mE\|\p_x^i\check{X}^{s,x,\sL_{\xi}}_t\|^4+\sum_{i=0}^3\mE\|\p_{\tx}^i\p_\mu\check{X}^{s,x,\sL_{\xi}}_t(\tx)\|^2\right)\leq C_T,\no\\
&&\sup_{0\leq s\leq t\leq T}\left(\sum_{i=1}^2\mE\|\p_x^i\p_\mu\check{X}^{s,x,\sL_{\xi}}_t(\tx)\|^2+\sum_{i=1}^2\mE\|\p_x^i\p_{\tx}\p_\mu\check{X}^{s,x,\sL_{\xi}}_t(\tx)\|^2 \right)\leq C_T,\no\\
&&\sup_{0\leq s\leq t\leq T}\sum_{i,j=0}^1\mE\|\p_{\cx}^i\p_{\tx}^j\p_\mu^2\check{X}^{s,x,\sL_{\xi}}_t(\tx,\cx)\|^2\leq C_T.
\label{pxmubarxcon}
\ee
For any $\varphi\in (C_b^{(1,3)}\cap C_b^{(2,2)})(\cP_2(\mR^n),\mR)$, we consider the following Cauchy problem on $[0,T]\times\cP_2(\mR^n)$ :
\be
\left\{\begin{array}{l}
\p_su^t(s,\sL_{\xi})+\mE[\bar{b}_{1,c}(\xi,\sL_{\xi})\cdot\p_\mu u^t(s,\sL_{\xi})(\xi)]\\
\qquad+\frac12\mE\left[Tr\big((\overline{\s_1\s_1^*})_c(\xi,\sL_{\xi})\cdot\p_{\tx}\p_\mu u^t(s,\sL_{\xi})(\xi)\big)\right]=0,\\
u^t(t,\sL_{\xi})=\varphi(\sL_{\xi}),\quad t\geq s.
\end{array}
\right.
\label{cauchyprocon}
\ee
The unique solution is given by
\be
u^t(s,\sL_{\xi})=\varphi(\sL_{\check{X}^{s,\xi}_t}),\quad t\geq s.\label{cauchyconsol}
\ee
Since $\varphi\in (C_b^{(1,3)}\cap C_b^{(2,2)})(\cP_2(\mR^n),\mR)$, it follows from (\ref{pxmubarxcon}) that, for any $T>0$,
\be\label{pmuutcon}
\sup _{0\leq s\leq t\leq T,\mu\in\cP_2(\mR^n),\tx\in\mR^n}\max\Big\{\|\p_{\mu}u^t(s,\mu)(\tx)\|,\|\p_{\tx}\p_{\mu}u^t(s,\mu)(\tx)\|\Big\}\leq C_T.
\ee

{\bf Proof of Theorem \ref{xbarweakcon}.}
Applying the It\^{o} formula to (\ref{cauchyconsol}) and using (\ref{cauchyprocon}), we derive that
\ce
&&\varphi(\sL_{\check{X}_t^\eps})-\varphi(\sL_{\check{X}_t})
=u^{t}(t,\sL_{\check{X}_t^\eps})-u^{t}(0,\sL_{\xi})\\
&=&\int_0^t\p_su^t(s,\sL_{\check{X}_s^\eps})\dif s+\mE\int_0^t\bar{b}_1(s/\eps,\check{X}_s^{\eps}, \sL_{\check{X}_s^{\eps}})\cdot\p_\mu u^t(s,\sL_{\check{X}_s^\eps})(\check{X}_s^\eps)\dif s\\
&&+\frac12\mE\int_0^tTr\big(\overline{\s_1\s_1^*}(s/\eps,\check{X}_s^{\eps}, \sL_{\check{X}_s^{\eps}})\cdot\p_{\tilde{x}}\p_\mu u^t(s,\sL_{\check{X}_s^\eps})(\check{X}_s^\eps)\big)\dif s\\
&=&\mE\int_0^t\left(\bar{b}_1(s/\eps,\check{X}_s^{\eps}, \sL_{\check{X}_s^{\eps}})-\bar{b}_{1,c}(\check{X}_s^{\eps}, \sL_{\check{X}_s^{\eps}})\right)\cdot\p_\mu u^t(s,\sL_{\check{X}_s^\eps})(\check{X}_s^\eps)\dif s\\
&&+\frac12\mE\int_0^tTr\Big(\left[\overline{\s_1\s_1^*}(s/\eps,\check{X}_s^{\eps}, \sL_{\check{X}_s^{\eps}})-(\overline{\s_1\s_1^*})_c(\check{X}_s^{\eps}, \sL_{\check{X}_s^{\eps}})\right]\\
&&\qquad\qquad\qquad\qquad\qquad\qquad\cdot\p_{\tilde{x}}\p_\mu u^t(s,\sL_{\check{X}_s^\eps})(\check{X}_s^\eps)\Big)\dif s.
\de
Arguing as in the proof of (\ref{barb1barb1c}), we have
\ce
&&\|\overline{\s_1\s_1^*}(t,x,\mu)-(\overline{\s_1\s_1^*})_c(x,\mu)\|\\
&\leq&C_{\beta}(1+|x|+\mu^{1/2}(|\cdot|^2))\Big[e^{-\beta\alpha t}+\Big(\int_0^te^{-2\beta\alpha(t-r)}\lambda^2(r)\dif r\Big)^{1/2}\Big],
\de
which together with (\ref{barb1barb1c}) and (\ref{pmuutcon}) yields that
\ce
&&\sup_{0\leq t\leq T}|\varphi(\sL_{\check{X}_t^\eps})-\varphi(\sL_{\check{X}_t})|\\
&\leq&C_{T,\beta}(1+\mE|\xi|^2)\int_0^T\Big[e^{-\beta\alpha s/\eps}+\Big(\int_0^{s/\eps}e^{-2\beta\alpha(s/\eps-r)}\lambda^2(r)\dif r\Big)^{1/2}\Big]\dif s\\
&\leq&C_{T,\beta}(1+\mE|\xi|^2)\eps\Bigg[1+\int_0^{T/\eps}\Big(\int_0^ue^{-2\beta\alpha(u-r)}\lambda^2(r)\dif r\Big)^{1/2}\dif u\Bigg].
\de
Combining this with (\ref{lxelbarxe}), we further obtain
\ce
&&\sup_{0\leq t\leq T}|\varphi(\sL_{X_t^{\eps}})-\varphi(\sL_{\check{X}_t})|\\
&\leq& C_{T,\beta}(1+\mE|\xi|^4+\mE|\varrho|^4)\eps\Bigg[1+\sup_{0\leq t\leq T}\Upsilon_{\g}(t/\eps)+\int_0^{T/\eps}\Big(\int_0^ue^{-2\beta\alpha(u-r)}\lambda^2(r)\dif r\Big)^{1/2}\dif u\Bigg].
\de
The proof is complete.

\section{Proofs of Theorem \ref{xbarstrper} and \ref{xbarweakper}}\label{periconvproo}

In this section, we prove Theorem \ref{xbarstrper} and \ref{xbarweakper}.

\subsection{Periodicity of the evolution system of measures}\label{evoper}

In this subsection we show that the evolution system of measures is periodic. We begin with a key lemma.

\bl[\cite{swx}]\label{fperlem}
Let $f$ be a bounded and $\tau$-periodic function on $\mR$. Then, for any $T>0$,
\ce
\sup_{a\in\mR}\left|\frac1T\int_a^{a+T}f(s)\dif s-\frac1\tau\int_0^\tau f(s)\dif s\right|\leq \frac{2\tau M}T,
\de
where $M=\sup\limits_{s\in[0,\tau]}|f(s)|$.
\el

\bl\label{nuxmuper}
Suppose that $(\mathbf{H}^3_{b_{2}, \s_{2}})$ and $(\mathbf{H}^6_{b_{2}, \s_{2}})$ hold. Then, the evolution system of measures $\{\nu_t^{x,\mu}\}_{t\in\mR}$ given in Lemma \ref{nuxmu} is $\tau$-periodic, i.e., for any $t\in\mR$, $\nu_{t+\tau}^{x,\mu}=\nu_t^{x,\mu}$.
\el

Since the proof of the above Lemma is similar to that for Lemma 5.2 in \cite{swx}, we omit it.

\br
By Lemma \ref{nuxmuper}, $\{\bar{b}_1(t,x,\mu)\}_{t\in\mR}$ and $\{\bar{\s}_1(t,x,\mu)\}_{t\in\mR}$ are also $\tau$-periodic. So, we can construct averaging equations whose coefficients are independent of $\eps$.
\er

\subsection{Proof of Theorem \ref{xbarstrper}}\label{straveper}
In this subsection, we require that $\s_1$ is independent of $y$ and prove Theorem \ref{xbarstrper}.

By (\ref{barb1lip}) and (\ref{barb1grow}), for any $x_1,x_2\in\mR^n, \mu_1,\mu_2\in\cP_2(\mR^n)$,
\be
|\bar{b}_{1,p}(x_1,\mu_1)-\bar{b}_{1,p}(x_2,\mu_2)|\leq C\left(|x_1-x_2|+\mW_2(\mu_1, \mu_2)\right),\label{barb1plip}
\ee
and for any $x\in\mR^n, \mu\in\cP_2(\mR^n)$,
\be
|\bar{b}_{1,p}(x,\mu)|\leq C\left(1+|x|+\mu^{1/2}(|\cdot|^2)\right).\label{barb1pgrow}
\ee
Consequently, Eq.(\ref{perstraveeq}) admits a unique solution $\{\tilde{X}_t\}_{t\geq0}$. Moreover, by (\ref{barb1pgrow}) and $(\mathbf{H}^1_{b_{1}, \s_{1}})$, for any $T>0$, there exists $C_T>0$ such that for any $0\leq s\leq t\leq T$,
\be
\mE|\tilde{X}_t-\tilde{X}_s|^2\leq C_T(1+\mE|\xi|^2)|t-s|,\label{barxtbarxsl2}
\ee
and
\be
\mE(\sup_{0\leq t\leq T}|\tilde{X}_t|^2)\leq C_T(1+\mE|\xi|^2).\label{barxtl2}
\ee

{\bf Proof of Theorem \ref{xbarstrper}.}
Note that
\ce
\bar{X}_t^{\eps}-\tilde{X}_t&=&\int_0^t\left[\bar{b}_1(r/\eps,\bar{X}_r^\eps,\sL_{\bar{X}_r^\eps})-\bar{b}_{1,p}(\tilde{X}_r,\sL_{\tilde{X}_r})\right]\dif r\\
&&+\int_0^t[\s_1(\bar{X}_r^\eps,\sL_{\bar{X}_r^\eps})-\s_1(\tilde{X}_r,\sL_{\tilde{X}_r})]\dif B_r.
\de
By (\ref{barb1lip}) and $(\mathbf{H}^1_{b_{1}, \s_{1}})$, for all $s\in[0,t]$,
\ce
&&\mE(\sup_{0\leq s\leq t}|\bar{X}_s^{\eps}-\tilde{X}_s|^2)\\
&\leq&C\mE\left(\sup_{0\leq s\leq t}\left|\int_0^s\left[\bar{b}_1(r/\eps,\bar{X}_r^\eps,\sL_{\bar{X}_r^\eps})-\bar{b}_1(r/\eps,\tilde{X}_r,\sL_{\tilde{X}_r})\right]\dif r\right|^2\right)\\
&&+C\mE\left(\sup_{0\leq s\leq t}\left|\int_0^s\left[\bar{b}_1(r/\eps,\tilde{X}_r,\sL_{\tilde{X}_r})-\bar{b}_{1,p}(\tilde{X}_r,\sL_{\tilde{X}_r})\right]\dif r\right|^2\right)\\
&&+C\mE\left(\sup_{0\leq s\leq t}\left|\int_0^s[\s_1(\bar{X}_r^\eps,\sL_{\bar{X}_r^\eps})-\s_1(\tilde{X}_r,\sL_{\tilde{X}_r})]\dif B_r\right|^2\right)\\
&\leq&C_T\mE\int_0^t\left|\bar{b}_1(r/\eps,\bar{X}_r^\eps,\sL_{\bar{X}_r^\eps})-\bar{b}_1(r/\eps,\tilde{X}_r,\sL_{\tilde{X}_r})\right|^2\dif r\\
&&+C\mE\left(\sup_{0\leq s\leq t}\left|\int_0^s\left[\bar{b}_1(r/\eps,\tilde{X}_r,\sL_{\tilde{X}_r})-\bar{b}_{1,p}(\tilde{X}_r,\sL_{\tilde{X}_r})\right]\dif r\right|^2\right)\\
&&+C_T\mE\int_0^t\left\|\s_1(\bar{X}_r^\eps,\sL_{\bar{X}_r^\eps})-\s_1(\tilde{X}_r,\sL_{\tilde{X}_r})\right\|^2\dif r\\
&\leq&C\mE\left(\sup_{0\leq s\leq t}\left|\int_0^s\left[\bar{b}_1(r/\eps,\tilde{X}_r,\sL_{\tilde{X}_r})-\bar{b}_{1,p}(\tilde{X}_r,\sL_{\tilde{X}_r})\right]\dif r\right|^2\right)+C_T\int_0^t\mE|\bar{X}_r^{\eps}-\tilde{X}_r|^2\dif r.
\de
An application of the Gronwall inequality yields
\ce
\mE(\sup_{0\leq s\leq T}|\bar{X}_s^{\eps}-\tilde{X}_s|^2)\leq C_T\mE\left(\sup_{0\leq s\leq T}\left|\int_0^s\left[\bar{b}_1(r/\eps,\tilde{X}_r,\sL_{\tilde{X}_r})-\bar{b}_{1,p}(\tilde{X}_r,\sL_{\tilde{X}_r})\right]\dif r\right|^2\right).
\de
Moreover,
\ce
&&\mE\left(\sup_{0\leq s\leq T}\left|\int_0^s\left[\bar{b}_1(r/\eps,\tilde{X}_r,\sL_{\tilde{X}_r})-\bar{b}_{1,p}(\tilde{X}_r,\sL_{\tilde{X}_r})\right]\dif r\right|^2\right)\\
&\leq& C\mE\left(\sup_{0\leq s\leq T}\left|\int_0^s\left[\bar{b}_1(r/\eps,\tilde{X}_r,\sL_{\tilde{X}_r})-\bar{b}_1(r/\eps,\tilde{X}_{r(\D)},\sL_{\tilde{X}_{r(\D)}})\right]\dif r\right|^2\right)\\
&&+C\mE\left(\sup_{0\leq s\leq T}\left|\int_0^s\left[\bar{b}_1(r/\eps,\tilde{X}_{r(\D)},\sL_{\tilde{X}_{r(\D)}})-\bar{b}_{1,p}(\tilde{X}_{r(\D)},\sL_{\tilde{X}_{r(\D)}})\right]\dif r\right|^2\right)\\
&&+C\mE\left(\sup_{0\leq s\leq T}\left|\int_0^s\left[\bar{b}_{1,p}(\tilde{X}_{r(\D)},\sL_{\tilde{X}_{r(\D)}})-\bar{b}_{1,p}(\tilde{X}_r,\sL_{\tilde{X}_r})\right]\dif r\right|^2\right)\\
&=:&\sum_{i=1}^3J_i(T),
\de
where $r(\D):=[r/\D]\D$ and $[r/\D]$ denotes the integer part of $r/\D$.

For $J_1(T)$ and $J_3(T)$, (\ref{barb1lip}), (\ref{barb1plip}) and (\ref{barxtbarxsl2}) yield that
\ce
J_1(T)+J_3(T)\leq C_T\int_0^T\mE|\tilde{X}_r-\tilde{X}_{r(\D)}|^2\dif r\leq C_T(1+\mE|\xi|^2)\D.
\de
We now estimate $J_2(T)$. Observe that
\ce
J_2(T)&\leq&C\mE\left(\sup_{0\leq s\leq T}\left|\sum_{k=0}^{[s/\D]-1}\int_{k\D}^{(k+1)\D}\left[\bar{b}_1(r/\eps,\tilde{X}_{k\D},\sL_{\tilde{X}_{k\D}})-\bar{b}_{1,p}(\tilde{X}_{k\D},\sL_{\tilde{X}_{k\D}})\right]\dif r\right|^2\right)\\
&&+C\mE\left(\sup_{0\leq s\leq T}\left|\int_{[s/\D]\D}^s\left[\bar{b}_1(r/\eps,\tilde{X}_{[s/\D]\D},\sL_{\tilde{X}_{[s/\D]\D}})-\bar{b}_{1,p}(\tilde{X}_{[s/\D]\D},\sL_{\tilde{X}_{[s/\D]\D}})\right]\dif r\right|^2\right)\\
&=:&J_{21}(T)+J_{22}(T).
\de
As for $J_{22}(T)$, by (\ref{barb1grow}), (\ref{barb1pgrow}) and (\ref{barxtl2}), it follows that
\ce
&&J_{22}(T)\\
&\leq&C\D\mE\left(\sup_{0\leq s\leq T}\int_{[s/\D]\D}^s\left(|\bar{b}_1(r/\eps,\tilde{X}_{[s/\D]\D},\sL_{\tilde{X}_{[s/\D]\D}})|^2+|\bar{b}_{1,p}(\tilde{X}_{[s/\D]\D},\sL_{\tilde{X}_{[s/\D]\D}})|^2\right)\dif r\right)\\
&\leq&C\D^2\left(1+\mE(\sup_{0\leq s\leq T}|\tilde{X}_{[s/\D]\D}|^2)\right)\\
&\leq&C_T(1+\mE|\xi|^2)\D^2.
\de
For $J_{21}(T)$, since $\bar{b}_1(\cdot,x,\mu)$ is $\tau$-periodic, combining (\ref{barb1grow}), (\ref{barxtl2}) and Lemma \ref{fperlem}, we obtain
\ce
&&J_{21}(T)\\
&=&C\D^2\mE\left(\sup_{0\leq s\leq T}\left|\sum_{k=0}^{[s/\D]-1}\left[\frac{\eps}{\D}\int_{\frac{k\D}{\eps}}^{\frac{(k+1)\D}{\eps}}\bar{b}_1(u,\tilde{X}_{k\D},\sL_{\tilde{X}_{k\D}})\dif u-\bar{b}_{1,p}(\tilde{X}_{k\D},\sL_{\tilde{X}_{k\D}})\right]\right|^2\right)\\
&\leq&C_T\D\mE\left(\sup_{0\leq s\leq T}\sum_{k=0}^{[s/\D]-1}
\left|\frac{\eps}{\D}\int_{\frac{k\D}{\eps}}^{\frac{(k+1)\D}{\eps}}\bar{b}_1(u,\tilde{X}_{k\D},\sL_{\tilde{X}_{k\D}})\dif u-\frac1\tau\int_0^\tau\bar{b}_1(u,\tilde{X}_{k\D},\sL_{\tilde{X}_{k\D}})\dif u\right|^2\right)\\
&\leq&C_T\frac{\eps^2}{\D}\sum_{k=0}^{[T/\D]-1}\mE(\sup_{0\leq u\leq \tau}|\bar{b}_1(u,\tilde{X}_{k\D},\sL_{\tilde{X}_{k\D}})|^2)\\
&\leq&C_T(1+\mE|\xi|^2)\frac{\eps^2}{\D^2}.
\de
Combining the above estimates gives
\ce
\mE(\sup_{0\leq s\leq T}|\bar{X}_s^{\eps}-\tilde{X}_s|^2)\leq C_T(1+\mE|\xi|^2)(\D+\frac{\eps^2}{\D^2}).
\de
Setting $\D=\eps^{2/3}$ yields
\ce
\mE(\sup_{0\leq s\leq T}|\bar{X}_s^{\eps}-\tilde{X}_s|^2)\leq C_T(1+\mE|\xi|^2)\eps^{2/3}.
\de
Combining this with (\ref{xebarxel2}), we conclude that
\ce
&&\mE(\sup_{0\leq t\leq T}|X_t^{\eps}-\tilde{X}_t|^2)\\
&\leq& C_T(1+\mE|\xi|^4+\mE|\varrho|^4)\eps^2\\
&&\cdot\left[\eps^{-4/3}+\sup_{0\leq t\leq T}\Upsilon^2_{\g}(t/\eps)+\sup_{0\leq t\leq T}\Upsilon^2(t/\eps)\left(\int_0^{T/\eps}\a(s)\dif s\right)^{1/2}+\int_0^{T/\eps}\a(s)\Upsilon^2(s)\dif s \right].
\de
The proof is complete.

\subsection{Proof of Theorem \ref{xbarweakper}}\label{weakaveper}
In this subsection, we prove Theorem \ref{xbarweakper}.

Since, for any $t\geq0$, $\bar{b}_1(t,\cdot,\cdot)\in (C_b^{4,(1,3)}\cap C_b^{2,(2,2)})(\mR^n\times\cP_2(\mR^n),\mR^n)$ and $\bar{\s}_1(t,\cdot,\cdot)\in (C_b^{4,(1,3)}\cap C_b^{2,(2,2)})(\mR^n\times\cP_2(\mR^n),\mR^{n\times n})$, the coefficients $\bar{b}_{1,p}$ and $\bar{\s}_{1,p}$ satisfy $\bar{b}_{1,p}\in (C_b^{4,(1,3)}\cap C_b^{2,(2,2)})(\mR^n\times\cP_2(\mR^n),\mR^n)$ and $\bar{\s}_{1,p}\in (C_b^{4,(1,3)}\cap C_b^{2,(2,2)})(\mR^n\times\cP_2(\mR^n),\mR^{n\times n})$. Consequently, Eq.(\ref{perweaveeq}) has a unique solution $\{\hat{X}_t\}_{t\geq0}$. Let $\hat{X}^{s,\xi}_t$ be the unique solution to Eq.\eqref{perweaveeq} with initial data $\xi\in L^2(\Omega, \sF_s, \mP; \mR^n)$ at time $s$, i.e., for $t\geq s$,
\ce
\dif \hat{X}^{s,\xi}_t=\bar{b}_{1,p}(\hat{X}^{s,\xi}_t,\sL_{\hat{X}^{s,\xi}_t})\dif t+\bar{\s}_{1,p}(\hat{X}^{s,\xi}_t,\sL_{\hat{X}^{s,\xi}_t})\dif \bar{B}_t,\quad \hat{X}^{s,\xi}_s=\xi.
\de
Consider the decoupled equation
\ce
\dif \hat{X}^{s,x,\sL_{\xi}}_t=\bar{b}_{1,p}(\hat{X}^{s,x,\sL_{\xi}}_t,\sL_{\hat{X}^{s,\xi}_t})\dif t+\bar{\s}_{1,p}(\hat{X}^{s,x,\sL_{\xi}}_t,\sL_{\hat{X}^{s,\xi}_t})\dif \bar{B}_t,\quad \hat{X}^{s,x,\sL_{\xi}}_s=x.
\de
Thus $\hat{X}^{s,\xi}_t=\hat{X}^{s,x,\sL_{\xi}}_t\big|_{x=\xi}$. We can verify that for any $T>0$,
\be
&&\sup_{0\leq s\leq t\leq T}\left(\sum_{i=1}^4\mE\|\p_x^i\hat{X}^{s,x,\sL_{\xi}}_t\|^4+\sum_{i=0}^3\mE\|\p_{\tx}^i\p_\mu\hat{X}^{s,x,\sL_{\xi}}_t(\tx)\|^2\right)\leq C_T,\no\\
&&\sup_{0\leq s\leq t\leq T}\left(\sum_{i=1}^2\mE\|\p_x^i\p_\mu\hat{X}^{s,x,\sL_{\xi}}_t(\tx)\|^2+\sum_{i=1}^2\mE\|\p_x^i\p_{\tx}\p_\mu\hat{X}^{s,x,\sL_{\xi}}_t(\tx)\|^2 \right)\leq C_T,\no\\
&&\sup_{0\leq s\leq t\leq T}\sum_{i,j=0}^1\mE\|\p_{\cx}^i\p_{\tx}^j\p_\mu^2\hat{X}^{s,x,\sL_{\xi}}_t(\tx,\cx)\|^2\leq C_T.
\label{pxmubarxper}
\ee
For any $\varphi\in (C_b^{(1,3)}\cap C_b^{(2,2)})(\cP_2(\mR^n),\mR)$, we consider the following Cauchy problem on $[0,T]\times\cP_2(\mR^n)$ :
\be
\left\{\begin{array}{l}
\p_su^t(s,\sL_{\xi})+\mE[\bar{b}_{1,p}(\xi,\sL_{\xi})\cdot\p_\mu u^t(s,\sL_{\xi})(\xi)]\\
\qquad+\frac12\mE\left[Tr\big((\overline{\s_1\s_1^*})_p(\xi,\sL_{\xi})\cdot\p_{\tx}\p_\mu u^t(s,\sL_{\xi})(\xi)\big)\right]=0,\\
u^t(t,\sL_{\xi})=\varphi(\sL_{\xi}),\quad t\geq s.
\end{array}
\right.
\label{cauchyproper}
\ee
The unique solution of (\ref{cauchyproper}) is given by
\be
u^t(s,\sL_{\xi})=\varphi(\sL_{\hat{X}^{s,\xi}_t}),\quad t\geq s.\label{cauchypersol}
\ee
According to $\varphi\in (C_b^{(1,3)}\cap C_b^{(2,2)})(\cP_2(\mR^n),\mR)$ and (\ref{pxmubarxper}), for any $0\leq s\leq t\leq T$, $u^t(s,\cdot)\in(C_b^{(1,3)}\cap C_b^{(2,2)})(\cP_2(\mR^n),\mR)$
and
\be\label{pspmuutper}
&&\sup _{0\leq s\leq t\leq T}\max\Big\{\|\p_s\p_{\mu}u^t(s,\mu)(\tx)\|,\|\p_s\p_{\tx}\p_{\mu}u^t(s,\mu)(\tx)\|\Big\}\no\\
&\leq& C_T\left(1+|\tx|+\mu^{1/2}(|\cdot|^2)\right).
\ee

{\bf Proof of Theorem \ref{xbarweakper}.}
Applying the It\^{o} formula, we derive that
\ce
&&\varphi(\sL_{\check{X}_t^\eps})-\varphi(\sL_{\hat{X}_t})
=u^{t}(t,\sL_{\check{X}_t^\eps})-u^{t}(0,\sL_{\xi})\\
&=&\int_0^t\p_su^t(s,\sL_{\check{X}_s^\eps})\dif s+\mE\int_0^t\bar{b}_1(s/\eps,\check{X}_s^{\eps}, \sL_{\check{X}_s^{\eps}})\cdot\p_\mu u^t(s,\sL_{\check{X}_s^\eps})(\check{X}_s^\eps)\dif s\\
&&+\frac12\mE\int_0^tTr\big(\overline{\s_1\s_1^*}(s/\eps,\check{X}_s^{\eps}, \sL_{\check{X}_s^{\eps}})\cdot\p_{\tilde{x}}\p_\mu u^t(s,\sL_{\check{X}_s^\eps})(\check{X}_s^\eps)\big)\dif s\\
&=&\mE\int_0^t\left(\bar{b}_1(s/\eps,\check{X}_s^{\eps}, \sL_{\check{X}_s^{\eps}})-\bar{b}_{1,p}(\check{X}_s^{\eps}, \sL_{\check{X}_s^{\eps}})\right)\cdot\p_\mu u^t(s,\sL_{\check{X}_s^\eps})(\check{X}_s^\eps)\dif s\\
&&+\frac12\mE\int_0^tTr\Big(\left[\overline{\s_1\s_1^*}(s/\eps,\check{X}_s^{\eps}, \sL_{\check{X}_s^{\eps}})-(\overline{\s_1\s_1^*})_p(\check{X}_s^{\eps}, \sL_{\check{X}_s^{\eps}})\right]\\
&&\qquad\qquad\qquad\qquad\qquad\qquad\cdot\p_{\tilde{x}}\p_\mu u^t(s,\sL_{\check{X}_s^\eps})(\check{X}_s^\eps)\Big)\dif s.
\de
Note that
\ce
&&\sup_{0\leq t\leq T}|\varphi(\sL_{\check{X}_t^\eps})-\varphi(\sL_{\hat{X}_t})|\\
&\leq&\sup_{0\leq t\leq T}\bigg|\mE\int_0^t\Big[\left(\bar{b}_1(s/\eps,\check{X}_s^{\eps}, \sL_{\check{X}_s^{\eps}})-\bar{b}_{1,p}(\check{X}_s^{\eps}, \sL_{\check{X}_s^{\eps}})\right)\cdot\p_\mu u^t(s,\sL_{\check{X}_s^\eps})(\check{X}_s^\eps)\\
&&\qquad\qquad-\left(\bar{b}_1(s/\eps,\check{X}_{s(\D)}^{\eps}, \sL_{\check{X}_{s(\D)}^{\eps}})-\bar{b}_{1,p}(\check{X}_{s(\D)}^{\eps}, \sL_{\check{X}_{s(\D)}^{\eps}})\right)\cdot\p_\mu u^t({s(\D)},\sL_{\check{X}_{s(\D)}^\eps})(\check{X}_{s(\D)}^\eps)\Big]\dif s\bigg|\\
&&+\sup_{0\leq t\leq T}\bigg|\mE\int_0^t\left(\bar{b}_1(s/\eps,\check{X}_{s(\D)}^{\eps}, \sL_{\check{X}_{s(\D)}^{\eps}})-\bar{b}_{1,p}(\check{X}_{s(\D)}^{\eps}, \sL_{\check{X}_{s(\D)}^{\eps}})\right)\cdot\p_\mu u^t({s(\D)},\sL_{\check{X}_{s(\D)}^\eps})(\check{X}_{s(\D)}^\eps)\dif s\bigg|\\
&&+\frac12\sup_{0\leq t\leq T}\bigg|\mE\int_0^t\Big\{Tr\Big(\left[\overline{\s_1\s_1^*}(s/\eps,\check{X}_s^{\eps}, \sL_{\check{X}_s^{\eps}})-(\overline{\s_1\s_1^*})_p(\check{X}_s^{\eps}, \sL_{\check{X}_s^{\eps}})\right]\cdot\p_{\tilde{x}}\p_\mu u^t(s,\sL_{\check{X}_s^\eps})(\check{X}_s^\eps)\Big)\\
&&\qquad\qquad\qquad\qquad-Tr\Big(\big[\overline{\s_1\s_1^*}(s/\eps,\check{X}_{s(\D)}^{\eps}, \sL_{\check{X}_{s(\D)}^{\eps}})-(\overline{\s_1\s_1^*})_p(\check{X}_{s(\D)}^{\eps}, \sL_{\check{X}_{s(\D)}^{\eps}})\big]\\
&&\qquad\qquad\qquad\qquad\qquad\qquad\qquad\qquad\cdot\p_{\tilde{x}}\p_\mu u^t({s(\D)},\sL_{\check{X}_{s(\D)}^\eps})(\check{X}_{s(\D)}^\eps)\Big)\Big\}\dif s\bigg|\\
&&+\frac12\sup_{0\leq t\leq T}\bigg|\mE\int_0^tTr\Big(\big[\overline{\s_1\s_1^*}(s/\eps,\check{X}_{s(\D)}^{\eps}, \sL_{\check{X}_{s(\D)}^{\eps}})-(\overline{\s_1\s_1^*})_p(\check{X}_{s(\D)}^{\eps}, \sL_{\check{X}_{s(\D)}^{\eps}})\big]\\
&&\qquad\qquad\qquad\qquad\qquad\qquad\qquad\qquad\qquad\qquad\cdot\p_{\tilde{x}}\p_\mu u^t({s(\D)},\sL_{\check{X}_{s(\D)}^\eps})(\check{X}_{s(\D)}^\eps)\Big)\dif s\bigg|\\
&=:&\sum_{i=1}^4\sR_i(T).
\de
For $\sR_1(T)$, by (\ref{barb1lip}), (\ref{barb1grow}), (\ref{barb1pgrow}), (\ref{barb1plip}), (\ref{pspmuutper}) and $u^t(s,\cdot)\in(C_b^{(1,3)}\cap C_b^{(2,2)})(\cP_2(\mR^n),\mR)$, it follows that
\ce
&&\sR_1(T)\\
&\leq&\sup_{0\leq t\leq T}\bigg|\mE\int_0^t\Big[\bar{b}_1(s/\eps,\check{X}_s^{\eps}, \sL_{\check{X}_s^{\eps}})-\bar{b}_{1,p}(\check{X}_s^{\eps}, \sL_{\check{X}_s^{\eps}})\\
&&\qquad\qquad-\left(\bar{b}_1(s/\eps,\check{X}_{s(\D)}^{\eps}, \sL_{\check{X}_{s(\D)}^{\eps}})-\bar{b}_{1,p}(\check{X}_{s(\D)}^{\eps}, \sL_{\check{X}_{s(\D)}^{\eps}})\right)\Big]\cdot\p_\mu u^t(s,\sL_{\check{X}_s^\eps})(\check{X}_s^\eps)\dif s\bigg|\\
&&+\sup_{0\leq t\leq T}\bigg|\mE\int_0^t\left(\bar{b}_1(s/\eps,\check{X}_{s(\D)}^{\eps}, \sL_{\check{X}_{s(\D)}^{\eps}})-\bar{b}_{1,p}(\check{X}_{s(\D)}^{\eps}, \sL_{\check{X}_{s(\D)}^{\eps}})\right)\\
&&\qquad\qquad\qquad\qquad\cdot\left(\p_\mu u^t(s,\sL_{\check{X}_s^\eps})(\check{X}_s^\eps)-\p_\mu u^t({s(\D)},\sL_{\check{X}_{s(\D)}^\eps})(\check{X}_{s(\D)}^\eps)\right)\dif s\bigg|\\
&\leq&C_T\mE\int_0^T\left[1+|\check{X}_{s(\D)}^{\eps}|+(\mE|\check{X}_{s(\D)}^{\eps}|^2)^{1/2}\right]\cdot\left[|\check{X}_s^{\eps}-\check{X}_{s(\D)}^{\eps}|+(\mE|\check{X}_s^{\eps}-\check{X}_{s(\D)}^{\eps}|^2)^{1/2}\right]\dif s\\
&&+C_T\D\int_0^T(1+\mE|\check{X}_s^{\eps}|^2+\mE|\check{X}_{s(\D)}^{\eps}|^2)\dif s.
\de
Arguing as for $\sR_1(T)$, we obtain
\ce
\sR_3(T)\leq C_T\int_0^T(\mE|\check{X}_s^{\eps}-\check{X}_{s(\D)}^{\eps}|^2)^{1/2}\dif s+C_T\D\int_0^T(1+\mE|\check{X}_s^{\eps}|^2)\dif s.
\de
In fact, for $\sR_2(T)$,
\ce
&&\sR_2(T)\\
&\leq&\sup_{0\leq t\leq T}\bigg|\sum_{k=0}^{[t/\D]-1}\mE\int_{k\D}^{(k+1)\D}\left(\bar{b}_1(s/\eps,\check{X}_{k\D}^{\eps}, \sL_{\check{X}_{k\D}^{\eps}})-\bar{b}_{1,p}(\check{X}_{k\D}^{\eps}, \sL_{\check{X}_{k\D}^{\eps}})\right)\cdot\p_\mu u^t(k\D,\sL_{\check{X}_{k\D}^\eps})(\check{X}_{k\D}^\eps)\dif s\bigg|\\
&&+\sup_{0\leq t\leq T}\bigg|\mE\int_{[t/\D]\D}^t\left(\bar{b}_1(s/\eps,\check{X}_{[t/\D]\D}^{\eps}, \sL_{\check{X}_{[t/\D]\D}^{\eps}})-\bar{b}_{1,p}(\check{X}_{[t/\D]\D}^{\eps}, \sL_{\check{X}_{[t/\D]\D}^{\eps}})\right)\\
&&\qquad\qquad\qquad\qquad\qquad\qquad\qquad\qquad\cdot\p_\mu u^t([t/\D]\D,\sL_{\check{X}_{[t/\D]\D}^\eps})(\check{X}_{[t/\D]\D}^\eps)\dif s\bigg|\\
&=:&\sR_{21}(T)+\sR_{22}(T).
\de
For $\sR_{21}(T)$, since $\bar{b}_1(\cdot,x,\mu)$ is $\tau$-periodic, combining $u^t(s,\cdot)\in(C_b^{(1,3)}\cap C_b^{(2,2)})(\cP_2(\mR^n),\mR)$ and Lemma \ref{fperlem}, one has
\ce
&&\sR_{21}(T)\\
&\leq&\sum_{k=0}^{[T/\D]-1}\mE\bigg[\bigg|\int_{k\D}^{(k+1)\D}\bar{b}_1(s/\eps,\check{X}_{k\D}^{\eps}, \sL_{\check{X}_{k\D}^{\eps}})\dif s-\D\bar{b}_{1,p}(\check{X}_{k\D}^{\eps}, \sL_{\check{X}_{k\D}^{\eps}})\bigg|\\
&&\qquad\qquad\qquad\qquad\qquad\cdot|\p_\mu u^t(k\D,\sL_{\check{X}_{k\D}^\eps})(\check{X}_{k\D}^\eps)|\bigg]\\
&\leq&C_T\D\sum_{k=0}^{[T/\D]-1}\mE\bigg|\frac{\eps}{\D}\int_{\frac{k\D}{\eps}}^{\frac{(k+1)\D}{\eps}}\bar{b}_1(u,\check{X}_{k\D}^{\eps}, \sL_{\check{X}_{k\D}^{\eps}})\dif u-\frac1\tau\int_0^\tau\bar{b}_1(u,\check{X}_{k\D}^{\eps}, \sL_{\check{X}_{k\D}^{\eps}})\dif u\bigg|\\
&\leq&C_T\eps\sum_{k=0}^{[T/\D]-1}\mE(\sup_{0\leq u\leq \tau}|\bar{b}_1(u,\check{X}_{k\D}^{\eps}, \sL_{\check{X}_{k\D}^{\eps}})|)\\
&\leq&C_T\frac{\eps}{\D}\big(1+\mE(\sup_{0\leq t\leq T}|\check{X}_t^{\eps}|^2)\big).
\de
For $\sR_{22}(T)$, invoking (\ref{barb1grow}), (\ref{barb1pgrow}) and $u^t(s,\cdot)\in(C_b^{(1,3)}\cap C_b^{(2,2)})(\cP_2(\mR^n),\mR)$, we conclude that
\ce
\sR_{22}(T)&\leq&\sup_{0\leq t\leq T}\mE\int_{[t/\D]\D}^t\Big(|\bar{b}_1(s/\eps,\check{X}_{[t/\D]\D}^{\eps}, \sL_{\check{X}_{[t/\D]\D}^{\eps}})|+|\bar{b}_{1,p}(\check{X}_{[t/\D]\D}^{\eps}, \sL_{\check{X}_{[t/\D]\D}^{\eps}})|\Big)\\
&&\qquad\qquad\qquad\qquad\qquad\qquad\cdot|\p_\mu u^t([t/\D]\D,\sL_{\check{X}_{[t/\D]\D}^\eps})(\check{X}_{[t/\D]\D}^\eps)|\dif s\\
&\leq&C_T\D\big(1+\mE(\sup_{0\leq t\leq T}|\check{X}_t^{\eps}|^2)\big).
\de
Similarly, it holds that
\ce
\sR_4(T)\leq C_T(\frac{\eps}{\D}+\D).
\de
Combining these estimates with (\ref{barxetsweak}) and (\ref{barxeweak}), and taking $\D=\eps^{2/3}$, we arrive at
\ce
\sup_{0\leq t\leq T}|\varphi(\sL_{\check{X}_t^\eps})-\varphi(\sL_{\hat{X}_t})|\leq C_T(1+\mE|\xi|^2)\eps^{1/3}.
\de
Together with (\ref{lxelbarxe}), this further implies that
\ce
\sup_{0\leq t\leq T}|\varphi(\sL_{X_t^{\eps}})-\varphi(\sL_{\hat{X}_t})|\leq C_T(1+\mE|\xi|^4+\mE|\varrho|^4)\eps\left[\sup_{0\leq t\leq T}\Upsilon_{\g}(t/\eps)+\eps^{-2/3}\right].
\de
The proof is complete.

\section{Examples}\label{example}

In this section, we present two examples to illustrate our results.

\bx\label{ex1}
Let $\a:[0,+\infty)\rightarrow(0,+\infty)$ satisfy (\ref{alphaposi}). Consider the slow-fast system
\be\left\{\begin{array}{l}
dX_t^{\eps} = [Y_t^{\eps}+\int_{\mR}z\sL_{X_t^\e}(\dif z)]dt+dB_t,\quad X_0^{\eps} = \xi, \\
dY_t^{\eps} = -\eps^{-1}\alpha(t/\eps)\left[2Y_t^{\eps}-\int_{\mR}z\sL_{X_t^\e}(\dif z)\right]dt+\left[\eps^{-1}\alpha(t/\eps)\right]^{1/2} dW_t,\quad Y_0^{\eps} = \varrho, \\
\end{array}\right.
\label{exorieq}
\ee
 where $(B_t)$ and $(W_t)$ are independent one-dimensional standard Brownian motions defined on the complete filtered probability space $(\Omega,\sF,\{\sF_t\}_{t \in[0, T]}, \mathbb{P})$. $\varrho$ and $\xi$ are $\sF_0$-measurable Gaussian random variables.

Note that
\ce
&&b_1(x,\mu,y)=y+\int_{\mR}z\mu(\dif z),\quad \s_1(x,\mu,y)=1,\\
&&b_2(t,x,\mu,y)=-\a(t)(2y-\int_{\mR}z\mu(\dif z)),\quad \s_2(t,x,\mu,y)=\a^{1/2}(t).
\de
Fix $\mu$ and denote
\ce
m:=\int_{\mR}z\mu(\dif z).
\de
The corresponding frozen fast equation is
\be
\dif Y_t=-\alpha(t)(2Y_t-m)\dif t+\alpha^{1/2}(t)\dif W_t,\qquad Y_s=y.\label{yex1}
\ee
The solution to Eq.(\ref{yex1}) admits the explicit representation
\ce
Y_t^{s,\mu,y}=\frac{m}2+e^{-2\int_s^t\a(r)\dif r}(y-\frac{m}2)+\int_s^te^{-2\int_u^t\alpha(r)\dif r}\a^{1/2}(u)\dif W_u.
\de
Consequently, $Y_t^{s,\mu,y}$ is Gaussian with distribution
\ce
Y_t^{s,\mu,y}\sim N\left(\frac{m}2+e^{-2\int_s^t\a(r)\dif r}(y-\frac{m}2),\frac14\bigl(1-e^{-4\int_s^t\alpha(r)\dif r}\bigr)\right).
\de
Letting $s\to-\infty$, we obtain
\ce
\nu_t^\mu:=N\left(\frac{m}2,\frac14\right),\qquad t\in\mathbb{R},
\de
is an evolution system of measures for Eq.(\ref{yex1}). The averaging drift is then given by
\ce
\bar b_1(t,x,\mu)=\int_{\mathbb{R}}\left(y+\int_{\mathbb{R}}z\,\mu(\dif z)\right)\nu_t^\mu(\dif y)=\frac32m.
\de
Therefore, the averaging equation takes the form
\ce
\dif\bar X^\e_t=\frac32\int_{\mR}z\sL_{\bar X^\e_t}(\dif z)\dif t+\dif B_t,\qquad \bar X^\e_0=\xi.
\de
A straightforward verification shows that $b_1, \s_1, b_2, \s_2$ satisfy $(\mathbf{H}^1_{b_{1}, \s_{1}})$-$(\mathbf{H}^3_{b_{1}, \s_{1}})$ and $(\mathbf{H}^3_{b_{2}, \s_{2}})$-$(\mathbf{H}^4_{b_{2}, \s_{2}})$. Hence, by Theorem \ref{xbarx2} and \ref{xbarweak} strong and weak averaging principles hold for the system (\ref{exorieq}).
\ex

\begin{remark}
When $\alpha(t)\equiv1$, the frozen fast equation admits an invariant measure and Example \ref{ex1} falls into the framework considered in \cite{rsx}. In this case, Theorem \ref{xbarx2} and \ref{xbarweak} reduce to the corresponding strong and weak averaging principles. Moreover, the strong averaging rate is of order $\varepsilon^{1/2}$, which coincides with the classical rate established in \cite{rsx}, while the weak averaging rate is of order $\varepsilon$.
\end{remark}

\bx\label{ex2}
Let $\phi:[0,+\infty)\rightarrow(0,+\infty)$ be bounded. Consider the slow-fast system
\ce
\left\{\begin{array}{l}
\dif X_t^\e=\Big[\sin(X_t^\e)+Y_t^{\eps}+\int_{\mR}\cos z\sL_{X_t^\e}(\dif z)\Big]\dif t+dB_t,\quad X_0^{\eps} = \xi, \\
dY_t^{\eps} = -\eps^{-1}\left[2Y_t^{\eps}-\phi(t/\eps)\int_{\mR}z\sL_{X_t^\e}(\dif z)\right]dt+\eps^{-1/2}X_t^\e dW_t,\quad Y_0^{\eps} = \varrho, \\
\end{array}\right.
\de
where $\left(B_t\right),\left(W_t\right)$ are $1$-dimensional standard Brownian motions, respectively, defined on the complete filtered probability space $(\Omega,\sF,\{\sF_t\}_{t \in[0, T]}, \mathbb{P})$ and are mutually independent. $\varrho, \xi$ are $\sF_0$-measurable Gaussian random variables.

The coefficients are given by
\ce
&&b_1(x,\mu,y)=\sin x+y+\int_{\mR}\cos z\mu(\dif z),\quad \s_1(x,\mu,y)=1,\\
&&b_2(t,x,\mu,y)=-\left(2y-\phi(t)\int_{\mR}z\mu(\dif z)\right),\quad \s_2(t,x,\mu,y)=x.
\de
Fix $\mu$ and set
\ce
m=\int_{\mR}z\mu(\dif z).
\de
For fixed $x$ and $\mu$, the associated frozen fast equation is
\be
\dif Y_t=-(2Y_t-\phi(t)m)\dif t+x\dif W_t,\qquad Y_s=y.\label{yex2}
\ee
The solution admits the explicit representation
\ce
Y_t^{s,x,\mu,y}=e^{-2(t-s)}y+m\int_s^t e^{-2(t-u)}\phi(u)\dif u+x\int_s^t e^{-2(t-u)}\dif W_u .
\de
Thus,
\ce
Y_t^{s,x,\mu,y}\sim N\left(e^{-2(t-s)}y+m\int_s^t e^{-2(t-u)}\phi(u)\dif u,\frac{x^2}4\bigl(1-e^{-4(t-s)}\bigr)\right).
\de
Letting $s\to-\infty$, we obtain
\ce
\nu_t^{x,\mu}:=N\left(m\int_{-\infty}^t e^{-2(t-u)}\phi(u)\dif u,\frac{x^2}4\right),\qquad t\in\mathbb{R}
\de
is an evolution system of measures for Eq.(\ref{yex2}). The averaging drift is given by
\ce
\bar b_1(t,x,\mu)=\int_{\mathbb{R}}\left(\sin x+y+\int_{\mR}\cos z\mu(\dif z)\right)\nu_t^{x,\mu}(\dif y)=\sin x+\int_{\mR}\cos z\mu(\dif z)+m\Phi(t),
\de
where $\Phi(t)=\int_{-\infty}^t e^{-2(t-u)}\phi(u)\dif u$. Therefore, the averaging equation takes the form
\ce
\dif\bar X_t^\eps=\left[\sin(\bar X_t^\eps)+\int_{\mR}\cos z\sL_{\bar X_t^\eps}(\dif z)+\int_{\mR}z\sL_{\bar X_t^\eps}(\dif z)\Phi(t/\eps)\right]\dif t+\dif B_t,\qquad \bar X_0^\eps=\xi.
\de
 It is straightforward to verify that $(\mathbf{H}^1_{b_{1}, \s_{1}})$-$(\mathbf{H}^3_{b_{1}, \s_{1}})$ and $(\mathbf{H}^3_{b_{2}, \s_{2}})$-$(\mathbf{H}^4_{b_{2}, \s_{2}})$ hold with $\a(t)\equiv1$ and $\Upsilon(t)\equiv1$. Consequently, the conclusions of Theorem \ref{xbarx2} and \ref{xbarweak} follow.

\textbf{The convergent case:} Assume that
\ce
\lim_{t\to +\infty} \phi(t) = 0.
\de
Then $(\mathbf{H}^5_{b_{2}, \s_{2}})$ holds with $\bar{b}_2(x,\mu,y)=-2y$ and $\bar{\s}_2(x,\mu,y)=x$. Note that
\ce
\Phi(t)=\int_{-\infty}^t e^{-2(t-u)}\phi(u)\dif u=\int_0^{+\infty}e^{-2s}\phi(t-s)\dif s.
\de
By the dominated convergence theorem, it follows that
\ce
\lim_{t\to +\infty} \Phi(t) = 0.
\de
The averaging equation is given by
\ce
\dif\bar X_t=\left[\sin(\bar X_t)+\int_{\mR}\cos z\sL_{\bar X_t}(\dif z)\right]\dif t+\dif B_t,\qquad \bar X_0=\xi.
\de
The conclusions of Theorems \ref{xbarx2con} and \ref{xbarweakcon} follow.

\textbf{The periodic case:}
Assume $\phi$ is $\tau$-periodic, that is, $\phi(t+\tau)=\phi(t)$ for all $t\in\mathbb{R}$. In this situation, it can be verified that $(\mathbf{H}^6_{b_{2}, \s_{2}})$ holds. The associated averaging equation is given by
\ce
\dif\tilde X_t=\left[\sin(\tilde X_t)+\int_{\mR}\cos z\sL_{\tilde X_t}(\dif z)+\frac{1}{\tau}\int_0^{\tau}\Phi(r)\dif r\int_{\mR}z\sL_{\tilde X_t}(\dif z)\right]\dif t+\dif B_t,\qquad \tilde X_0=\xi.
\de
Therefore, the results in Theorem \ref{xbarstrper} and \ref{xbarweakper} hold.
\ex

\begin{remark}
We mention that in \cite{slg} Shi, Liu and Gao established strong and weak averaging principles for nonautonomous multiscale McKean-Vlasov SDEs with almost periodic coefficients. However, they didn't give explicit convergence rates. Here we provide strong and weak averaging rates of order $\varepsilon^{1/3}$.
\end{remark}


\begin{thebibliography}{99}
\bibitem{br}
V. Barbu and M. R\"{o}ckner: From nonlinear Fokker-Planck equations to solutions of distribution dependent SDE, {\it Ann. Probab.}, 48(2020)1902-1920.

\bibitem{br1}
V. Barbu and M. R\"{o}ckner: Uniqueness for nonlinear Fokker-Planck equations and weak uniqueness for McKean-Vlasov SDEs, {\it Stoch PDE: Anal Comp.}, 9(2021)702-713.

\bibitem{bm}
M. Branicki and A. J. Majda: Quantifying uncertainty for predictions with model error in non-Gaussian systems with intermittency, {\it Nonlinearity}, 25(2012)2543-2578.

\bibitem{blpr}
R. Buckdahn, J. Li, S. Peng and C. Rainer: Mean-field stochastic differential equations and associated PDEs, {\it Ann. Probab.}, 45(2017)824-878.

\bibitem{card}
P. Cardaliaguet: Notes on mean field games, P.-L. Lions Lectures at College de France, https://www.researchgate.net/publication/228702832. Notes on Mean Field Games.

\bibitem{carm}
R. Carmona and F. Delarue: Probability Theory of Mean Field Games with Applications I, Vol. 83. Springer, Cham, 2018.

\bibitem{cl}
S. Cerrai and A. Lunardi: Averaging principle for nonautonomous slow-fast systems of stochastic reaction-diffusion equations: the almost periodic case, {\it SIAM J. Math. Anal.}, 49(2017)2843-2884.

\bibitem{csg}
M. D. Chekroun, E. Simonnet and M. Ghil: Stochastic climate dynamics: Random attractors and time dependent invariant measures, {\it Phys. D}, 240 (2011)1685-1700.

\bibitem{csx}
M. Cheng, X. Sun and Y. Xie: Strong averaging principle for nonautonomous multi-scale SPDEs with fully local monotone and almost periodic coefficients, {\it J. Funct. Anal.}, 290(2026)111411.

\bibitem{dr}
G. Da Prato and M. R\"{o}ckner: A note on evolution systems of measures for time-dependent stochastic differential equations, in: Seminar on Stochastic Analysis, Random Fields and Applications V, 2008, 115-122.

\bibitem{dn}
P. Del Moral and A. Niclas: A Taylor expansion of the square root matrix function, {\it J. Math. Anal. Appl.}, 465(2018)259-266.

\bibitem{dq1}
X. Ding and H. Qiao: Euler-Maruyama approximations for stochastic McKean-Vlasov equations with non-Lipschitz coefficients, {\it J. Theoret. Probab.}, 34(2021)1408-1425.

\bibitem{dq2}
X. Ding and H. Qiao: Stability for stochastic McKean-Vlasov equations with non-Lipschitz coefficients, {\it SIAM J. Control Optim.}, 59(2021)887-905.

\bibitem{gw}
M. Galtier and G. Wainrib: Multi-scale analysis of slow-fast neuronal learning models with noise, {\it J. Math. Neurosci.}, 2(2012)1-64.

\bibitem{hll}
W. Hong, S. Li and W. Liu: Strong convergence rates in averaging principle for slow-fast McKean-Vlasov SPDEs, {\it J. Differential Equations}, 316(2022)94-135.

\bibitem{hwx}
Y. Hou, L. Wang and L. Xie: Diffusion approximation for non-autonomous multi-scale stochastic systems, {\it Commun. Pure Appl. Anal.}, 2026, https://doi.org/10.3934/cpaa.2026048.

\bibitem{hw}
X. Huang and F.-Y. Wang: Distribution dependent SDEs with singular coefficients, {\it Stochastic Process. Appl.}, 129(2019)4747-4770.

\bibitem{rk}
R. Khasminskii: On the principle of averaging the It\^{o} stochastic differential equations, {\it Kibernetika}, 4(1968)260-279.

\bibitem{lswx}
Y. Li, X. Sun, Z. Wang and Y. Xie: Strong averaging principle for nonautonomous slow-fast SPDEs driven by $\a$-stable processes, https://arxiv.org/abs/2507.07538v1.

\bibitem{lx}
Y. Li and L. Xie: Functional law of large numbers and central limit theorem for slow-fast McKean-Vlasov equations, {\it Discrete Contin. Dyn. Syst. Ser. S}, 16(2023)846-877.

\bibitem{lrsx}
W. Liu, M. R\"{o}ckner, X. Sun and Y. Xie: Averaging principle for slow-fast stochastic differential equations with time dependent locally Lipschitz coefficients, {\it J. Differential Equations}, 268(2020)2910-2948.

\bibitem{hpm}
H. P. McKean: A class of Markov processes associated with nonlinear parabolic equations, Proc. Natl. Acad. Sci. U.S.A., 56(1966)1907-1911.


\bibitem{pix}
B. Pei, Y. Inahama and Y. Xu: Averaging principle for fast-slow system driven by mixed fractional Brownian rough path, {\it J. Differential Equations}, 301(2021)202-235.

\bibitem{q}
H. Qiao: Large deviations of multiscale multivalued McKean-Vlasov stochastic systems, {\it Stoch. Anal. Appl.}, 2026, https://doi.org/10.1080/07362994.2026.2670720.

\bibitem{qw}
H. Qiao and W. Wei: Strong approximation of nonlinear filtering for multiscale McKean-Vlasov stochastic systems, to appear in {\it Acta Math. Sci. Ser. B}, 2026.

\bibitem{qw2}
H. Qiao and W. Wei: Weak approximation of nonlinear filtering for multiscale McKean-Vlasov stochastic systems, {\it Nonlinear Differ. Equ. Appl.}, 33(2026)67.

\bibitem{rw}
P. Ren and F.-Y. Wang: Space-Distribution PDEs for path independent additive functionals of McKean-Vlasov SDEs, {\it Infin. Dimens. Anal. Quantum Probab. Relat. Top.}, 23(2020)2050018.

\bibitem{rsx}
 M. R\"{o}ckner, X. Sun and Y. Xie: Strong convergence order for slow-fast McKean-Vlasov stochastic differential equations, {\it Ann. Inst. H. Poincar\'{e} Probab. Statist.}, 57(2021)547-576.

\bibitem{rz}
M. R\"{o}ckner and X. Zhang: Well-posedness of distribution dependent SDEs with singular drifts, {\it Bernoulli}, 27(2021)1131-1158.

\bibitem{sxw}
G. Shen, J. Xiang and J.-L. Wu: Averaging principle for distribution dependent stochastic differential equations driven by fractional Brownian motion and standard Brownian motion, {\it J. Differential Equations}, 321(2022)381-414.

\bibitem{syw}
G. Shen, J. Yin and J.-L. Wu: Stochastic averaging principle for two-time-scale SDEs with distribution dependent coefficients driven by fractional Brownian motion, {\it Commun. Math. Stat.}, 13(2025)1445-1479.

\bibitem{slg}
Y. Shi, H. Liu and H. Gao: Weak and strong averaging principle for non-autonomous slow-fast McKean-Vlasov SDEs with almost periodic coefficients, {\it Math. Methods Appl. Sci.}, 48(2025)15111-15131.

\bibitem{swx0}
X. Sun, J. Wang and Y. Xie: Averaging principles for time-inhomogeneous multi-scale SDEs with partially dissipative coefficients, https://arxiv.org/abs/2506.18558v1.

\bibitem{swx}
X. Sun, J. Wang and Y. Xie: Averaging principles for time-inhomogeneous multi-scale SDEs via nonautonomous Poisson equations, https://arxiv.org/abs/2412.09850.

\bibitem{Uda}
K. Uda: Averaging principle for stochastic differential equations in the random periodic regime, {\it Stochastic Process. Appl.}, 139(2021)1-36.

\bibitem{gW}
G. Wainrib: Double averaging principle for periodically forced slow-fast stochastic systems, {\it Electron. Commun. Probab.}, 18(2013)1-12.

\bibitem{Wang}
F.-Y. Wang: Distribution dependent SDEs for Landau type equations, {\it Stochastic Process. Appl.}, 128(2018)595-621.

\bibitem{xq1}
J. Xiang and H. Qiao: Asymptotic behaviors of multiscale McKean-Vlasov stochastic systems, {\it Discrete Contin. Dyn. Syst. Ser. S}, 21(2026)63-87.

\bibitem{xq2}
J. Xiang and H. Qiao: Averaging principles and central limit theorems for multiscale McKean-Vlasov stochastic systems, https://arxiv.org/abs/2501.11853.

\bibitem{xw}
Y. Xu and R. Wang: Averaging principles for non-autonomous two-time-scale stochastic reaction-diffusion equations with jump, {\it Complexity}, 2020(2020)9864352.

\end{thebibliography}
\end{document}